\documentclass[a4paper,12pt]{article}

\begin{document}

\title{Quasi-conformal functions of quaternion and octonion variables,
their integral transformations.}
\author{S.V. Ludkovsky.}

\date{12 March 2007}
\maketitle

\begin{abstract}
The article is devoted to holomorphic and meromorphic functions of
quaternion and octonion variables. New classes of quasi-conformal
and quasi-meromorphic mappings are defined and investigated.
Properties of such functions such as their residues and argument
principle are studied. It is proved, that the family of all
quasi-conformal diffeomorphisms of a domain form a topological group
relative to composition of mappings. Cases when it is a
finite-dimensional Lie group over $\bf R$ are studied. Relations
between quasi-conformal functions and integral transformations of
functions over quaternions and octonions are established. For this,
in particular, noncommutative analogs of the Laplace and Mellin
transformations are studied and used. Examples of such functions are
given. Applications to problems of complex analysis are
demonstrated.
\end{abstract}

\section{Introduction}
\par Complex analysis is one of the corner-stones of mathematics.
On the other hand, natural generalizations of complex numbers
obtained by subsequent doubling procedures with the help of
generators were introduced in the second half of the 19-th century.
The most important among them are quaternions invented by W.R.
Hamilton and their generalizations such as octonions and
Cayley-Dickson numbers are known \cite{baez,kansol}. The problem of
developing analysis over quaternions ${\bf H}={\cal A}_2$ and
octonions ${\bf O}={\cal A}_3$ was posed by Hamilton and Yang and
Mills for the needs of celestial mechanics and quantum field theory
\cite{hamilt,guetze,rothe}. \par Quaternions and octonions were used
in quantum mechanics and quantum field theory and even by J.C.
Maxwell, but mainly algebraically, because theory of functions of
quaternion and octonion variables was little developed
\cite{emch,guetze,lawmich}. Noncommutative analysis is being
developed in recent years for the needs of mathematics and
theoretical physics \cite{berez,brdeso,connes,dewitt,khren,oystaey},
but it remains very little promoted in comparison with classical
analysis, especially its non super-commutative part. Derivations of
abstract algebras are widely used and a work with functions on
algebras is frequently related with their representations by words
and phrases \cite{bourbal,razmus}.

\par In preceding works of the author
super-differentiable (in another words holomorphic) functions of
Cayley-Dickson variables were investigated
\cite{luoyst,lusmfcdv,luhcnfcdv,lufscdvm} such that they generalize
the theory of complex holomorphic functions. In the particular case
of complex functions the notion of super-differentiability reduces
to the usual complex differentiability. In these publications
super-differentiability was defined as derivation of an algebra and
taking into account specific features of the Cayley-Dickson algebra.
In view of the noncommutativity of the Cayley-Dickson algebra ${\cal
A}_r$ with $2\le r$, the theory of functions over them is not only
the usual theory of functions, but it also bears the algebraic
structure and certainly is related with representations of functions
with the help of words and phrases over Cayley-Dickson variables.
\par It is necessary to note, that there are natural embeddings $\theta
^r_k$ of ${\cal A}_r$ into ${\cal A}_k$ for each $1\le r<k\in \bf N$
associated with the subsequent doubling procedure, but besides them
there are others algebraic embeddings. The algebra ${\cal A}_{\infty
}$ obtained by completion of the strict inductive limit $str-ind \{
{\cal A}_r, \theta ^r_k, {\bf N} \} $ relative to the $l_2$ norm has
no any internal anti-automorphism $z\mapsto z^*$, since it is
external, where $zz^*=|z|^2$. Therefore, it is natural to consider
holomorphic functions of ${\cal A}_r$ variables with $2\le r<\infty
$ as restrictions of functions of ${\cal A}_{\infty }$ variables on
the corresponding domains. Though the hypercomplex Cayley-Dickson
algebra ${\cal A}_{\infty }$ is noncommutative and non-associative,
but with respect to the absence of the internal anti-automorphism
$z\mapsto z^*$ it resembles by such property the complex field $\bf
C$. Then operator theory over Cayley-Dickson algebras on the base of
this function theory was studied in \cite{lusmalop,ludanavf}.
Super-differentiable functions are locally analytic by their
Cayley-Dickson variables, but series for them are more complicated
in comparison with the complex case due to noncommutativity for
$r\ge 2$ or non-associativity for $r\ge 3$ of ${\cal A}_r$. Then the
noncommutative analog of the Laplace transformation was studied in
\cite{lusmtslt,luhcnlt}. In particular, pseudo-conformal mappings
over quaternions and octonions having properties closer to that of
complex holomorphic functions were defined and studied in
\cite{ludagpd,lusmnfgpcd}. Pseudo-conformal mappings over the
quaternion skew field $\bf H$ or over the octonion algebra $\bf O$
are analogous to complex conformal functions, but in the
noncommutative setting, so generally pseudo-conformal functions may
be non-isometric (see Definition 2.1).
\par This work continues these investigations using preceding results.
Professor Fred Van Oystaeyen has formulated in 2002 the problem of
developing analysis over quaternions and octonions for the needs of
mathematics and theoretical physics, particularly, of noncommutative
geometry and their potential applications, as well as for problems
of complex analysis such as the Riemann's hypothesis \cite{luoyst}.
Natural extensions of complex holomorphic functions are introduced
over quaternions and octonions such that a notion of quasi-conformal
mappings is defined. They form a different class than that of
pseudo-conformal mappings. Quasi-conformal mappings on domains $U$
in ${\cal A}_b$ are formed from pseudo-conformal functions on
domains $W$ in the Cayley-Dickson subalgebra ${\cal A}_r$, $1\le
r<b\le 3$, with the help of operators to which rotations of the real
shadow ${\bf R}^{2^b}$ correspond.
\par In the second section new classes of quasi-conformal and
quasi-meromorphic mappings are defined and investigated. Properties
of such functions such as their residues and argument principle are
studied. It is proved, that the family of all quasi-conformal
diffeomorphisms of a domain form a topological group relative to
composition of mappings. Cases when it is a finite-dimensional Lie
group over $\bf R$ are studied. \par In the third section relations
between quasi-conformal functions and integral transformations of
functions over quaternions and octonions are established. For this,
in particular, noncommutative analogs of the Laplace and Mellin
transformations are studied and used. Examples of such functions are
given. An effectiveness of analysis over quaternions and octonions
is demonstrated for problems of complex analysis.
\par Many results of this paper are obtained for the first time.

\section{Quasi-conformal functions}
\par {\bf 1. Definitions and Notation.}
Let ${\cal A}_r$ denotes the Cayley-Dickson
algebra of dimension $2^r$ over $\bf R$, where in particular
${\bf C}={\cal A}_1$ is the field of complex numbers,
${\bf H}={\cal A}_2$ is the skew field of quaternions, ${\bf O}={\cal
A}_3$ is the algebra of octonions. Suppose that $U$ is an open subset
in ${\cal A}_r$, $2\le r\le 3$. A function $f$ on $U$ we call
pseudo-conformal at a point
$\xi $ in $U$, if $f$ is ${\cal A}_r$ holomorphic (super-differentiable)
in a neighborhood of $\xi $ and satisfies Conditions
$(P1-3)$: \par $(P1)$ $\partial f(z)/\partial {\tilde z}=0$ for $z=\xi
$; \par $(P2)$ $Re \{ [(\partial f(z)/
\partial z).h_1] [(\partial f(z))/ \partial z).h_2]^{*} \} |h_1|
|h_2|$ \par $ = |(\partial f(z)/\partial z).h_1| |(\partial f(z)/\partial
z).h_2| Re (h_1{\tilde h}_2)$ for $z=\xi $ for each $h_1, h_2\in
{\cal A}_r$,
\par $(P3)$ $(\partial f(z)/\partial z)|_{z=\xi }.h\ne 0$ for each
$h\ne 0$ in ${\cal A}_r$, where ${\tilde z}=z^*$ denotes the adjoint
number of $z\in {\cal A}_r$ such that $z {\tilde z} = |z|^2$, $Re
(z) := (z+{\tilde z})/2$; for $f$ it is used either the shortest
phrase compatible with these conditions or in the underlying real
space (shadow)  $\bf R^4$ or $\bf R^8$ non-proper rotations $[f']$
associated with $f'$ are excluded. That is, $[f']\in SO(2^b,{\bf
R})$, where $b=2$ or $b=3$, $SO(n,{\bf R})$ denotes the special
orthogonal group of ${\bf R}^n$, $[f']$ denotes the operator in the
real shadow corresponding to the super-derivative $f'$ over ${\cal
A}_r$. For short $f(z,{\tilde z})$ is written as $f(z)$ due to the
bijectivity between $z\in {\cal A}_r$ and $\tilde z$.
\par If $f$ is pseudo-conformal at
each point $\xi \in U$, then it is called pseudo-conformal in the
domain $U$. \par For mappings of complex numbers, $r=1$, each holomorphic
function satisfying Condition $(P3)$ fits to this definition, so we can
include this case also.
\par We say that a function $\phi $ at $\zeta $ or on $V$ is
$p$-pseudo-conformal, if $\phi (z)=f(z^p)$ and $f$ is
pseudo-conformal at $\xi $ or on $U$, where $\zeta ^p=\xi $ or $U =
\{ z^p: z\in V \} $, $p\in \bf N$. A function of several variables
$\mbox{ }_1z,...,\mbox{ }_nz$ is called pseudo-conformal or
$(p_1,...,p_n)$-pseudo-conformal, if it is pseudo-conformal or
$p_j$-pseudo-conformal by $\mbox{ }_jz$ for each $j=1,...,n$.
\par Let now $U$ be open in ${\cal A}_b$ and $W=U\cap {\cal A}_r\ne
\emptyset $ be open in ${\cal A}_r$ and non-void, where $1\le r<b\le
3$. Consider the natural embedding of ${\cal A}_r$ into ${\cal A}_b$
associated with the standard doubling procedure.
\par Suppose that $f$ is a holomorphic function on $U$
with values in ${\cal A}_b$ satisfying the following conditions:
\par $(Q1)$ the function $g(z) := f(y_0+z)$ has the
$p$-pseudo-conformal restriction $g|_{W-y_0}$ on
$W-y_0 := \{ z: z=x-y_0, x\in W \} $ for some marked point $y_0\in W$
and $g(W-y_0)\subset {\cal A}_r$,
\par $(Q2)$ there exists a family of automorphisms
${\hat R}_{z,x}={\hat R}^f_{z,x}: {\cal A}_b\to {\cal A}_b$ for each
$z \in U-y_0$ and $x\in W-y_0$ with $Re (z)=Re (x)$ such that to
each ${\hat R}_{z,x}$ a proper rotation $T=[{\hat R}_{z,x}]\in
SO(2^b,{\bf R})$ of the real shadow ${\bf R}^{2^b}$ corresponds such
that for each $z\in U-y_0$ there exists $x\in W-y_0$ for which
$z={\hat R}_{z,x}x$, where $SO(n,{\bf R})$ denotes the special
orthogonal group of the Euclidean space ${\bf R}^n$,
\par $(Q3)$ ${\hat R}_{z,x}|_{\bf R}=id|_{\bf R}$ for each $z\in U-y_0$
and every $x\in W-y_0$, that is, $T=[{\hat R}_{z,x}]\in
SO_{\bf R}(2^b,{\bf R})$, where $SO_{\bf R}(n,{\bf R}) :=
\{ T: T\in SO(n,{\bf R}); T|_{\bf R}=I \} $,
\par $(Q4)$ ${\hat R}_{z,x}=id$ for each $z\in W-y_0$ and every
$x\in W-y_0$,
\par $(Q5)$ ${\hat R}_{z,x}$ depends ${\cal A}_b$ holomorphically on
$z \in U-y_0$ and ${\cal A}_r$ holomorphically on $x\in W-y_0$ in a
suitable $(z,x)$-representation,
\par $(Q6)$ $g(z)={\hat R}_{z,x}g(x)$ for each $x \in W-y_0$
and every $z\in U-y_0$ such that $Re (z)=Re (x)$ and $z={\hat
R}_{z,x}x$,
\par $(Q7)$ $g'({\hat R}_{z,y}y).({\hat R}_{z,y}h):=
g'(\eta ).w|_{(\eta ={\hat R}_{z,y}y, w={\hat R}_{z,y}h)}= {\hat
R}_{z,y}[g'(y).h]$ for each $z\in U-y_0$ and $y\in W-y_0$ such that
$Re (z)=Re (y)$ and $z={\hat R}_{z,y}y$ and every $h\in {\cal A}_r$,
where $g'(z)$ is the (super)derivative operator over ${\cal A}_b$.
\par We call such function $(p,r,b)$-quasi-conformal.
If a function $f$ is ${\cal A}_b$ holomorphic on $U$ and satisfies
$(Q1-Q6)$ on $U$ and $f$ is $(p,r,b)$-quasi-conformal on $U\setminus
{\cal S}_A$, where ${\cal S}_A:= \bigcup \{ z+y_0: z={\hat
R}_{z,x}x; x\in A-y_0, z\in U-y_0, Re (z)=Re (x) \} $, while
$A:=A_f:= \{ y+y_0: y\in W-y_0, g'(y)=0 \} $ is a discrete subset in
$W$ consisting of isolated points such that for each $y_1\in A$
there exists $\delta
>0$ for which $\inf_{y\in A, y\ne y_1} |y-y_1|\ge \delta $, then we
call $f$ the $(p,r,b)$-quasi-regular function on $U$. In the
particular latter case of $U={\cal A}_b$ we call $f$ the
$(p,r,b)$-quasi-integral function.
\par If $f$ is a $(p,r,b)$-quasi-regular function on $U\setminus {\cal
S}_C$, where $C$ is a discrete set of isolated points in $W$ at
which $f$ has poles (of finite orders), then we call $f$ the
$(p,r,b)$-quasi-meromorphic function on $U$.
\par For $p=1$ to shorten the notation we shall write that $f$ is
$(r,b)$-quasi-conformal or $(r,b)$-quasi-regular on $U$ or
$(r,b)$-quasi-integral correspondingly. A function of several
variables $\mbox{ }_1z,...,\mbox{ }_nz$ is called
$(p_1,...,p_n;r,b)$-quasi-conformal or
$(p_1,...,p_n;r,b)$-quasi-regular on $U$ open in ${\cal A}_b^n$ or
$(p_1,...,p_n;r,b)$-quasi-integral, if it is
$(p_j,r,b)$-quasi-conformal or $(p_j,r,b)$-quasi-regular or
$(p_j,r,b)$-quasi-integral by $\mbox{ }_jz$ for each $j=1,...,n$. If
$M$ is canonically closed, $M$ is the closure of $U$, then $f$ is
quasi-conformal or quasi-regular if it is such on $U$ and
$f(z)|_{\partial M}$ and $f'(z)|_{\partial M}$ are the continuous
limits of $f$ and $f'$ in $U$, where $\partial M$ is the boundary of
$M$ such that $M\cap {\cal A}_r^n$ is a pseudo-conformal manifold.
\par {\bf 2. Examples.} Let a function $g$ be ${\cal A}_r$ holomorphic
and hence locally analytic. It has a local series expansion of
$g|_{W-y_0}$ with coefficients in ${\cal A}_r$ and the variable
$x\in W-y_0$ such that this series converges on an open ball
$B({\cal A}_r,\xi ,R^-) := \{ x\in {\cal A}_r: |x-\xi |<R \}$ for
each $\xi \in W-y_0$, where $0<R=R(\xi )\le \infty $. Then the
operator ${\hat R}_{z,x}$ acts on $g(x)$ throughout a local series
expansion of $g|_{W-y_0}$ with coefficients in ${\cal A}_r$ and the
variable $x\in W-y_0$, since ${\hat R}_{z,x}$ is the automorphism of
the Cayley-Dickson algebra ${\cal A}_b$. Therefore, $g$ has the
${\cal A}_b$ extensions by the variable $z\in B({\cal A}_b,\xi
,R^-)$ such that $U = \bigcup_{\xi \in W-y_0}B({\cal A}_b,\xi
,R^-)$. Though this extension satisfies Conditions $(Q1-Q6)$, but it
need not be satisfying $(Q7)$ in general.
\par Each $z\in {\cal A}_b$ has the polar decomposition $z=|z|
\exp (Arg (z))$, where $Arg (z)\in {\cal I}_b := \{ y\in {\cal A}_b:
Re (y)=0 \} $ (see Section 3 in \cite{luoyst,lusmfcdv,luhcnfcdv}).
Fix a branch of $Arg (z)$ choosing one definite branch of $Ln$ over
${\cal A}_b$ such that $Arg (z) = Ln (z/|z|)$ and $Arg (z)=0$ for
each real $z\ge 0$. Then
\par $(A)$ $Arg (z)=M\phi =M(z)\phi (z)$, where $M\in {\cal I}_b$,
$|M|=1$, $\phi \in {\bf R}$, $|\phi | = |Arg (z)|$, $M=M(z)$, $\phi
= \phi (z)$.
\par Take without loss of generality $y_0=0$.
For the pair $({\bf C},{\cal A}_b)$ with $2\le b\le 3$ using the
polar decomposition $z- Re (z)=|z- Re (z)|\exp (M\psi )$ for $z\in
{\cal A}_b$, where $Re (M)=0$, $M=M(z-Re (z))\in {\cal A}_b$,
$|M|=1$, $\psi =\psi (z)=\phi (z- Re (z))\in \bf R$, gives a family
of automorphisms ${\hat R}_{z,x}$ for each $z\in {\cal A}_b$ and
every $x\in \bf C$ satisfying the equality
\par $(1)$ ${\hat R}_{z,x}(i_1)=R_{z,x}(i_1)=M$ with $M=M(z-x - Re (z-x))$
and $\phi = \phi (z-x - Re (z-x))$ given by Equation $(A)$ and
${\hat R}_{z,x}(u)=R_{z,x}(u)=u$ for each $u\in \bf R$ and
$R_{z,y}y=y$ for each $z\in \bf C$ and every $y\in {\bf C}$, hence
\par ${\hat R}_{z,x}(x)=|x|\exp (M\phi )$ for each
$x=|x|\exp ({\bf i}\phi )\in \bf C$,\\
where $\phi =\phi (x)\in \bf R$, $z\in {\cal A}_b\setminus {\bf C}$.
Indeed, the algebra isomorphic with ${\cal A}_b$ can be constructed
by the subsequent doubling procedure starting from $M$ as well
instead of ${\bf i}=i_1$ choosing a generator $M_2$ instead of $i_2$
orthogonal to $M$ and taking ${\hat R}_{z,x}(i_2)=M_2$, ${\hat
R}_{z,x}(i_3)=MM_2$, where $M_2$ depends holomorphically on $z$ and
$x$ (see also Proposition 3.2 and Corollary 3.5
\cite{luhcnfcdv,lusmfcdv}), where $\{ i_0, i_1,...,i_{2^r-1} \} $
are generators of ${\cal A}_r$ such that $i_0=1$, $i_j^2=-1$ and
$i_ji_k=-i_ki_j$ for each $1\le j\ne k\le 2^r-1$. Then for $b=3$
take the doubling generator $L\in {\cal I}_b$ orthogonal to $M$,
$M_2$ and $MM_2$ such that $L$ depends holomorphically on $z$ and
$x$ and put $R_{z,x}(i_4)=L$ (see in details below). Write $z\in
{\cal A}_b$ in the form
\par $(2)$ $z=\sum_{s\in \bf b}w_ss$, \\ where $w_s\in \bf R$
for each $s\in {\bf b}:= \{ 1,i_1,...,i_{2^b-1} \} $, $\bf b$ is the
basis of generators of ${\cal A}_b$, put ${\hat b}:={\bf b}\setminus
\{ 1 \} $, hence
\par $(3)$ $z^*=(2^b-2)^{-1} \{ -z +\sum_{s\in \hat b}s(zs^*) \} $.
Therefore,
\par $(4)$ $|z|=(zz^*)^{1/2}=
[z (2^b-2)^{-1} \{ -z +\sum_{s\in \hat b}s(zs^*) \} ]^{1/2}$ and
\par $(5)$ $Re (z) = (z+z^*)/2= \{ (1- (2^b-2)^{-1})z + (2^b-2)^{-1}
\sum_{s\in \hat b}s(zs^*) \} $ are the holomorphic functions on
${\cal A}_b\setminus \{ 0 \} $ in these $z$-representations $(4,5)$.
Then
\par $(6)$ $M(z)\phi (z)=Ln (z/|z|)$ for $z\ne 0$ and for $\phi (z)>0$
with $z\in {\cal A}_b\setminus {\bf R}$ we have
\par $(7)$ $M(z)=Ln (z/|z|)/|Ln (z/|z|)|$ \\
is implied to be written in the $z$-representation with the help of
Formula $(4)$, putting $\phi (z)=0$ for each real non-negative $z$.
In view of Condition $(Q2)$ it is sufficient to consider $\phi
(z)>0$ in the half-space of ${\cal A}_b\setminus {\bf R}$. The
logarithmic function $Ln (z)$ is holomorphic on ${\cal A}_b
\setminus \{ 0 \} $ with the noncommutative non-associative analog
of the Riemann surface described in Section 3.7
\cite{luoyst,lusmfcdv,luhcnfcdv}. In view of Formulas $(4,5)$ the
automorphism ${\hat R}_{z,x}$ given by Equations $(1)$ becomes
holomorphic by $z\in {\cal A}_b$ and by $x\in \bf C$ in the
$(z,x)$-representation.
\par For the pair $({\cal A}_q,{\cal A}_{q+1})$, where $1\le q\in \bf
N$, using the iterated exponent
\par $(8)$ $\exp (\mbox{ }_3M\phi (\xi ))=
\exp \{ \mbox{ }_2M\phi _1(\xi ) \exp (-N\phi _2(\xi )
\exp (-\mbox{ }_2M\phi _3(\xi ))) \} $, \\
where $\xi =|\xi |\exp (\mbox{ }_3M\phi (\xi ))$, $\mbox{ }_3M=
M(\xi )\in {\cal I}_{q+1}$, $\phi _1(\xi ), \phi _2(\xi ), \phi
_3(\xi )\in \bf R$, $z=z_1+i_{2^q}z_2$, $z\in {\cal A}_{q+1}$, $z_1,
z_2\in {\cal A}_q$, $|\mbox{ }_3M| = |\mbox{ }_2M| = |N|=1$, $\xi =
z-x - Re (z-x)$, $N=N(\xi )\perp \mbox{ }_2M=\mbox{ }_2M(\xi )$,
that is, $Re (\mbox{ }_2MN)=0$; $N$ and $\mbox{ }_2M \in {\cal
I}_{q+1}$, $l=i_{2^q}$. Consider $(8)$ for $q=1$ and then for $q=2$.
This gives the family of automorphisms ${\hat R}_{z,x}$ for each
$z\in {\bf O}= {\cal A}_3$ and every $x\in {\bf H}= {\cal A}_2$ such
that
\par $(9)$ ${\hat R}_{z,x}(u) = R_{z,x}(u) =u$ for each $u\in \bf R$,
${\hat R}_{z,x}(i_1)= R_{z,x}(i_1) = \mbox{ }_2M$ and ${\hat
R}_{z,x}(i_2) = R_{z,x}(i_2) =N$ and ${\hat R}_{z,x}(i_3)= R_{z,x}
(i_3) = \mbox{ }_2MN$, $R_{z,y}y=y$ for each $y\in \bf H$ and every
$z\in \bf H$. \par  Indeed, the algebra isomorphic with ${\cal A}_3$
can be constructed starting with $M, N, MN$ instead of $i_1, i_2,
i_3$ and using the doubling procedure and choosing $L\perp {\bf
R}\oplus {\bf R}M\oplus {\bf R}N\oplus {\bf R}MN$, $|L|=1$ (see also
Note 2.4 \cite{luhcnlt,lusmfcdv}).
\par Since $e^M = \cos |M| + M (\sin |M|)/|M|$ for each $M\in {\cal
I}_b\setminus \{ 0 \} $, $e^0=1$, then Equation $(8)$ gives
\par $\mbox{ }_3M\phi (\xi )= \mbox{ }_2M\phi _1\cos \phi _2 +
N \phi _1 \sin \phi _2 \sin \phi _3 +N \mbox{ }_2M\phi _1\sin \phi
_2\cos \phi _3$, hence $\phi _1 = \phi (\xi )$ can be taken. Then
\par $(10)$ $w_s= (-zi_s + i_s(2^b-2)^{-1} \{ -z
+\sum_{k=1}^{2^b-1} i_k (zi_k^*) \} )/2$ for each $s=1,...,2^b-1$.
\\ With the initial conditions $\mbox{ }_2M(0)=i_1$ and $N(0)=i_3$
this gives a family of solutions depending ${\cal
A}_3$-holomorphically on $z$ and ${\cal A}_2$-holomorphically on $x$
with the help of Equations $(3-10)$. For example, take $\mbox{ }_2M
\in i_1{\bf R}\oplus i_5{\bf R}\oplus i_7{\bf R}$, $N\in i_3{\bf
R}\oplus i_4{\bf R}\oplus i_6{\bf R}$. Then choose the doubling
generator $L\in {\cal I}_b$ orthogonal to $\mbox{ }_2M$ and $N$ and
$\mbox{ }_2MN$ such that $L$ depends holomorphically on $z$ and $x$.
In view of Formulas $(2-7)$ the automorphism ${\hat R}_{z,x}$ is
holomorphic by $z\in \bf O$ and by $x\in \bf H$ in the
$(z,x)$-representation. \par This is possible for each $(r,b)$ pair
using the sequence of embeddings ${\cal A}_r\hookrightarrow {\cal
A}_{r+1}\hookrightarrow ... \hookrightarrow {\cal A}_b$ and
considering with the help of $(2-7)$ subsequent holomorphic
solutions of $(8)$ for ${\cal A}_q\hookrightarrow {\cal A}_{q+1}$ in
the corresponding $(z,x)$-representation for each $q=r,...,b-1$. If
$R_{z,x}(i_{2^q})$ are specified for $q=0,1,...,b-1$, then their
multiplication in ${\cal A}_b$ gives $R_{z,x}(i_j)$ for each $1\le
j\le 2^b-1$ (see also \cite{kansol}). This is evident, since ${\cal
A}_b= \{ z\in {\cal A}_b: \exists x\in {\cal A}_{b-1}$ $\mbox{and}$
$\exists T\in SO_{\bf R}(2^b,{\bf R})$ $\mbox{such that}$ $[z]=T[x]
\} $ and $SO_{\bf R}(2^b,{\bf R})$ is the real analytic Lie group
isomorphic with $SO(2^b-1,{\bf R})$, where $[x]\in {\bf R}^{2^b}$,
$[x]=(x_0,x_1,...,x_{2^{b-1}-1},0,0,...)$,
$x=x_0i_0+x_1i_1+...+x_{2^{b-1}-1}i_{2^{b-1}-1}$,
$[z]=(z_0,...,z_{2^b -1})$, $x_j, z_j\in \bf R$ for each $j$, $2\le
b\in \bf N$.
\par From the construction of $R_{z,x}$ it follows, that for the
$({\bf C},{\bf H})$ and $({\bf C},{\bf O})$ pairs, that is, $r=1$
and $b=2, 3$, there exists $R_{z,x}$ satisfying conditions:
\par $(11)$ $R_{vz,wy}=R_{z,y}$ for each $v$ and $w\in {\bf
R}\setminus \{ 0 \} $ such that $vw>0$ and \par $(12)$ $R_{z,y} =
R_{a,x}$ for each $Im (z) =Im (a)$ and $Im (y) = Im (x)$, where $Im
(z) := z - Re (z)$. Therefore, if $R_{z,y}y=z$, then $R_{{\tilde
z},{\tilde y}}{\tilde y}={\tilde z}$.

\par Henceforth, up to an ${\cal A}_b$-pseudo-conformal
diffeomorphism $\xi $ of a domain $U$ such construction of the
family ${\hat R}_{z,x}$ will be implied for $1\le r<b\le 3$,
$R_{z,x}\mapsto {\hat R}_{\xi (z), \xi (x)}$, where $\xi (U)=U$,
$\xi (W)=W$, $\xi ({\bf R}\cap W)={\bf R}\cap W$, $R_{z,x}$ is the
family of this example given by Equations $(1-10)$.
\par Each ${\cal A}_r$-pseudo-conformal (particularly,
complex holomorphic) function with real expansion
coefficients of a power series converging by $x\in W-y_0$ evidently
has an $(r,b)$-quasi-conformal extension due to Condition $(Q3)$.
\par {\bf 2.1. Definition.} For each $p\in {\bf H}={\cal A}_2$ let
\par $(1)$ $E_2(p) := E(p) := p_0+p_1i_1\exp (-p_2i_3\exp (-p_3i_1))$,
while for each $p\in {\bf O}={\cal A}_3$ put \par $(2)$ $E_6(p) :=
E(p)=p_0+p_1i_1\exp (-p_2i_3\exp (-p_3i_1$\\ $\exp (-p_4i_7\exp
(p_5i_1\exp (-p_6i_3\exp (-p_7i_1))...)$,
\\ where $p=p_0i_0+p_1i_1+...+p_{2^b-1}i_{2^b-1}$, $p\in {\cal
A}_b$, $p_0,...,p_{2^b-1}\in \bf R$, $2\le b\le 3$, $i_1i_2=i_3$,
$i_1i_4=i_5$, $i_2i_4=i_6$, $i_3i_4=i_7$, $i_1i_6= -i_7$,
$i_1i_7=i_6$, $i_2i_5=i_7$, $i_2i_7=-i_5$, $i_3i_5=-i_6$,
$i_3i_6=i_5$, $i_5i_6= - i_3$, $i_5i_7=i_2$, $i_6i_7= -i_1$,
$i_ki_l= -i_li_k$ for each $1\le k<l$, $i_k^2=-1$ for each $1\le k$,
$i_0=1$, $z(zy)=(z^2)y$ and $(yz)z=y(z^2)$ for each $z, y\in {\cal
A}_3$; $i_0,...,i_{2^b-1}$ are the standard generators of ${\cal
A}_b$, $\bf R$ is the center of ${\cal A}_b$. It is supposed that
$E_2(p)$ and $E_6(p)$ are written in the $p$-representations over
${\cal A}_2$ and ${\cal A}_3$ respectively with the help of Formulas
2.$(2-5)$.

\par If $f^s$ is an
${\cal A}_b$-holomorphic function on a domain $V$ and $V=E^{-1}(U)$,
where $U$ is a domain in ${\cal A}_b$, $f=f^s\circ E^{-1}$ is
$(p,1,b)$-quasi-conformal or quasi-regular (or quasi-integral) or
quasi-meromorphic on $U$, then we call $f^s$ the
$(p,1,b)$-quasi-conformal or quasi-regular (or quasi-integral for
$V={\cal A}_b$) or quasi-meromorphic function in spherical ${\cal
A}_b$-coordinates on $V$ respectively.
\par Certainly, in Formulas $(1,2)$ other choice of basic
generators or some other order in the iterated exponent can be, but
these formulas provide canonical spherical ${\cal A}_b$-coordinates.
\par {\bf 3. Theorem.} {\it  Let $U\subset {\cal A}_r^n$ be an open
subset, let also $F=(\mbox{ }_1f,...,\mbox{ }_mf): U\to {\cal
A}_r^m$ be a holomorphic mapping, where either $2\le r\in \bf N$ or
$r=\Lambda $, $card (\Lambda )\ge \aleph _0$, $1\le m\le n\in \bf
N$. If $z_0\in U$, $F(z_0)=0$ and the operator $(\partial \mbox{
}_kf/\partial \mbox{ }_jz)_{1\le j, k\le m }$ is invertible at
$z_0$, where $z=(\mbox{ }_1z,...,\mbox{ }_nz)$, $\mbox{ }_jz\in
{\cal A}_r$ for each $j=1,...,n$, then there exist an open
neighborhood $W$ of a point $x_0$ in ${\cal A}_r^m$ and a
neighborhood $V$ of a point $y_0\in {\cal A}_r^{n-m}$ with $W\times
V\subset U$ and a holomorphic mapping $G=(\mbox{ }_1g,...,\mbox{
}_mg): V\to {\cal A}_r^m$ such that $W\cap \{ z\in U: F(z)=0 \} = \{
z=(G(y),y): y\in V \} $ and $g(x_0)=y_0$, where $z_0=(x_0,y_0)$.}
\par {\bf Proof.} Consider the mapping $H=(\mbox{ }_1f,...,
...,\mbox{ }_mf,\mbox{ }_{m+1}z,...,\mbox{ }_nz): U \to {\cal
A}_r^n$. Then the operator $L(z) := (\partial \mbox{ }_kh/\partial
\mbox{ }_jz)_{1\le j, k\le n }$ is invertible at $z_0$, hence it is
invertible in a neighborhood $U_0$ of $z_0$, since $L(z)$ is
super-differentiable, where $(\mbox{ }_1h,...,\mbox{ }_nh)=H$.
Therefore, $L^{-1} (z)$ is super-differentiable in $U_0$, since
$L^{-1}(z)L(z)=L(z)L^{-1}(z)=I$ for each $z\in U_0$, where $I$ is
the unit operator. Then the operators $A(z) := (\partial \mbox{
}_kf/\partial \mbox{ }_jz)_{1\le j, k\le m }$ and $A^{-1}(z)$ are
locally analytic in a neighborhood of $z_0$. Consider the mapping
$q_y(x):=x-A^{-1}(z_0)F(x,y)$ in a neighborhood of $z_0$, where
$(x,y)=z$, $x=(\mbox{ }_1z,...,\mbox{ }_mz)$, $y=(\mbox{
}_{m+1}z,...,\mbox{ }_nz)$.
\par Without loss of generality using shifts we can consider $z_0=0$.
Then $q_y(x)=x$ if and only if $F(x,y)=0$. We have the identity:
$\partial q_y(x)/\partial x=I-A^{-1}(0)(\partial F(x,y)/\partial x)=
A^{-1}(0)(A(0)-\partial F(x,y)/\partial x)$.
From the continuity of $\partial F(x,y)/\partial x$ it follows, that
there exist $a>0$ and $b>0$ such that $ \|\partial q_y(x)/\partial x \|
\le \| A^{-1}(0) \| \| A(0)-\partial F(x,y)/\partial x \| <1/2$
for each $z=(x,y)$ with $ \| x \| <a$ and $ \| y \| <b$.
\par Applying the fixed point theorem to this contracting mapping
$q_y(x)$ we get a solution $G(y)$ in a neighborhood of $0$
(see also the general implicit function
theorem in \S X.7 \cite{zorich} and Theorems II.IV.4.2, 5.1 and 6.1
\cite{grauert}). Then the solution
is locally analytic by $(z,{\tilde z})$, since
$f(z)$ and $A^{-1}(z)$ and $L(z)$ are locally analytic.
Thus in a neighborhood of $(\mbox{ }_{m+1}z_0,...,\mbox{ }_nz_0)$
there are satisfied the identities $\mbox{ }_kf(\mbox{ }_1g,...,
\mbox{ }_mg,\mbox{ }_{m+1}z,...,\mbox{ }_nz)=0$ for $k=1,...,m$
and they are $(z,{\tilde z})$-differentiable and the differentiation by
$\mbox{ }_j{\tilde z}$ gives: \\
$\sum_{l=1}^m (\partial \mbox{ }_kf/\partial \mbox{ }_lz).
(\partial \mbox{ }_lg/\partial \mbox{ }_j{\tilde z}).h+
\sum_{l=1}^m (\partial \mbox{ }_kf/\partial \mbox{ }_l{\tilde z}).
(\partial \mbox{ }_lg/\partial \mbox{ }_jz)^*.h+
(\partial \mbox{ }_kf/\partial \mbox{ }_j{\tilde z}).h=0$ \\
for each $h\in {\cal A}_r$, but $\partial \mbox{ }_kf/\partial
\mbox{ }_l{\tilde z}=\partial \mbox{ }_kf/\partial \mbox{ }_j{\tilde
z}=0$, since $f$ is ${\cal A}_r$-holomorphic and \\ $(\partial
\mbox{ }_kf/\partial \mbox{ }_jz)_{1\le j, k\le m}$ is invertible by
the condition of this theorem, where $z^*={\tilde z}$ denotes the
adjoint of $z$ in the Cayley-Dickson algebra ${\cal A}_r$.
Therefore, $\partial \mbox{ }_lg/\partial \mbox{ }_j{\tilde z}=0$
for each $l=1,...,m$ and $j=m+1,...,n$, consequently, $G$ is
holomorphic.
\par {\bf 4. Corollary.} {\it Let $U$ be an open subset in ${\cal A}_p^n$,
$1\le m\le n\in \bf N$, $F=(\mbox{ }_1f,...,\mbox{ }_mf): U\to {\cal
A}_p^m$ be a $(r,p)$-quasi-conformal mapping, where $1\le r<p\le 3$.
If $z_0\in U$, $F(z_0)=0$ and the operator $(\partial \mbox{
}_kf/\partial \mbox{ }_jz)_{1\le j, k\le m }$ is invertible at
$z_0$, where $z=(\mbox{ }_1z,...,\mbox{ }_nz)$, $\mbox{ }_jz\in
{\cal A}_r$ for each $j=1,...,n$, then there exist an open
neighborhood $W_p$ of a point $x_0$ in ${\cal A}_p^m$ and a
neighborhood $V=V_p$ of a point $y_0$ in ${\cal A}_p^{n-m}$ such
that $(W_p\times V_p)\subset U$ and a holomorphic mapping $G=(\mbox{
}_1g,..., \mbox{ }_mg): V_p\to {\cal A}_p^m$ such that $W_p\cap \{
z\in U: F(z)=0 \} = \{ z=(G(y),y): y\in V_p \} $ with $g(x_0)=y_0$.}
\par {\bf Proof.} Since $F$ is $(r,p)$-quasi-conformal, then
it is holomorphic on $U$ satisfying Conditions $(Q1-Q7)$ with
$W=U\cap {\cal A}_r^n\subset U$, where ${\cal A}_r\hookrightarrow
{\cal A}_p$ is the natural embedding. In view of Theorem 3 there
exists a holomorphic solution of this theorem.
\par {\bf 5. Corollary.} {\it Let $F$ satisfies conditions of Corollary
4 and $n=2$ and $m=1$. Then $G$ is $(r,p)$-quasi-conformal in a
neighborhood of $y_0$ at each point $y\in V=V_p$ such that
$F(G(y),y)=0$.}
\par {\bf Proof.} In view of Corollary 4 $G$ is holomorphic, hence
satisfies Condition $(P1)$. We have that \par $(1)$ $G'(y).h = -
({F'}_x(x,y))^{-1}.[({F'}_y(x,y)).h]$ for all $x=G(y)$ and each
$h\in {\cal A}_p$, when $F(G(y),y)=0$, since the quaternion skew
field ${\bf H}={\cal A}_2$ is associative and the octonion algebra
${\bf O}={\cal A}_3$ is alternative. \\
The restriction of $F$ on $U\cap {\cal A}_r^2$ is pseudo-conformal,
hence ${F'}_x(x,y)$ and ${F'}_y(x,y)$ for $(x,y)\in U\cap {\cal
A}_r^2$ satisfy Conditions $(P2,P3)$. In view of Theorem 2.4
\cite{lusmnfgpcd} \par $(2)$ ${F'}_x(x,y).h=a(x,y)hb(x,y)$ and
${F'}_y(x,y).h=c(x,y)he(x,y)$ \\ for each $h\in \bf H$ and each
$(x,y)\in U\cap {\cal A}_r^2$, for $r=2$, where $a(x,y), b(x,y),
c(x,y), e(x,y)$ are non-zero ${\cal A}_r$-holomorphic functions on
$U\cap {\cal A}_r^2$. For $r=1$, over $\bf C$, evidently due to the
commutativity of $\bf C$ we take as usually $b=1$ and $e=1$.
\par Therefore, from Equation $(1)$ it follows, that the restriction of
$G'(y)$ on $V\cap {\cal A}_r$ satisfies Conditions $(P2,P3)$, since
the quaternion skew field $\bf H$ and the complex field $\bf C$ are
associative. But ${\hat R}_{z,x}$ are automorphisms of ${\cal A}_p$
such that Conditions $(Q1-Q6)$ are satisfied. For simplicity of the
notation take the zero marked point. We have $\bigcup_{s\in V} \{
{\hat R}_{s,y}y: y\in V\cap {\cal A}_r \} =V$, hence $\{ q\in V:
\exists y\in V\cap {\cal A}_r$ $\mbox{such that}$ ${\hat R}_{q,y}y=q
\} =V$ (see also $(Q2)$). Let $q\in V\setminus {\cal A}_r$ and $y\in
V\cap {\cal A}_r$ be such that ${\hat R}_{q,y}y=q$, then ${\hat
R}_{q,y}F(x,y)=F({\hat R}_{q,y}x,q)$, but $F(\zeta ,q)=0$ in $W$ is
equivalent to $(\zeta ,q)=(G(q),q)$ and $q\in V$. Therefore, ${\hat
R}_{q,y} F(G(y),y) = F(G(q),q)=0$ and ${\hat
R}_{q,y}G'(y).h=G'(q).({\hat R}_{q,y}h)$ due to $(1)$ for each $h\in
{\cal A}_r$ and every $y\in V\cap {\cal A}_r$ and every $q\in V$
such that ${\hat R}_{q,y}y=q$ and $Re (y)= Re (q)$, consequently,
${\hat R}_{q,y}G(y)=G(q)$ for each $(q,y)\in V\times (V\cap {\cal
A}_r)$ such that ${\hat R}_{q,y}y=q$ and $Re (y)= Re (q)$, since
${F'}_x(x,y)$ is invertible for $x=G(y)$ and $F$ is locally analytic
and using expansion of $F(x,y)$ by $(x,y)$ with $x=G(y)$.
\par Put $H=G'(y).h$, then $h=(G'(y))^{-1}.H$, since $F$ is
pseudo-conformal and $ker F'_x(x,y)=0$ and $ker F'_y(x,y)=0$ for
each $(x,y)\in U\cap {\cal A}_r^2$, while ${\cal A}_p$ is the finite
dimensional algebra over $\bf R$. If the operator ${F'}_x(x,y):
{\cal A}_p\to {\cal A}_p$ satisfies Condition $(Q7)$ and has the
inverse, then its inverse also satisfies $(Q7)$, since the
restriction ${F'}_x(x,y)|_{{\cal A}_r}$ has Form $(2)$ and each
non-zero number in ${\cal A}_p$ is invertible, where $1\le p\le 3$.
Since the right hand side of Equation $(1)$ satisfies Condition
$(Q7)$, then the left hand side of it satisfies $(Q7)$ as well.
\par {\bf 6. Corollary.} {\it If $f: U\to {\cal A}_p$ is a
$(r,p)$-quasi-conformal function, where $U$ is open in ${\cal A}_p$
and $f(U)=V$ is open in ${\cal A}_p$ and $f$ is bijective on $U$,
then $f^{-1}: V\to {\cal A}_p$ is $(r,p)$-quasi-conformal.}
\par {\bf Proof.} Take the function $F(x,y)=f(x)-y$, then it is
$(r,p)$-quasi-conformal and satisfies Conditions of Lemma 5. Since
$f(U)=V$ and $f: U\to V$ is bijective, then there exists $g=f^{-1}:
V\to U$, which is $(r,p)$-quasi-conformal due to Lemma 5.
\par {\bf 7. Theorem.} {\it Let $f$ and $g$ be $(p,r,b)$- and
$(q,r,b)$-quasi-conformal mappings on neighborhoods $U$ of $z_0$ and
$V$ of $y_0$ respectively such that $f(U)\supset V$ and
$f(z_0)=y_0$, where $1\le r< b\le 3$, $y_0$ and $z_0\in {\cal A}_r$,
then their composition $g\circ f$ is $(pq,r,b)$-quasi-conformal on a
neighborhood $W$ of $z_0$.}
\par {\bf Proof.} The composition of pseudo-conformal mappings is
pseudo-conformal, so in accordance with Definition 1 it is
sufficient to take the neighborhood $W$ of $z_0$ such that $W =
f^{-1}(V)$ is open, since $f$ is continuous (see \cite{ludagpd} and
Theorem 2.6 \cite{lusmnfgpcd}). Therefore, $g\circ f$ is
$pq$-pseudo-conformal at each point in $W\cap {\cal A}_r$, since
$y_0$ and $z_0\in {\cal A}_r$. The composition of holomorphic
mappings is holomorphic, the composition $T_1T_2$ of proper elements
$T_1, T_2\in SO(2^b,{\bf R})$ is the proper element $T_1T_2\in
SO(2^b,{\bf R})$, since $SO(2^b,{\bf R})$ is the special orthogonal
group. If rotations $T_1$ and $T_2$ have a common axis, then their
composition preserves this axis, hence $SO_{\bf R}(n,{\bf R})$ is
the subgroup of $SO(n,{\bf R})$. Take the families of automorphisms
${\hat R}^g$ and ${\hat R}^f$ for $g$ and $f$ correspondingly in
accordance with Definition 1. Therefore, the composition ${\hat
R}^{g\circ f}_{z,x} := {\hat R}^g_{({\hat
R}^f_{z,x}(f(x+z_0)-y_0)),(f(x+z_0)-y_0)}$ is defined for each $x\in
(W-z_0)\cap {\cal A}_r$ and every $z \in W-z_0$ and it gives the
restriction ${\hat R}_{z,x}=id$ for each $x$ and $z\in (W-z_0)\cap
{\cal A}_r$, since $f(U\cap {\cal A}_r)\subset {\cal A}_r$. Thus the
family of operators ${\hat R}^{g\circ f}(z,x)$ satisfies Conditions
$(Q2-Q5)$. Therefore, $g\circ f(z+z_0)={\hat R}^{g\circ f}_{z,x}
h\circ f(x+z_0)$ for each $x\in (W-z_0)\cap {\cal A}_r$ and every
$z\in W-z_0$ with $Re (z)=Re (x)$ and $z={\hat R}_{z,x}x$, where
$h(y):=g(y+y_0)$, consequently, $g\circ f$ satisfies $(Q6)$. Since
$(g\circ f)'(z+z_0).h= g'(f(z+z_0)).(f'(z+z_0).h)$ for each $z\in
W-z_0$ and $h\in {\cal A}_b$, while $g$ and $f$ satisfy $(Q6,Q7)$
with ${\hat R}^g$ and ${\hat R}^f$ respectively, then $(g\circ f)'$
satisfies $(Q7)$ with ${\hat R}^{g\circ f}$ and inevitably $g\circ
f$ is $(pq,r,b)$-quasi-conformal on $W$.
\par {\bf 8. Corollary.} {\it Let $U$ be an open domain in
${\cal A}_b$ with a marked point $y_0\in U\cap {\cal A}_r$, $1\le
r<b\le 3$, then the family of all $(r,b)$-quasi-conformal
diffeomorphisms $f$ of $U$ onto $U$ preserving a marked point $y_0$
has the group structure.}
\par {\bf Proof.} In Accordance with Theorem 7 compositions
of $(r,b)$-quasi-conformal mappings $g, f$ are
$(r,b)$-quasi-conformal, since $f(y_0)=y_0$ and $g(y_0)=y_0$. In
view of Corollary 6 the inverse mapping of $f$ is also
$(r,b)$-quasi-conformal. Evidently, the identity mapping $id (x)=x$
for each $x\in U$ is pseudo-conformal, hence it is
$(r,b)$-quasi-conformal. Since $f\circ id=id\circ f=f$ for each
homeomorphism $f: U\to U$, then $id=e$ is the unit element of the
family of $(r,b)$-quasi-conformal diffeomorphisms.
\par {\bf 8.1. Remark.} Topologize the family $H(M,P)$ of ${\cal A}_b$
holomorphic mappings from a domain $M$ in ${\cal A}_b^n$ into a
domain $P$ in ${\cal A}_b^k$ by the compact-open topology of locally
analytic mappings as in the proof of Theorem 3.18 \cite{lusmnfgpcd},
where $n, k\in \bf N$. This topology on $H(M,P)$ induces the
topology on the group of ${\cal A}_b$ holomorphic diffeomorphisms
$DifH(M)$ of $M$. For the family of $(r,b)$-quasi-conformal
diffeomorphisms $f$ of $M$ suppose that $f(M\cap {\cal A}_r^n)=M\cap
{\cal A}_r^n$.
\par {\bf 8.2. Theorem.} {\it The family of all $(r,b)$-quasi-conformal
diffeomorphisms $DifQ(M)$ of a compact canonical closed domain $M$
in ${\cal A}_b^n$ preserving a marked point $y_0$ (see 8.1), $1\le
r<b\le 3$, $n\in \bf N$, form the topological metrizable group,
which is complete relative to its metric and locally compact. The
group $DifQ(M)$ is the analytic Lie group over $\bf R$.}
\par {\bf Proof.} In view of Theorems 3.24,25  \cite{lusmnfgpcd}
the group $DifP(M\cap {\cal A}_r^n)$ of all pseudo-conformal
diffeomorphisms of $M\cap {\cal A}_r^n$ is the topological
metrizable locally compact analytic Lie group over $\bf R$. Consider
$y_0=0$ without loss of generality. On the other hand, each $f\in
DifQ(M)$ is obtained from the corresponding $q\in DifP(M\cap {\cal
A}_r^n)$ with the help of operators ${\hat R}_{z,x}$ in accordance
with Conditions $(Q1-Q7)$ of Definition 1. In its turn each ${\hat
R}_{z,x}$ is the automorphism of ${\cal A}_b$ depending ${\cal A}_b$
and ${\cal A}_r$ holomorphically on $z$ and $x$ respectively (see
Example 2). Thus $f^{-1}$ is obtained from $q^{-1}$ with the help of
${\hat R}_{\zeta ,y}$, $\zeta =f(z)$, $y=q(x)$, but $q^{-1}$ is
pseudo-conformal, hence $f^{-1}$ is $(r,b)$-quasi-conformal due to
Conditions $(Q6,Q7)$ for $f$, since $(f^{-1})'=(f')^{-1}$ and each
${\hat R}_{z,x}$ is invertible. Since for each $f_1, f_2\in DifQ(M)$
we have the corresponding $q_1, q_2\in DifP(M\cap {\cal A}_r^n)$ and
$DifP(M\cap {\cal A}_r^n)$ is the group, then $f_1\circ f_2\in
DifQ(M)$ due to $(Q1-Q7)$. The group $DifP(M\cap {\cal A}_r)$ is the
finite-dimensional locally compact analytic Lie group over $\bf R$
and the family ${\hat R}_{z,x}$ also form the finite-dimensional
analytic family over $\bf R$, hence $DifQ(M)$ is the
finite-dimensional analytic Lie group over $\bf R$ and inevitably it
is locally compact, metrizable and complete.
\par {\bf 9.1. Proposition.} {\it If $q_1$ and $q_2$ are holomorphic
functions on a domain $W$ in ${\bf C}$,  $q_1$ and $q_2$ have
$(1,b)$-quasi-conformal extensions $f_1$ and $f_2$ on $U$, $1<b\le
3$, with the same family ${\hat R}_{z,x}$ and the same marked point
$y_0\in W$ for $f_1$ and $f_2$ (see $(Q1-Q7)$ in Definition 1) and
$(q_1q_2)(x)'\ne 0$ at each point $x\in W$, then their product
$q_1q_2$ has the $(1,b)$-quasi-conformal extension $f_1f_2$ on $U$.}
\par {\bf Proof.} Conditions $(Q1-Q6)$ are evidently satisfied
for $f_1f_2$, since ${\hat R}_{z,x}$ is the automorphism of ${\cal
A}_b$ for each $z\in U-y_0$ and $x\in W-y_0$. On the other hand, the
product of complex holomorphic functions is complex holomorphic, the
product of ${\cal A}_b$ holomorphic functions is ${\cal A}_b$
holomorphic. Conditions $(Q6,Q7)$ are satisfied for $f_1$ and $f_2$,
hence
\par $(g_1g_2)'({\hat R}_{z,y}y).({\hat R}_{z,y}h):=
g_1(z)[(g_2)'({\hat R}_{z,y}y).({\hat R}_{z,y}h)] + [(g_1)'({\hat
R}_{z,y}y).({\hat R}_{z,y}h)]g_2(z)$ \\  $= {\hat
R}_{z,y}[(g_1g_2)'(y).h]$
\\ for each $z\in U-y_0$ and $y\in W-y_0$ such that $Re (z)=Re (y)$
and $z={\hat R}_{z,y}y$ and every $h\in {\cal A}_r$, where $g'(z)$
is the (super)derivative operator over ${\cal A}_b$, consequently,
$(Q7)$ is satisfied for $f_1f_2$.
\par {\bf 9.2. Corollary.} {\it Let $q_1$ and $q_2$ be holomorphic
functions on a domain $W$ in $\bf C$ with isolated zeros of $q_1'$
and $q_2'$, $q_1$ and $q_2$ have $(1,b)$-quasi-regular extensions
$f_1$ and $f_2$ on $U$, $1<b\le 3$, with the same family ${\hat
R}_{z,x}$ and the same marked point $y_0\in W$ for $f_1$ and $f_2$
(see $(Q1-Q7)$ in Definition 1), then their product $q_1q_2$ has the
$(1,b)$-quasi-regular extension $f_1f_2$ on $U$.}
\par {\bf Proof.} In view of Theorem 9.1 $f_1f_2$ is
$(1,b)$-quasi-conformal on $U\setminus S_A$, where $S_A = \{ z+y_0:
z={\hat R}_{z,x}x,$ $z\in U-y_0, x\in W-y_0, q_1'(x+y_0)=0$
$\mbox{or}$  $q_2'(x+y_0)=0 \} $, $S_A\cap W$ is the discrete subset
consisting of isolated points in $W$. Conditions $(Q1-Q6)$ are
satisfied for $f_1f_2$, since each ${\hat R}_{z,x}$ is the
automorphism of ${\cal A}_b$. Thus $f_1f_2$ is $(1,b)$-quasi-regular
on $U$.
\par {\bf 9.3. Remark.} In general Theorem 9.1 and Corollary 9.2
may be not true for $((a_1f_1)(a_2f_2))$ instead of $f_1f_2$, when
$f_j$ are taken with constant non-real multipliers $a_j\in {\bf
C}\setminus \bf R$ (see also Notes 13).
\par {\bf 9.4. Theorem.} {\it Let $q_n$ be a sequence of complex
holomorphic functions on an open connected convex domain $W$ in $\bf
C$ such that the series $\sum_{n=1}^{\infty }q_n'(y)$ converges
uniformly on $W$ to a function $q'(y)$ with $q'(y)\ne 0$ for each
$y\in W$ (or $y\in W\setminus A$ with discrete subset $A$ consisting
of isolated points in $W$) and $\sum_{n=1}^{\infty } q_n(y_0)$
converges at a marked point in $W$ to $q(y_0)$ while each $q_n$ has
a $(1,b)$-quasi-conformal (or $(1,b)$-quasi-regular) extension $f_n$
on a domain $U$ in ${\cal A}_b$ with the same family $\{ {\hat
R}_{z,x}: z\in U-y_0, x\in W-y_0 \} $. Then the series
$\sum_{q=1}^{\infty }q_n(y)$ converges on $W$ to a function $q(y)$
which has a $(1,b)$-quasi-conformal (or $(1,b)$-quasi-regular
correspondingly) extension $f$ on $U$.}
\par {\bf Proof.} In view of Theorem XVI.3.4 \cite{zorich}
the series $\sum_{q=1}^{\infty }q_n(y)$ converges on $W$ to a
function $q(y)$ and this convergence is uniform on compact subsets
of $W$. Since $q'(y)=\sum_{n=1}^{\infty }q_n'(y)$ on $W$, then there
exist $\partial q(y)/\partial y_1$ and $\partial q(y)/\partial y_2$,
where $y=y_1+{\bf i} y_2$, $y_1, y_2\in \bf R$, ${\bf
i}=(-1)^{1/2}$. Consequently, there exists $\partial q(y)/\partial
{\bar y}=\sum_{n=1}^{\infty }\partial q_n(y)/\partial {\bar y} =0$
on $W$, since each $q_n$ is holomorphic on $W$. Since ${\hat
R}_{z,x}$ is the automorphism of ${\cal A}_b$ depending
holomorphically on $z\in U-y_0$ and $x\in W-y_0$, then the series
$\sum_{n=1}^{\infty }f_n(y)$ and $\sum_{n=1}^{\infty }f_n'(y)$
converge on $U$ to $f(y)$ and $f'(y)$ respectively and this
convergence is uniform on $P$ and $P\times B$ for each compact
subset $P$ in $U$, where $B=B({\cal A}_b,0,1) := \{ z\in {\cal A}_b:
|z|\le 1 \} $. Consequently, there exists $\partial f(y)/\partial
{\tilde y}=0$, since $\partial f_n(y)/\partial {\tilde y}=0$ for
each $n\in \bf N$ and inevitably $f(y)$ is ${\cal A}_b$ holomorphic
on $U$. Conditions $(Q1-Q6)$ are satisfied for each $f_n$ on $U$ and
$(Q7)$ on $U\setminus S_{A_n}$, where $S_{A_n}=\emptyset $ in the
$(1,b)$-quasi-conformal case, hence $(Q1-Q6)$ are satisfied for
$f(y)$ on $U$ and $(Q7)$ on $U\setminus S_A$, where $A=\emptyset $
in the $(1,b)$-quasi-conformal case, since the series
$\sum_{n=1}^{\infty }f_n(y)$ and $\sum_{n=1}^{\infty }f_n'(y)$
converge on $U$ to $f(y)$ and $f'(y)$ correspondingly and
$f'(y).h\ne 0$ for each $(y,h)\in (U\setminus S_A)\times ({\cal
A}_b\setminus \{ 0 \} )$.
\par {\bf 9.5. Examples. 1.} Let $q_n(y)=c_n(y-y_0)^n$, coefficients
be real $c_n\in \bf R$, $\sum_{n=1}^{\infty }|c_n|R^n <\infty $ for
each $n\in \bf N$, where $0<R<\infty $, $W= \{ y\in {\bf C}:
|y-y_0|<R \} $, $y_0\in \bf C$ such that $\sum_{n=1}^{\infty
}c_nn(y-y_0)^{n-1}\ne 0$ on $W$, put $U= \{ z\in {\cal A}_b:
|z-y_0|< R \} $ and take $ \{ {\hat R}_{z,x} \} $ from Examples 2.
Then conditions of Theorem 9.4 are satisfied, since ${\hat
R}_{z,x}c_n=c_n$ for each $n$. In particular, take $c_n=1$ for each
$n$, $0<R<1$, $y_0=0$, then $f(y)=1/(1-y)$ and $q'(y)=(1-y)^{-2}\ne
0$ on $W$. If $q(y)=\sin (y)$ or $q(y)=\cos (y)$, then $q$ has the
$(1,b)$-quasi-integral extension. Evidently, if $a=const\ne 0$,
$a\in {\bf R}\setminus \{ 0 \} $, $f(z)$ is $(r,b)$-quasi-regular on
$U$ or $(r,b)$-quasi-integral, then $q(z)=f(az)$ is
$(r,b)$-quasi-regular on $U/a = \{ z/a: z\in U \} $ or
$(r,b)$-quasi-integral correspondingly.

\par {\bf 2.} Take $q_n(y)=c_n\exp (a_n(y-y_0))$, where
$c_n\in \bf R$, $a_n\in \bf R$, $a_n\ne 0$ for each $n\ge 2$,
$\sum_{n=1}^{\infty }c_n$ converges and $\sum_{n=1}^{\infty }
|a_nc_n|\exp (|a_n|R) <\infty $ such that $\sum_{n=1}^{\infty } a_n
c_n \exp (a_n(y-y_0))\ne 0$ on $W := \{ y\in {\bf C}: |y-y_0|<R \}
$, $0<R<\infty $, $y_0\in \bf C$, $U:= \{ z\in {\cal A}_b: |z-y_0|<R
\} $ and ${\hat R}_{z,x}$ is from Examples 2. Then Conditions of
Theorem 9.4 are satisfied, since ${\hat R}_{z,x}c_n=c_n$ and ${\hat
R}_{z,x}a_n=a_n$ for each $n$, while $\exp (a(y-y_0))$ is ${\cal
A}_b$-pseudo-conformal on ${\cal A}_b$ for $a\ne 0$ in ${\cal
A}_b\setminus \{ 0 \} $ (see \cite{ludagpd,lusmnfgpcd}). \par In
particular, for $c_n=1$ and $q_n(y)=n^{-y}=\exp (-y\ln n)$ for each
 $n\ge 1$, $y_0=0$, the series $\sum_{n=1}^{\infty } q_n(y)$ and
$\sum_{n=1}^{\infty }q_n'(y)$ converge uniformly on $W_R := \{ y\in
{\bf C}: Re (y)>R \} $ for $1<R<\infty $ to the holomorphic function
$\zeta (y)$ and put $U_R := \{ z\in {\cal A}_b: Re (z)>R \} $.
Therefore, take $W= \{ y\in W_R: \zeta '(y)\ne 0, 1<R<\infty \} $
and $U=U_1\setminus \bigcup \{ {S'}_y: y\in W_1\setminus W \} $,
where ${S'}_y := \{ z: z={\hat R}_{z,y}y; z\in U_1, Re (z)=Re (y) \}
$ for $y\in W_1\setminus W$. Thus, $\zeta (y)$ has the
$(1,b)$-quasi-conformal extension on $U$ from $W$ for $b=2$ and for
$b=3$. Since the derivative $\zeta '(y)$ is holomorphic on $W_1$
with isolated zeros, then $\zeta (y)$ has the $(1,b)$-quasi-regular
extension on $U_1$ from $W_1$ for $b=2$ and for $b=3$. \par If take
\par $(1)$ $q_n(y)= c_n\exp (v_n E(t_n(y-y_0)))$\\ with $v_nt_n=a_n$,
where $v_n, t_n\in {\bf R}\setminus \{ 0 \} $, then these examples
provide $(1,b)$-quasi-conformal or quasi-regular extensions in
spherical ${\cal A}_b$-coordinates (see Definition 2.1). For this
choose the family $R_{E(t(p-y_0)), E(t(y-y_0))}$ for the $({\bf
C},{\cal A}_b)$ pair in Example 2 independent from $t\in {\bf
R}\setminus \{ 0 \} $. This is possible due to additional Conditions
$2.(11,12)$, since $\exp (p_0+p_SS)=\exp (p_0) (\cos (p_S) + S \sin
(p_S))$ for each $p\in {\cal A}_b$ and $sign (\sin (p_st) \sin
(p_yt))=sign (p_sp_y)$ for each $t\in {\bf R}\setminus \{ 0 \} $,
where $p_0=Re (p-y_0)$, $p_SS= Im (p-y_0):=p-y_0 - Re (p-y_0)$,
$p_0, p_S\in \bf R$, $S\in {\cal I}_b$, $|S|=1$, $sign (t)=1$ for
$t>0$, $sign (t)=-1$ for $t<0$ and $sign (0)=0$. Indeed, the group
$SO_{\bf R}(2^b,{\bf R})$ is isomorphic with $SO(2^b-1,{\bf R})$ and
$E(t(p-y_0))\in {\bf R}\oplus S{\bf R}$ for each $t\in \bf R$.
Therefore,
\par $(2)$ $R_{E(t(p-y_0)), E(t(y-y_0))}= R_{E(p-y_0), E(y-y_0)}$ \\
for each $p\in {\cal A}_b$, $y\in \bf C$ and each real nonzero $t\in
{\bf R}\setminus \{ 0 \} $, where $y_0\in \bf C$ is the marked
point. Then
\par $(3)$
$R_{E(p-y_0), E(y-y_0)}\sum_{n=1}^{\infty }q_n(y)=$ \\
$\sum_{n=1}^{\infty }R_{E(t_n(p-y_0)), E(t_n(y-y_0))} \exp (v_n
E(t_n(y-y_0)))=$ \\
$\sum_{n=1}^{\infty }\exp (v_n E(t_n(p-y_0)))=f(p)$ \\
for each $Re (E(p-y_0))=Re (E(y-y_0))$ with $R_{E(p-y_0), E(y-y_0)}
E(y-y_0)=E(p-y_0)$, since $R_{z,x}(tx)=tR_{z,x}x$ for each $t\in \bf
R$ and $E(y)=y$ for each $y\in \bf C$. Thus $q(y)=\sum_{n=1}^{\infty
}q_n(y)$ has the $(1,b)$-quasi-conformal in spherical ${\cal
A}_b$-coordinates extension $f$ with $q$ on $V$ such that $E(V)=W$
and $f$ on $P$ such that $E(P)=U$ choosing the corresponding
branches of $Ln$.
\par Henceforth, it is supposed that Condition $(2)$ is satisfied
in spherical ${\cal A}_b$-coordinates.
\par {\bf 3.} It is known, that the gamma function
$\Gamma (z)$ is holomorphic on ${\bf C}\setminus \{ 0, -1, -2,
-3,... \} $ having poles of the first order at points $0, -1, -2,
-3,...$ with residues $res_{z=-n} \Gamma (z)= (-1)^n/n!$, $n=0, 1,
2,...$. Moreover, \par $(1)$ $1/\Gamma (z+1)=e^{Cz}
\prod_{k=1}^{\infty } [(1+z/k) e^{-z/k}]$ and this product converges
everywhere on $\bf C$, where $C=\lim_{n\to \infty }(\sum_{k=1}^n1/k
- \ln n)=0.5772157...$ is the Euler constant (see \S VII.1
\cite{lavrshab}). In the product of Formula $(1)$ all coefficients
are real. It is possible to consider for this function different
holomorphic extensions over ${\cal A}_b$ (see Section 4 in
\cite{lusmfcdv}). Applying operators ${\hat R}_{z,y}$ from Example 2
with $y_0=0$ and Proposition 9.1 and Theorem 9.4 to Equation $(1)$
provides the $(1,b)$-quasi-meromorphic extension of $\Gamma (z)$,
which is $(1,b)$-quasi-conformal  on ${\cal A}_b\setminus \{ S_0,
S_{-1}, S_{-2},... \} $ for $2\le b\le 3$, where
$S_{-n}=S_{-n}^{\Gamma }$. In particular, for $y_0=0$ we have
$S_{-n}^{\Gamma } = \{ - n \} $, since rotations are around the real
axis, $T={\hat R}_{z,y}\in SO_{\bf R}(2^b,{\bf R})$. Moreover,
$1/\Gamma (z)$ is $(1,b)$-quasi-integral, hence $\Gamma (z)$ has not
zeros in ${\cal A}_b$.

\par {\bf 10. Definition.} Let $a_1,...,a_n, z\in {\cal A}_r$,
put $Exp_1(a_1;z) := \exp (a_1z)$, $Exp_n(a_1,...,a_n;z) :=
Exp_{n-1}(a_1,...,a_{n-1};Exp_1(a_n;z))$ for $n>1$, where $2\le r$.
For $a_1\ne 0,...,a_n\ne 0, z\ne 0$ put $Ln_1(a_1;z):=a_1^{-1}Ln
(z)$, $Ln_n(a_1,...,a_n;z):=Ln_{n-1}(a_1,...,a_{n-1};Ln_1(a_n;z))$
for $n>1$, where $Exp_0(z):=id(z)=z$ and $Ln_0(z):=id(z)=z$ for each
$z\in {\cal A}_r$, $1\le r$. Here $a_1,...,a_{n-1}$ can be
constants, but more generally ${\cal A}_r$-pseudo-conformal
functions $a_1(z)\ne 0$,...,$a_{n-1}(z)\ne 0$, $a_n\ne 0$ is a
constant in ${\cal A}_r$, $2\le r$.
\par Suppose that $\gamma (t) := z_0 + \rho Exp_n(a_1,...,a_n;\xi (t))$
is a curve in an open domain $U$ in ${\cal A}_r$ and $f$ is a
holomorphic function $f: U\to {\cal A}_r$, where $a_1=a_1(t)\ne
0,...,a_{n-1}=a_{n-1}(t)\ne 0$ are constants or pseudo-conformal
functions with values in ${\cal A}_r$ on an open domain $V_a\supset
[0,1]$ or $V_a\supset \bf R$ in ${\cal A}_r$, $a_n\ne 0$, $\xi (t)$
is a rectifiable curve in ${\cal A}_r$, $t\in [0,1]\subset \bf R$,
$f(z)\ne 0$ for each $z=\gamma (t)$, where $0<\rho <\infty $. Then
put \par $\Delta _{\gamma } Arg_nf := \Delta _{\gamma }
Arg_n(a_1,...,a_n;f) := \int_{z\in \gamma }
dLn_n(a_1,...,a_{n-1},1;f(z))$ \\ with a chosen branch of $Ln$.
\par {\bf 10.1. Note.} For $n=1$ and $\xi (t)=t$ with $M\in
{\cal A}_r$, $Re (M)=0$, $|M|=1$, $a_n=2\pi M$, $\{ \gamma (t): t\in
[0,1] \} $ is the circle. If $n=1$ and $a_1=2\pi M$, then $\Delta
_{\gamma }Arg_1(f) = \Delta _{\gamma }Arg (f)$ is the usual change
of the argument of a function $f$ along a curve $\gamma $ (see also
Section 3 in \cite{luoyst,lusmfcdv,luhcnfcdv} and Theorem 2.23
\cite{lusmnfgpcd}).
\par {\bf 10.2. Proposition.} {\it The logarithmic function
$Ln$ on ${\cal A}_r\setminus \{ 0 \} $, where $1\le r\le \infty $,
has a countable number of branches.}
\par {\bf Proof.} For $r=1$ we have ${\cal A}_1=\bf C$ and in this
case the statement of this proposition  is well-known. So consider
$2\le r\le \infty $. \par Each nonzero $z\in {\cal A}_r\setminus \{
0 \} $ can be written in the polar form \par $(1)$ $z=|z|\exp (M\phi
+2\pi nM)$, \\ where $M\in {\cal I}_r := \{ z\in {\cal A}_r: Re
(z)=0 \} $, $|M|=1$, $\phi \in [0,2\pi )$, $n\in \bf Z$, $Arg (z)=
M\phi +2\pi nM$ (see Section 3 in \cite{luoyst,lusmfcdv,luhcnfcdv}).
If $K\in {\cal I}_r$, $|K|=1$, $K$ is not parallel to $M$, that is,
$|Re (MK^*)|<1$, then $M$ and $K$ do not commute. When $0<\phi <\pi
$, then $\exp (M\phi + \pi Ks )\ne \exp (M\phi +\pi nM)$ for each $s
\ne 0$, $s \in {\bf Z}\setminus \{ 0 \} $, and each $n\in \bf Z$,
since $\exp (M\phi + \pi Ks ) = \cos |M\phi + \pi Ks | + (M\phi +
\pi Ks ) (\sin |M\phi + \pi Ks|)/ |M\phi + \pi Ks |$ while $\exp
(M\phi +\pi nM) = \cos |\phi + \pi n| + M(\phi + \pi n) (\sin |\phi
+ \pi n|)/ |\phi + \pi n|$ and $|M\phi + \pi Ks |^2=\phi ^2 +(\pi
s)^2 +2 Re (MK^*)\phi \pi s$ and $(\phi +\pi n)\notin \pi \bf Z$ and
$K$ is not parallel to $M$, where ${\bf Z} := \{..., -2, -1, 0, 1,
2,... \} $. On the other hand, $Im (z) := z- Re (z)$ is parallel to
$M$ in Equation $(1)$, hence the only solutions of $(1)$ are $Arg
(z) = M(\phi +2\pi n)$, where $\phi \in [0,2\pi )$, $n\in \bf Z$,
$M$ is parallel to $Im (z)$. Therefore, $Ln$ has only countable
number of branches which can be enumerated by $n\in \bf Z$.
\par In more details it is possible to construct the following
noncommutative analog of the Riemann $2^r$-dimensional surface $\cal
R$ of $Ln$ such that $Ln :{\cal A}_r\setminus \{ 0 \} \to \cal R$ is
the univalent mapping, where $2\le r\le \infty $, $2^{\infty
}=\infty $. Consider copies $({\cal A}_r, n i_1 )$ of ${\cal A}_r$
embedded into ${\cal A}_r^2$, where $ni_1\in i_1\bf R$ is in the
second multiple ${\cal A}_r$, $n\in \bf Z$. Put $P_j := \{ z\in
{\cal A}_r: z_0< 0, z_j=0, z_k\in {\bf R} \forall 0< k\ne j \} $ and
consider sections of ${\cal A}_r$ (of the first multiple) by $P_j$
for each $1\le j\in \bf Z$, where $z=\sum_{j=0}^{2^r-1}z_ji_j$,
$z_j\in \bf R$, $i_j$ are standard generators of ${\cal A}_r$. Then
the set $\{ z\in {\cal A}_r: z_0<0 \} $ is partitioned into the
subsets $S(k_1,k_2,...)$ corresponding to definite combinations of
signs of $z_j$: either $z_j\ge 0$ or $z_j\le 0$ with $k_j=1$ and
$k_j=-1$ respectively. For finite $r$ the number of such parts is
$2^q$ with $q=2^r-1$, since $j=1,...,2^r-1$, for $r=\infty $ their
family is infinite and uncountable of the cardinality ${\sf c} =
2^{\aleph _0}$.
\par Then embed each partitioned copy of $({\cal A}_r, ni_1 )$ into
${\cal A}_r^2$ and bend slightly each subset $(\{ z\in {\cal A}_r:
z_0<0 \}, ni_1)$ in directions $\nu _1, \nu _2,...$ perpendicular to
$(i_1,ni_1), (i_2,ni_1),...$ using the imaginary part ${\cal I}_r$
of the second multiple such that after this procedure $\{ (\mbox{
}_1z,\mbox{ }_2z) \in (S(k_1,k_2,... ),ni_1) \cap
(S(l_1,l_2,...),mi_1): \mbox{ }_1z_0<0 \} =\emptyset $ for each
either $n\ne m$ with arbitrary $k, l$ or $n=m$ with $k\ne l$, where
$l=(l_1,l_2,...)$, $z=(\mbox{ }_1z,\mbox{ }_2z)\in {\cal A}_r^2$,
$\mbox{ }_1z$ and $\mbox{ }_2z\in {\cal A}_r$. Then identify faces
$Q_j := P_j\setminus (\bigcup_{m, m\ne j}P_m)$ of two copies $n$ and
$n+1$ of $(S(k),ni_1)$ and $(S(k-2e_j),(n+1)i_1)$ by the
corresponding straight rays of two copies of $(Q_j,ni_1)$ and
$(Q_j,(n+1)i_1)$, where $k=(k_1,k_2,...)$, $k_1, k_2,... \in \{ -1,
1 \} $, $k_j=1$, $e_j=(0,...,0,1,0,...)\in {\bf R}^{2^r-1}$,
$2^{\infty }-1= {\sf c}$. Do this equivalence relation for all $n\in
\bf Z$, each $1\le j\in \bf Z$ and each $k$ with $k_j=1$. Consider
after this identification that $Q_j$ is the part of $({\cal
A}_r,(n+1)i_1)$. Denote by $\cal L$ the $2^r$-dimensional surface in
${\cal A}_r^2$ obtained by such procedure. \par To each
perpendicular transition through the face $Q_j$ from $(S(k),ni_1)$
to $(S(k-2e_j),(n+1)i_1)$ attach the change $2\pi i_j$ of the
argument of the Cayley-Dickson number, where $k_j=1$, $1\le j\in \bf
Z$. To the perpendicular to $Q_j$ transition in the opposite
direction from $(S(k-2e_j),(n+1)i_1)$ to $(S(k),ni_1)$ with $k_j=1$
attach the opposite change of the argument $-2\pi i_j$. \par
Consider the spherical coordinates $(a,\theta _1,...,\theta _m)$ in
the Euclidean space ${\bf R}^{m+1}$ which are related with the
Cartesian coordinates $x_1,...,x_{m+1}\in \bf R$ of a vector
$x=(x_1,...,x_{m+1})\in {\bf R}^{m+1}$ by the equations:
\par $(2)$ $x_1=a\cos (\theta _1)$, $x_2=a\sin (\theta _1)\cos
(\theta _2)$,...,$x_m=a\sin (\theta _1)...\sin (\theta _{m-1})\cos
(\theta _m)$, $x_{m+1}=a \sin (\theta _1)...\sin (\theta _m)$, \\
where $0\le a=|x|<\infty $, $0\le \theta _1\le 2\pi $, $0\le \theta
_2\le \pi $,...,$0\le \theta _m\le \pi $ (see \S XII.1
\cite{zorich}). Then this gives the spherical coordinates in ${\cal
A}_r$ taking $x_{j+1}=z_j$ for each $j=0,1,2,...,2^r-1$ and
$m=2^r-1$, where $z=\sum_{j=0}^{2^r-1}z_ji_j\in {\cal A}_r$.
Comparing Equations $(1)$ and $(2)$ gives:
\par $(3)$ $M=i_1\cos (\theta _2) + i_2\sin (\theta _2)\cos (\theta
_3) +... +i_{2^r-2}\sin (\theta _2)...\sin (\theta _{2^r-2})\cos
(\theta _{2^r-1})$ \\ $ + i_{2^r-1}\sin (\theta _2)...\sin (\theta
_{2^r-1})$ and $\theta _1=\phi $. \\
For ${\cal A}_{\infty }$ the limit of $(2)$ when $r$ tends to the
infinity gives spherical coordinates in ${\cal A}_{\infty }$, since
for each $z\in {\cal A}_{\infty }$ the norm
$|z|:=(\sum_{j=0}^{\infty }z_j^2)^{1/2}<\infty $ is finite.
Therefore, each non-zero $z=|z|\exp (M\phi _1)$ is periodic
(invariant) under substitutions $\theta _j\mapsto \theta _j+2\pi
m_j$ for each $j$, moreover, $z$ is invariant relative to the
pairwise substitutions: $\theta _j\mapsto 2\pi m_j- \theta _j$ and
$\theta _{j+1}\mapsto \theta _{j+1}+(2m_{j+1}+1)\pi $ for each
marked $j$, where $m_1,...,m_{2^r-1}\in \bf Z$.
\par To each spherical coordinates $(\theta _1+2\pi m_1,\theta
_2+\pi m_2,...,\theta _{2^r-1}+\pi m_{2^r-1})=:\psi $ attach two
vectors $m^+=(m_1^+,m_2^+,...)$ and $m^-=(m_1^-,m_2^-,...)$, where
$m_j^+ := \max (0,m_j)$, $m_j^- := \min (0,m_j)$, $|\psi |^2:=
\sum_j\psi _j^2<\infty $, the set $\{ j: m_j\ne 0 \} $ is finite,
since only rectifiable curves in ${\cal A}_r$ are considered. Then
$m=m^+ + m^-$ and put $n^+:=\sum_jm_j^+$, $n^-:= \sum_jm_j^-$, hence
$0\le n^+\in \bf Z$ and $0\ge n^-\in \bf Z$. Therefore, to each such
$\psi $ there corresponds a unique $Arg (z)$ and $z$ is uniquely
characterized by two points $(y_1,y_2)$, where $y_1$ in $({\cal
L},1)$ belongs to $({\cal A}_r,n^+)$ and $y_2$ in $({\cal L},2)$
belongs to $({\cal A}_r,n^-)$ whose spherical coordinates are
$(|z|,\psi ^+)$ and $(|z|,\psi ^-)$ correspondingly, where $\psi
_1^+=\theta _1+2\pi m_1^+$, $\psi _1^-=\theta _1+2\pi m_1^-$, $\psi
_j^+=\theta _j+\pi m_j^+$ and $\psi _j^-=\theta _j^- +\pi m_j^-$ for
each $j\ge 2$, $({\cal L},1)$ and $({\cal L},2)$ are two copies of
$\cal L$. Then embed ${\cal R}:= \{ (y_1,y_2): y_1\in ({\cal L},1)$
$\mbox{with}$ $n_1\ge 0, y_2\in ({\cal L},2)$ $\mbox{with}$ $n_2\le
0 \} $ into ${\cal A}_r^2$, which is possible, since ${\cal
L}\subset {\cal A}_r\times {\cal I}_r$. Points $y_1$ and $y_2$ are
equivalent if and only if $n^+ = n^- =0$. Then ${\cal R}$ is the
noncommutative for $2\le r$ and non-associative for $3\le r$ analog
of the Riemann surface such that $Ln : {\cal A}_r\setminus \{ 0 \}
\to \cal R$ is the univalent mapping and $Ln $ has the countable
number of branches such that $Ln (z)=ln |z| +Arg (z)$, where $ln :
(0,+\infty )\to \bf R$ is the usual real natural logarithm.
\par {\bf 11. Lemma.} {\it If $U$ and $V$ are open domains in
${\cal A}_r$, $2\le r\le 3$, $f: U\to {\cal A}_r$, $\psi : V\to
{\cal A}_r$, $f$ is holomorphic on $U$ and $\psi $ is a holomorphic
diffeomorphism of $V$ on $U$, $\gamma $ is a rectifiable curve in
$U$, where $\gamma (t)=z_0+Exp_n(a_1,...,a_n;t)$ for each $t\in
[0,1]$, $a_1,...,a_{n-1}\in {\cal A}_r\setminus \{ 0 \} $ are
nonzero constants or ${\cal A}_r$ pseudo-conformal functions,
$a_n=const \in {\cal A}_r\setminus \{ 0 \} $, $Re (a_n)=0$,
$|a_n|=2\pi $, $f(\gamma (t))\ne 0$ for each $t\in [0,1]$, then
$\Delta _{\gamma } Arg_n f= \Delta _{\eta }Arg_nf\circ \psi $ and
$\Delta _{\gamma }Arg_nf$ is independent from $a_1,...,a_{n-1}$,
when $n\ge 2$, where $\phi (z) := Ln_n(a_1,...,a_n;\psi ^{-1}(z))$
on $U$ and $\eta (t) := Exp_n(a_1,...,a_n;\phi (\gamma (t))$ for
each $t\in [0,1]$.}
\par {\bf Proof.} Compositions of pseudo-conformal functions
are pseudo-conformal, the inverse of a pseudo-conformal mapping is
pseudo-conformal (see \cite{ludagpd} and Theorem 2.6
\cite{lusmnfgpcd}). Since $\exp $ and $Ln $ are pseudo-conformal,
the mappings $z\mapsto az$ and $z\mapsto za$ are pseudo-conformal
for $a\ne 0$, then $Exp_n$ and $Ln_n$ are also pseudo-conformal for
$a_1\ne 0$,...,$a_n\ne 0$. Choose a branch of the logarithmic
function (see Proposition 10.2) and consider $\phi (z) :=
Ln_n(a_1,...,a_n;\psi ^{-1}(z))$ and put $\eta (t) :=
Exp_n(a_1,...,a_n;\phi (\gamma (t))$, hence $\psi (\zeta )=\gamma
(t)=z$ if and only if $\zeta =\eta (t)$. On the other hand, $\phi
(z)$ is the holomorphic mapping as the composition of holomorphic
mappings. Thus, $\eta $ is the rectifiable curve, since $\gamma $ is
the rectifiable curve. A rectifiable curve $\gamma $ is compact in
${\cal A}_r$, hence it can be covered with a finite number of balls
on each of which $f$ has not zeros, since $f$ is continuous and has
not zeros on $\gamma $. Since $\int_a^bdg(z)=g(b)-g(a)$ for a
holomorphic function on a ball $W$ in ${\cal A}_r$ (see Theorem 2.18
\cite{lusmnfgpcd}), then
\par $\Delta _{\gamma }Arg_nf=\int_{\gamma }dLn_n(a_1,...,a_{n-1},1;
f(z)) = \int_{\eta }dLn_n(a_1,...,a_{n-1},1;f\circ \psi (\zeta ))
=\Delta _{\eta } Arg_nf\circ \psi $. \\
Since $Ln_{n-1}(a_1,...,a_{n-1};z)$ is the inverse function of
$Exp_{n-1}(a_1,...,a_{n-1};y)$, then $\Delta _{\gamma }Arg_nf$ is
independent from $a_1,...,a_{n-1}$, when $n\ge 2$. This is valid for
each phrase $\mu $ representing $f$ and for each branch of the line
integral, for example, specified with the help of the left or the
right algorithm (see Lemma 2.16 and Theorems 2.17, 2.18 in
\cite{lusmnfgpcd} and \cite{lusmfcdv}). Phrases corresponding to $f$
are consistent for canonical (analytic) elements which are analytic
extensions of each other in the domain due to the monodromy Theorem
2.1.5.4 \cite{lusmgdlcm} and 2.45 \cite{lusmnfgpcd}.
\par {\bf 12. Lemma.} {\it Let $f$ be a $(p,r,b)$-quasi-conformal
function on an open connected domain $U$ in ${\cal A}_b$ with a zero
$z_0\in U\cap {\cal A}_r$, $f(z_0)=0$, where $1\le r<b\le 3$,
$0<p\in \bf Z$. Then $f$ has a connected surface
$S=S_{z_0}$ in ${\cal A}_b$ of zeros of $f$ such that $z_0\in S$ and
its dimension over $\bf R$ is $dim_{\bf R}S=2^b-2^r$.}
\par {\bf Proof.} Since $f(z+y_0)={\hat R}_{z,x}f(y_0+x)$ for each
$x\in (U-y_0)\cap {\cal A}_r$ and $z\in U-y_0$ such that $Re(z)=Re
(x)$ and $z={\hat R}_{z,x}x$, where $y_0$ is a marked point in
$U\cap {\cal A}_r$, then $f(z)={\hat R}_{z-y_0,z_0-y_0}f(z_0)=0$ for
each $z\in U$ such that $Re (z)=Re (z_0)$ and $z-y_0={\hat
R}_{z-y_0,z_0-y_0}(z_0-y_0)$, since $f(z)={\hat
R}_{z-y_0,z_0-y_0}f(z_0)$ is obtained from $f(z_0)$ with the help of
the automorphism ${\hat R}_{z-y_0,z_0-y_0}$ of the Cayley-Dickson
algebra ${\cal A}_b$, which is the quaternion skew field $\bf H$ for
$b=2$ or the octonion algebra $\bf O$ for $b=3$ (see Definition 1).
The family of automorphisms ${\hat R}_{z,x}$ is holomorphic and
satisfies Conditions $(Q2-Q5)$ such that when $z$ tends to a point
$\zeta $ in $(U-y_0)\cap {\cal A}_r$, then ${\hat R}_{z,x}$ tends to
the unit operator for a given $x\in (U-y_0)\cap {\cal A}_r$,
consequently, $S_{z_0}:=S:= \{ z: z-y_0={\hat
R}_{z-y_0,z_0-y_0}(z_0-y_0), z\in U, Re (z)=Re (z_0) \} $ is
connected. In particular, ${\hat
R}_{z_0-y_0,z_0-y_0}(z_0-y_0)=z_0-y_0$, since $z_0\in U\cap {\cal
A}_r$, hence $z_0\in S$. Its dimension over $\bf R$ is: $dim_{\bf
R}S =2^b-2^r$, since $dim_{\bf R}{\cal A}_b=2^b$, $dim_{\bf R} {\cal
A}_b\ominus {\cal A}_r=2^b-2^r$ for ${\cal A}_b\ominus {\cal A}_r$
considered as the $\bf R$-linear space.
\par {\bf 13. Notes.} Generally the product of
$(r,b)$-quasi-conformal functions (with a prescribed order of
multiplication for $b=3$) even for $r=1$, where $1\le r<b\le 3$ need
not be $(r,b)$-quasi-conformal, since the derivative operator of the
product is the sum of operators (see Definition 2.2(11)
\cite{luoyst,lusmfcdv,luhcnfcdv}). Indeed, the sum of
pseudo-conformal or quasi-conformal functions may be
non-pseudo-conformal or non-quasi-conformal respectively even when a
derivative operator is non-zero, especially for $r=2$, since there
are projection operators $\pi _j$ from ${\cal A}_b$ into $i_j\bf R$
for each $j=0, 1,...,2^b-1$ and every $2\le b$, where $\pi
_0(z)=(z+(2^b-2)^{-1} \{ -z +\sum_{k=1}^{2^b-1}i_k(zi_k)^*) \} )/2$,
$\pi _j(z)=(-zi_j+i_j(2^b-2)^{-1}\{ -z +\sum_{k=1}^{2^b-1}
i_k(zi_k)^*) \} )/2$ for each $j\in \{ 1,...,2^b-1 \} $, $\pi
_j(z)=z_j$ for each $z\in {\cal A}_b$, where
$z=z_0i_0+...+z_{2^b-1}i_{2^b-1}$. This is the effect of the
noncommutativity of the Cayley-Dickson algebras for $2\le b$.
Moreover, starting from complex constants $a=a_0+{\bf i}a_1$ with
$a_0, a_1\in \bf R$ gives ${\hat R}_{z,y}a=a_0+Ma_1$, where $M\in
{\cal A}_b$ depends on $z\in {\cal A}_b$, $y\in {\cal A}_r$, as it
is described by Formulas $2(4,6)$ in the $z$-representation, hence
\par $(1)$ $(\partial {\hat R}_{z,y}a/\partial z).h=ha_1$ for each
$h\in {\cal I}_b$,
$z\in {\cal A}_b\setminus {\cal A}_r$.
\par By the same reasoning the sum of pseudo-conformal or quasi-conformal
mappings may be a non-pseudo-conformal or non-quasi-conformal
mapping respectively for $2\le b$ even when its derivative is
non-zero, especially for $r=2$. Each complex holomorphic function
$f$ on $\bf C$ (integral function) can be decomposed in accordance
with the Weierstrass Theorem V.72 \cite{lavrshab} as $f(z)=z^m\exp
[g(z)]\prod_{n=1}^{\infty } (1-z/a_n) \exp (z/a_n +
(z/a_n)^2/2+...+(z/a_n)^{p_n}/p_n)$, where $p_n$ is a sequence of
natural numbers and $g(z)$ is an integral function, $m$ is the
multiplicity of $z=0$ as the zero of $f(z)$. But its extension with
the help of automorphisms satisfying Conditions $(Q2-Q5)$ may be
non-quasi-conformal function in view of obstacles described above.
\par When a family ${\hat R}_{z,x}$ is given as ${\hat R}_{\xi
(z),\xi (x)}=R_{\xi (z), \xi (x)}$ (see Example 2), then the phrases
corresponding to canonical (analytic) elements of $f$ which are
analytic extensions of each other in the domain are defined
consistently due to the monodromy Theorem 2.1.5.4 \cite{lusmgdlcm}
and 2.45 \cite{lusmnfgpcd}, where $\xi $ is a pseudo-conformal
diffeomorphism of $U$. This is usually simpler, when $\xi =id$ is
the identity mapping. If $f$ is $(1,b)$-quasi-regular, $2\le b\le
3$, then $q=g|_{W-y_0}$ is the complex holomorphic function and
phrases of analytic elements of $q$ are commutative over $\bf C$.

\par {\bf 14. Remarks and Definitions. 1.}
Zeros and poles of complex holomorphic
functions are defined classically in the standard way.
For an ${\cal A}_r$ $p$-pseudo-conformal function $f$ on an open
domain $V$ we call a point $z_0\in V$ the zero of $f$,
if $f(z_0)=0$, where $2\le r\le 3$, $1\le p\in \bf Z$.
In view of Theorem 2.5 \cite{lusmnfgpcd} and \S 1
its order is $p$.
For an ${\cal A}_r$ $p$-pseudo-conformal function $f$
on $V\setminus \{ z_0 \} $, where $2\le r\le 3$, $V$ is open in
${\cal A}_r$ and $z_0\in {\cal A}_r$ is a point, then we call
$z_0$ the pole of $f$, if $g(z) := 1/f(z)$ is ${\cal A}_r$
$p$-pseudo-conformal on $V\setminus \{ z_0 \} $ and $z_0$ is the
zero of $g$.
\par Define $f$ to be holomorphic in a neighborhood $V$ of $\infty $
or $(p,r)$-pseudo-conformal or $(p,r,b)$-quasi-conformal at $\infty $
if and only if $g(z):=f(1/z)$ is holomorphic in $U:=\{ z: 1/z\in
V \} $ or $(p,r)$-pseudo-conformal or
$(p,r,b)$-quasi-conformal at zero respectively.
We say that $y=\infty $ is a zero or a pole of $f$
if and only if $g(z)=f(1/z)$ has a zero or pole at $z=0$ respectively.

\par {\bf 2.} Consider the following situation having
the natural embedding of ${\cal A}_r$ into ${\cal A}_b$ associated
with the standard doubling procedure, where $1\le r<b\le 3$. Suppose
that $U$ is an open connected domain in ${\cal A}_b$ and $W=U\cap
{\cal A}_r$ is an open connected domain in ${\cal A}_r$ such that
$U$ is pseudo-conformally diffeomorphic with a domain $V$, where $V$
is obtained from $W$ with the help of all rotations in all planes
${\bf R}i_v\oplus {\bf R}i_u$ with $v=1,...,2^r-1$ and
$u=2^r,...,2^b-1$ on angles $\phi \in (0,2\pi )$ with the real
rotation axis, since each operator $T\in SO(2^b,{\bf R})$ can be
presented as the finite product of one-parameter subgroups and here
is considered its subgroup $SO_{\bf R}(2^b,{\bf R})$ of operators
restrictions of which on the real axis $\bf R$ is the identity.

Let $f$ be a $(p,r,b)$-quasi-conformal or $(p,r,b)$-quasi-regular
mapping on $U$ may be besides a finite number of surfaces $S_y := \{
y_0+z: z={\hat R}_{z,y-y_0}(y-y_0), z\in U-y_0, Re (z)=Re (y-y_0) \}
$ of poles $y\in W$, which may only be points in $W$, where $0<p\in
\bf Z$.

\par Put for each $z_0\in S_y$ by the definition: \par $(1)$ $(2\pi
)^{-1}\int_{\gamma }f(z)dz =: res(z_0,f).N$ the residue operator of
$f$ at $z_0$, where $N\in {\cal A}_b$, $Re (N)=0$, $|N|=1$, $\gamma
(t)= z_0+\rho \exp (2\pi tN)\subset V$, $t\in [0,1]$, $\rho >0$ is
sufficiently small such that $f|_{B\setminus \{ z_0 \} }$ is locally
analytic and $\gamma $ does not encompass another poles of $f$ in
the set $ \{ y_0+z: z\in {\hat R}_{q,y-y_0}{\cal A}_r; |z-y|<\rho
+\epsilon \} $ for some $0<\epsilon <\infty $ and some $q\in {\cal
A}_b$ such that $N\in {\hat R}_{q,y-y_0}{\cal A}_r$ with $Re (q)=Re
(y-y_0)$, $B=B({\bf R}\oplus N{\bf R},z_0,2\rho ^-)$, $B(X,z,R^-) :=
\{ x\in X: d_X(x,z)<R \} $ denotes the open ball in a metric space
$X$ with a metric $d_X$. For $a\in \bf R$ put $res(z_0,f).(aN) := a~
res(z_0,f).N$ and $res(z_0,f).0=0$.

\par Suppose also that $z_0\in W$ is a zero or a pole of $f$ and
$S=S_{z_0}$ be the surface corresponding to $z_0$ from Lemma 12,
where $W$ may contain only finite number of zeros or poles $z_l$,
which may only be points. For a subset $G$ in ${\cal A}_r$, let $\pi
_{s,q,t}(G):= \{ u: z\in G, z=\sum_{v\in \bf b}w_vv,$ $u=w_ss+w_qq,
w_v\in {\bf R} \forall v\in {\bf b} \} $ for each $s\ne q\in \bf b$,
where $t:=\sum_{v\in {\bf b}\setminus \{ s, q \} } w_vv \in {\cal
A}_{r,s,q}:= \{ z\in {\cal A}_r:$ $z=\sum_{v\in \bf b} w_vv,$
$w_s=w_q=0 ,$ $w_v\in \bf R$ $\forall v\in {\bf b} \} $. That is,
geometrically $\pi _{s,q,t}(G)$ is the projection on the complex
plane ${\bf C}_{s,q}$ of the intersection of $G$ with the plane
${\tilde \pi }_{s,q,t}\ni t$, ${\bf C}_{s,q}:=\{ as+bq:$ $a, b \in
{\bf R} \} $, since $sq^*\in {\hat b}$, where ${\bf b}:= \{
i_0,i_1,...,i_{2^r-1} \} $ is the family of the standard generators
of the Cayley-Dickson algebra ${\cal A}_r$, ${\hat b}:= {\bf
b}\setminus \{ i_0 \} $, $i_0=1$.

\par Suppose that $\omega $
is a rectifiable loop, that is, a closed curve, $\omega (0)=\omega
(1)$, in an open sub-domain $J$, $z_0\in J\subset W$ in ${\cal A}_r$
such that $\omega $ encompasses $z_0$, where $J$ does not contain
any other zero or pole of $f$, $J$ is $(2^r-1)$-connected and $\pi
_{s,q,t}(J)$ is simply connected in $\bf C$ for each $t\in {\cal
A}_{r,s,q}$ and $u\in {\bf C}_{s,q}$, $s=i_{2k}$ and $q=i_{2k+1}$,
$k=0,1,...,2^{r-1}-1$ for which there exists $\zeta =u+t\in J$.

\par {\bf 15. Theorem.} {\it  Let $f$ be a pseudo-conformal function on
$V\setminus \{ y \} $ with a pole at $y$ in $V$ and let also $F$ be
a univalent branch of its $(r,b)$-quasi-conformal extension in
$W\setminus S_y$ relative to a marked point $y_0$, where $W$ is an
open subset in ${\cal A}_b$ such that $W\cap {\cal A}_r=V$. Then the
residue operators $res(y,f)$ and $res(z,F)$ are such that
$res(z,F).M={\hat R}_{z-y_0,y-y_0}[res(y,f).N]$ for each $z\in
S_y\cap W$ and every $N\in {\cal A}_r$ with $Re (N)=0$, where
$M={\hat R}_{z-y_0,y-y_0}N$. Moreover, $res(z_0,f).N$ is $\bf R$
homogeneous and ${\cal A}_b$ additive by $f$.}
\par {\bf Proof.}  If $y$ is a finite point,
then $z\in S_y$ is a finite point and \par $(1)$ $(2\pi
)^{-1}\int_{\gamma }f(z)dz = res(y,f).N$, where $\gamma (t)=y+\rho
\exp (2\pi tN)$ and $\rho >0$ is sufficiently small such that
$\gamma $ does not encompass another poles of $f$ in the set $ \{
y_0+z: z\in {\hat R}_{q,y-y_0}{\cal A}_r; |z-y|<\rho +\epsilon \} $
for some $0<\epsilon <\infty $ and some $q\in {\cal A}_b$ with $Re
(q)=Re (y-y_0)$ such that $N\in {\hat R}_{q,y-y_0}{\cal A}_r$. Using
Conditions $(Q6,Q7)$ the action of ${\hat R}_{z-y_0,y-y_0}$ on both
sides of Equation $(1)$ gives
\par $(2)$ ${\hat R}_{z-y_0,y-y_0}[res (y,f).N]=(2\pi )^{-1}
\int_{\eta }F(s)d s$ for $Re (z)=Re (y)$,\\ where $\eta (t)= {\hat
R}_{z-y_0,y-y_0}\gamma (t)$ for each $t$.

\par Since each $z\in S_y$ is the pole of $F$ restricted on the
corresponding subalgebra ${\hat R}_{z-y_0,y-y_0}({\cal A}_r)$ in
${\cal A}_b$, then there is defined the ${\cal A}_b$-additive and
$\bf R$-homogeneous operator $res(z,F).M=(2\pi )^{-1}\int_{\eta
}F(s)d s$ for $M\in {\cal A}_b$, $Re (M)=0$ (see Theorem 3.23 \cite{
luoyst,lusmfcdv,luhcnfcdv}).

Therefore, the first statement of this theorem follows from Equation
$(2)$ and Conditions $(Q6,Q7)$, \S 14 above and Theorem 2.5
\cite{lusmnfgpcd}.
\par Since $\int_{\gamma }(a_1f_1(z)+a_2f_2(z))dz=
a_1\int_{\gamma }f_1(z)+a_2\int_{\gamma }f_2(z))dz$ for each $a_1,
a_2\in \bf R$ and ${\cal A}_b$ holomorphic functions $f_1$ and $f_2$
on the domain $U$ containing a rectifiable curve $\gamma $ (see
Theorem 2.7  \cite{luoyst,lusmfcdv,luhcnfcdv}), then $res(z_0,f).N$
is $\bf R$-homogeneous and ${\cal A}_b$-additive by $f$:
$res(z_0,a_1f_1+a_2f_2)=a_1 res (z_0,f_1)+ a_2 res (z_0,f_2)$.

\par If $y=\infty $, then consider $g(z)=f(1/z)$
and $g$ has the pole at $z=0$, hence in this case the statement of
this theorem follows from the the first part of the proof.
\par {\bf 15.1. Example.} If a function $f$ can be written in
the form $f(z)=(a(z)((b(z)1/(z-y))c(z)))e(z)$ in a neighborhood of
$y\in {\cal A}_r$, where $a(z)$, $b(z)$, $c(z)$ and $e(z)$ are
${\cal A}_r$-holomorphic and $a(y)\ne 0$, $b(y)\ne 0$, $c(y)\ne 0$
and $e(y)\ne 0$, $2\le r\le 3$. Then $res (y,f).N = (2\pi )^{-1}
\lim_{0<\rho \to 0} \int_{\gamma }(a(z)((b(z)(1/(z-y))c(z)))e(z) dz
= (a(y)((b(y)N)c(y))e(y)$. At the same time for $\gamma $ from
Definition 10 with $a_n=2\pi M$, $M\in {\cal A}_b$, $Re (M)=0$,
$|M|=1$, $z_0=0$, $\xi (t)=t$ for each $t$, we have $\Delta _{\gamma
} Arg_n \gamma =2\pi M$.

\par {\bf 16. Theorem.} {\it Let $U$ be a proper open subset
in ${\cal A}_b$, let also $f_1$ and $f_2$ be two continuous
functions from the closure ${\bar U} := cl (U)$ of $U$ into ${\hat
{\cal A}}_b$ such that on a topological boundary $Fr (U)$ of $U$
they satisfy the inequalities $|f_1(z)|<|f_2(z)|<\infty $ for each
$z\in Fr (U)$, where ${\hat {\cal A}}_b := {\cal A}_b\cup \{ \infty
\} $ is the one-point (Alexandroff) compactification of ${\cal
A}_b$. Suppose $q_2:=f_2$ and $q_1:=f_1+f_2$ are
$(p,r,b)$-quasi-meromorphic functions in $U$ and zeros and poles of
$f_j|_W$ are isolated, where $W=U\cap {\cal A}_r$, $1\le r<b\le 3$,
$j=1, 2$. Let also $q_j$ be $(p,r,b)$-quasi-conformal in a
neighborhood $U_{z_0}$ in $U$ of each its zero $z_0$ and $1/q_j(z)$
be $(p,r,b)$-quasi-conformal in $U_{z_0}\setminus \{ z_0 \} $ for
each pole $z_0$, where $p\in \bf N$ may depend on $z_0$, for $j=1$
and $j=2$. Suppose also that $\gamma $ from Definition 10 is a loop,
where $\gamma (0)=\gamma (1)$ is a loop, does not cross any
$S_y(q_j)$ for any zero or pole $y$ of $q_j$ for $j=1$ and $j=2$,
where $\gamma \subset Fr (U)$, $1\le n$. Then $\Delta _{\gamma }
Arg_nq_1=\Delta _{\gamma }Arg_nq_2$.}
\par {\bf Proof.} Put without loss of generality $z_0=0$ and $\rho =1$
for those of Definition 10. If $n>1$ consider $h_j:=q_j\circ
Exp_{n-1}(a_1,...,a_{n-1};z)$ instead of $q_j$ for $j=1$ and $j=2$,
since the compositions of mappings are associative in the set
theoretic sense. On the other hand, $Exp_{n-1}(a_1,...,a_{n-1};z)$
is the pseudo-conformal mapping for $a_1\ne 0$,...,$a_{n-1}\ne 0$.
In view of Theorem 8.2 $h_j$ satisfy suppositions of this theorem.
Substituting $q_j$ on $h_j$ we can reduce the proof to the $n=1$
case, since $Ln_{n-1}(a_1,...,a_{n-1};\gamma (t))=\xi (t)$ for a
branch of $Ln$ such that $Ln (1)=0$. For example, it is possible to
take $a_n=2\pi M$ and $\xi (t)=t$, $t\in [0,1]$, $M\in {\cal
A}_b\ominus {\cal A}_r$, $|M|=1$, $Re (M)=0$ in Definition 10.
Therefore, consider $q_j$ for $n=1$, where $j=1, 2$.
\par The curve $\gamma $ is rectifiable, hence compact.
Zeros and poles of $f_j|_W$ are isolated consequently, there exists
a sequence $\{ \psi _m: m\in {\bf N} \} $ of rectifiable loops in
$U$ converging to $\gamma $ uniformly when $m$ tends to the infinity
and such that each $\psi _m$ does not cross any $S_y(q_j)$ for $j=1$
and $j=2$ for zero or pole $y$ of $f_j$. Thus consider the integral
$\int_{\gamma }d Ln (q_j(z))$ along $\gamma $ as the limit of
$\int_{\psi _m}d Ln (q_j(z))$ when $m$ tends to the infinity, since
$f_j$ is continuous in a neighborhood $V$ of $Fr (U)$ in $\bar U$
and $\psi _m\subset U$ for each $m$ and $j$, where $V$ does not
contain any zero or pole of $q_1$ and $q_2$. The latter $V$ exists,
since $|f_1(z)|< |f_2(z)|<\infty $ for each $z\in Fr (U)$.
\par If $z_0$ is a pole of $q_j$ at $z_0$, then $1/q_j(z)$ has a zero at
$z_0$.  There are not any zeros or poles of $q_1$ and $q_2$ on $Fr (
U)$, since $|f_1(z)|<|f_2(z)|<\infty $ on $Fr (U)$.
\par  For each chosen branch of the logarithmic function there is the
equality $Ln (1/q_j)=-Ln (q_j)$ (see Proposition 10.2). Moreover,
$q_1=f_1+f_2=f_2+f_2[(1/f_2)f_1]= f_2(1+(1/f_2)f_1)$, since ${\cal
A}_3=\bf O$ is alternative, and $|(1/f_2)f_1|<1$ on $Fr (U)$.
Consider the triangle formed by the vectors $q_1(z)$, $q_2(z)$ and
$q_1(z)-q_2(z)=f_1(z)$ for $z\in Fr (U)$, then $Arg~ q_1(z)= Arg~
q_2(z) + \phi (z)$ such that $|\phi (z)|<\pi /2$ for each $z\in Fr
(U)$ for a chosen branch of $Ln$, where $Arg~q_j$ and $\phi (z)\in
{\cal I}_b$. Therefore, $\Delta _{\xi } Arg~ q_1 = \Delta _{\xi }
Arg~ q_2+\Delta _{\xi }\phi = \Delta _{\xi }Arg~ q_2$, since $\xi
(0)=\xi (1)$ and $|\phi (z)|<\pi /2$ for each $z\in \xi
([0,1])\subset Fr(U)$ such that $\Delta _{\xi }\phi =0$. Indeed, the
point $w(z)=(1/f_2(z))f_1(z)$ is within the unit open ball $ B({\cal
A}_b,0,1^-) := \{ w \in {\cal A}_b: |w|<1 \} $. Therefore, the
vector $v=1+w$ can not rotate on $2\pi $ around zero. Thus the
winding numbers of $q_1$ and $q_2$ around the zero are the same.
From the relations of $h_j$ with $q_j$ for $n>1$ the statement of
this theorem follows for $n>1$ as well: $\Delta _{\gamma }Arg_nq_1=
\Delta _{\gamma }Arg_nq_2$.
\par {\bf 16.1. Theorem.} {\it  Let suppositions of Remark 14.2
be satisfied, when \\ $W=B({\cal A}_r,y_0,R^-)\setminus A$, where $A
:= \{ y\in B({\cal A}_r,y_0,R^-): f'(y)=0 \} $ consists of isolated
points, $0<R<\infty $. Suppose also that $W$ contains either zeros
or poles of a $(q,r,b)$-quasi-regular function $f$, but not zeros
and poles simultaneously, $1\le r<b\le 3$. Then for each rectifiable
curve $\omega $ in $J$ encompassing $z_0$ and each $2^r\le n\le
2^b-1$ there exists a family of rectifiable curves $\gamma $ in $U$
encompassing $S=S_{z_0}$ such that $\gamma \cap S_{z_l}=\emptyset $
for each zero or pole $z_l$ of $f$ in $W$ and such that $\gamma $ is
not contained in ${\cal A}_r$ and
\par $\Delta _{\gamma }Arg_nf= p K \Delta _{\omega }Arg_1f$ \\
for some $K\in {\cal A}_b$, $|K|=1$, $Re (K)=0$, $K=K(\gamma )$,
$1\le p\in \bf Q$.}
\par {\bf Proof.} The zero or pole $z_0$ of $f$ is isolated in $W$,
hence $f(z)\ne 0$ in $Y\setminus \{ z_0 \} $ for a sufficiently
small neighborhood $Y$ of $z_0$ in ${\cal A}_r$. If $z_0$ is a pole
of $f$, then $z_0$ is a zero of $1/f$ and vice versa. If $z\in
Y\setminus \{ z_0 \} $, then $f(z)\ne 0$ and apply automorphisms
${\hat R}_{\zeta -y_0,z-y_0}$ to $f(z)$ by all $\zeta \in U$ with
$Re (\zeta )=Re (z)$ and $(\zeta -y_0)={\hat R}_{\zeta
-y_0,z-y_0}(z-y_0)$. Take without loss of generality $y_0=0$, since
$z\mapsto z+y_0$ is the bijective pseudo-conformal mapping from
${\cal A}_b$ on ${\cal A}_b$. This gives the closed surface
$S_{f(z)}$ analogous to $S=S_{z_0}$. In accordance with Lemma 12
$dim_{\bf R}S_{z_0}= dim_{\bf R}S_{f(z)}=2^b-2^r$.

\par If $z_k\in {\cal A}_r$ is obtained from
$z_l\in {\cal A}_r$ with $l\ne k$ due to a rotation around the real
axis in a plane $\pi ^{k,l}$ contained in ${\cal A}_b \ominus {\cal
A}_r$, which corresponds to the one-parameter over $\bf R$ subgroup
of rotations in $SO_{\bf R}(2^b,{\bf R})$, then $z_k$ and $z_l$ are
both either zeros or poles due to the conditions of this theorem,
since $f$ is $(q,r,b)$-quasi-regular.
\par There are the following decompositions of algebras
as the $\bf R$-linear spaces due to the doubling procedure: ${\bf
H}={\bf C}\oplus i_2{\bf C}$, ${\bf O}={\bf H}\oplus i_4{\bf H}$ and
${\bf O}= {\bf C}\oplus i_2{\bf C}\oplus i_4{\bf C} \oplus i_6{\bf
C}$ corresponding to $(r,b)$ pairs equal to $(1,2)$, $(2,3)$ and
$(1,3)$ respectively, that gives the embedding of geometry in ${\cal
A}_r$ into geometry in ${\cal A}_b$. Consider an intersection of the
surface $S$ with the plane $\pi $ containing $z_0$ and perpendicular
to the real axis $\bf R$, $\pi =z_0 + i_s{\bf R}\oplus i_q{\bf R}$,
where $2^r\le s<q\le 2^b-1$. Then $\eta := S\cap \pi $ is a
rectifiable loop containing $z_0$ and $\eta $ has the winding number
1 for each internal point in the domain $P_{\eta }$ encompassed by
$\eta $ in $\pi $ with the boundary $\partial P_{\eta }=\eta $. \par
Consider a rectifiable loop $\gamma $ consisting of the following
parts: the loop $\gamma _+$ outside $P_{\eta }$, the loop $\gamma
_-$ inside $P_{\eta }$, $\psi $, where $\gamma (t)=\gamma _+(3t)$
for $0\le t\le 1/3$, $\gamma (t)=\psi (6t-2)$ for $1/3<t<1/2$,
$\gamma (t)=\gamma _-(3t-3/2)$ for $1/2\le t\le 5/6$, $\gamma
(t)=\psi (6-6t)$ for $5/6 <t\le 1$ such that $\psi $ joins $\gamma
_+$ with $\gamma _-$ such that $\psi $ is gone along twice in one
and the opposite direction, $\gamma _+$ and $\gamma _-$ are in $\pi
$ for which $|\gamma _+(t)|> |\eta (t)|$ and $|\gamma _-(1-t)|<|\eta
(t)|$ and $|\gamma _+(t)-\eta (t)|<\delta $ and $|\gamma
_-(1-t)-\eta (t)|<\delta $ for each $t\in [0,1]$ and $\gamma (t)\in
U$ and $\gamma (t)$ is not zero or pole of the function $f$ for each
$t\in [0,1]$, $\delta >0$ is a sufficiently small constant such that
$\gamma _-$ and $\eta $ encompass the same zeros and poles besides
those belonging to $S_{z_0}$, $\gamma _+$ and $\gamma _-$ have
opposite orientations (see also Theorem 16 and Equations $(1,2)$
below for more details). Since the set $A$ in $W$ consists of
isolated points, then the loop $\gamma $ can be chosen such that
$\gamma ([0,1])\cap (\bigcup \{ S_y: y\in A \} )=\emptyset $. This
encompassment is subordinated in ${\cal A}_b$ to properties of $Ln$
(see Theorems 2.23 and 2.24 \cite{lusmnfgpcd}). Take $0<\rho _+
-\rho _-$ sufficiently small and use the approximation $f(z_0+h) =
f(z_0) +f'(z_0).h + O(h^2)$ and Properties $1.(P1-P3)$ in a
neighborhood of a zero $z_0$ of $f$ in ${\cal A}_r$ and Properties
$1.(Q1-Q7)$ in a neighborhood of $z_0$ in ${\cal A}_b$, where
$f(z_0)=0$ for a zero $z_0$ of $f$ (see also Theorems 2.4 and 2.5
\cite{lusmnfgpcd}).
\par Choose $\psi $ in a plane $\pi _1$ containing $\bf R$
and a point $\zeta '\in {\cal A}_b\ominus {\cal A}_r$ such that
$\psi $ does not intersect any $S_{z_l}$. Hence $\gamma $ does not
intersect any $S_{z_l}$ and encompasses $S=S_{z_0}$. The direction
of $\gamma $ is natural such that in the plane $\pi $ the loop
$\gamma _+$ is gone counter-clock-wise and $\gamma _-$ is gone
clock-wise as seen from the positive axis of $M_{\pi }{\bf R}_+$,
where ${\bf R}_+ := (0,\infty )$, $M_{\pi }\in {\cal I}_b$, $M_{\pi
}\perp \pi $, $M_{\pi }$ corresponds to a vector being the vector
product in the real shadow of basis vectors of the plane $\pi $.
Though, instead of description of an orientation it is sufficient to
write analytic formulas for curves, that is done below.
\par  For $r=1$ if
${\bar z}_0$ is the zero or pole together with $z_0$, then $\gamma $
encompasses $z_0$ and ${\bar z}_0$ symmetrically, since ${\bar z}_0$
is obtained from $z_0$ by rotation on the angle $\pi $ around the
real axis. For $r=2$ if $z_j=\sum_{k=0}^3z_{j,k}i_k$ and $z_l\in
S_{z_0}$, then $z_j$ and $z_l$ are both either zeros or poles due to
the condition imposed in Remark 14.2 and in this theorem, for
example, when $z_{j,k}=-z_{l,k}$ for some $1\le k\le 3$, where
$z_{j,k}\in \bf R$ for each $j, k$.

\par Using the iterated exponent choose $\gamma (t)$ up to
an ${\cal A}_b$-pseudo-conformal diffeomorphism of $U$ in the form
\par $(1)$ $\gamma _+(t)=Re (z_0)+ \rho _+Exp_{n-1}(a_1,...,a_{n-1};
\xi (t))$,
\par $(2)$ $\gamma _-(t)=Re (z_0)+ \rho _-
Exp_{n-1}(a_1,...,a_{n-1};\xi (1-t))$, \\
where $0<\rho _- <\rho _+ <R$, $1\le n\le 2^b-2^r$,
$a_1,...,a_{n-1}\in {\cal A}_b$ are nonzero constants, $\xi ([0,1])$
is a loop, for example, $\xi (t)=\exp (2\pi Mt)$, $t\in [0,1]$,
$M\in {\cal A}_b\ominus {\cal A}_r$, $|M|=1$, $Re (M)=0$. Since $S$
is the smooth $C^{\infty }$ compact manifold having the $C^{\infty
}$ real shadow, then it has $2^b-2^r$ local coordinates. Moreover,
$S$ is homeomorphic with the rotation surface such that $S$ can be
parameterized with angles $\theta _1,...,\theta _m$, where $0\le
\theta _1\le 2\pi $ and $0\le \theta _j\le \pi $ for $j=2,...,m$,
$m=2^b-2^r$, such that $z\in S$ is the function $z=z(\theta
_1,...,\theta _m)$ (see Formulas 10.2(2)).
\par In view of Corollary 3.5 \cite{lusmfcdv,luhcnfcdv} the sphere
$S({\cal A}_b,y_0,R)$ in ${\cal A}_b$ of radius $0<R<\infty $ with
the center at $y_0$ can be parameterized with the help of the
iterated exponential functions. Let $ \{ i_0, i_1, i_2, i_3 \} $ be
the standard generators of the quaternion algebra $\bf H$, where
$i_0=1$, $i_1^2=i_2^2=i_3^2=-1$, $i_1i_2= -i_2i_1=i_3$,
$i_2i_3=-i_3i_2=i_1$, $i_3i_1=-i_1i_3=i_2$, then
\par $(3)\quad \exp (i_1(p_1t+\zeta _1)\exp (-i_3(p_2t
+\zeta _2) \exp (-i_1(p_3t+\zeta _3)))) = \exp (i_1(p_1t+\zeta
_1)\exp (- (p_2t+\zeta _2)(i_3\cos (p_3t+\zeta _3) - i_2\sin
(p_3t+\zeta _3))))$
\par $= \exp (i_1(p_1t+\zeta _1)(\cos (p_2t+\zeta _2) -
\sin (p_2t+\zeta _2)(i_3\cos (p_3t+\zeta _3) - i_2\sin (p_3t+\zeta
_3))))$
\par $= \exp ((p_1t+\zeta _1)(i_1\cos (p_2t+\zeta _2) + i_2
\sin (p_2t+\zeta _2)\cos (p_3t+\zeta _3) + i_3\sin (p_2t+\zeta
_2)\sin (p_3t+\zeta _3)))= \cos (p_1t+\zeta _1)+ i_1\sin (p_1t+\zeta
_1)\cos (p_2t+\zeta _2) + i_2 \sin (p_1t+\zeta _1)\sin (p_2t+\zeta
_2)\cos (p_3t+\zeta _3) + i_3\sin (p_1t+\zeta _1)\sin (p_2t+\zeta
_2)\sin (p_3t+\zeta _3)$, \\
where $p_j, \zeta _j\in \bf R$ for each $j$.
\par Further by induction the equality is accomplished:
\par $(4)\quad \exp (\mbox{ }_{q+1}M(p,t;\zeta ))=$ \\
$\exp \{\mbox{ }_qM((i_1p_1+_...+i_{2^q-1}p_{2^q-1}),t; (i_1\zeta
_1+...+i_{2^q-1} \zeta _{2^q-1})\exp (-i_{(2^{q+1}-1)}($ \\
$p_{2^q}t +\zeta _{2^q}) \exp (-\mbox{
}_qM((i_1p_{2^q+1}+...+i_{2^q-1}p_{2^{q+1}-1}),t; (i_1\zeta
_{2^q+1}+...+i_{2^q-1} \zeta _{2^{q+1}-1}))) \} $, \\
where $i_{2^q}$ is the generator of the doubling of the algebra
${\cal A}_{q+1}$ from the algebra ${\cal A}_q$, such that
$i_ji_{2q}=i_{2^q+j}$ for each $j=0,...,2^q-1$, the function
$M(p,t;\zeta )$ is written with the lower index $\mbox{ }_qM$
and it is given by the equation
\par $(5)$ $M(p,t)=M(p,t;\zeta ) = (p_1t+\zeta _1)[ i_1 \cos
(p_2t +\zeta _2) + i_2 \sin (p_2t+\zeta _2)$ \\
$\cos (p_3t+\zeta _3) +...+ i_{2^q-2} \sin (p_2t+\zeta _2) ...\sin
(p_{2^q-2}t+\zeta _{2^q-2}) \cos (p_{2^q-1}t+\zeta _{2^q-1})$ \\ $+
i_{2^q-1}\sin (p_2t+\zeta _2)...\sin (p_{2^q-2}t+\zeta _{2^q-2})
\sin (p_{2^q-1}t+ \zeta _{2^q-1})]$ \\ for the Cayley-Dickson
algebra with $2\le q<\infty $, where  $\zeta = \zeta _1i_1+...+\zeta
_{2^q-1}i_{2^q-1}\in {\cal A}_q$ is the parameter of an initial
phase, $\zeta _j\in \bf R$ for each $j=0,1,...,2^q-1$. When $t=1$
and $p_j$ are variables $p_1\in [0,2\pi ]$ and $p_j\in [0,\pi )$ for
each $j=2,...,2^q-1$, then the image of the iterated exponent given
by Equation $(4)$ for $2\le q\le 3$ or, in particular, by Formula
$(3)$ is the unit sphere in ${\cal A}_q$, where $\zeta _j$ is fixed
for each $j=1,...,2^q-1$ and can also be taken particularly zero.
This gives the one-sheeted covering of the unit sphere in ${\cal
A}_q$. If $p_j$ and $\zeta _j$, $j=1,...,2^q-1$ are fixed and $t$ is
the variable, then Formulas $(3)$ or $(4,5)$ give the curve in
${\cal A}_q$. It reduces to the loop, when $p_1=2\pi $ and $p_j=0$
or $p_j=\pi $ for each $j=2,...,2^q-1$, $t\in [0,1)$. Particularly,
if $\zeta _j=0$ for each $j>n$, $\zeta _1\ne 0$,...,$\zeta _n\ne 0$
and $p_k=0$ for each $k\ne n$, then the iterated exponent in
Formulas $(3)$ or $(4,5)$ reduces to $Exp_n$.
\par Then $S$ is diffeomorphic with the intersection
$S({\cal A}_b,y_0,R)\cap (i_{2^r}{\bf R}\oplus i_{2^r+1}{\bf
R}\oplus ... \oplus i_{2^b-1}{\bf R})$.
\par In accordance with the Riemann mapping Theorem 4.12.40 over $\bf C$
\cite{shabat} or Theorems 2.1.5.7 and 2.47 \cite{lusmnfgpcd,lusmgdlcm}
over $\bf H$ and $\bf O$ if $P$ is an open
subset in ${\cal A}_q$, $q=1$ or $q=2$ or $q=3$, satisfying conditions of
Remark 14 and with a boundary $\partial P$ consisting more,
than one point, then $P$ is pseudo-conformally equivalent with the
open unit ball in ${\cal A}_q$.
\par Within each sub-domain $P$ of $U$ in ${\cal A}_b$ satisfying
conditions of Remark 14 is applicable the homotopy Theorem 2.11
\cite{luoyst,lusmfcdv,luhcnfcdv} for the line integral over ${\cal
A}_b$. In view of Lemma 11 above we can consider a domain $U$ and
hence a curve $\gamma $ in it up to a pseudo-conformal
diffeomorphism. Therefore, the rest of the proof up to a
pseudo-conformal diffeomorphism is with balls and spheres due to
Conditions $(P1-P3)$.
\par Therefore, there exist $n$, $a_1,...,a_{n-1}$ and $\xi $ such that
$\gamma _+$ and $\gamma _-$ are given by Formulas $(1,2)$ and
$\Delta Arg_1\xi \ne 0$, consequently, $\Delta _{\gamma }Arg_n\gamma
\ne 0$, for example, $a_n=2\pi M$, $M\in {\cal A}_b\ominus {\cal
A}_r$, $|M|=1$, $Re (M)=0$. This is applicable both to $S_{z_0}$ and
$S_{f(z)}$ obtained from $z_0$ and $f(z)$ by families of
automorphisms. If we take $\chi =\pi \cap S_z$ with $Re (z)=Re
(z_0)$, $|Im (z)|>|Im (z_0)|$, $z\ne z_0$ such that $S_z$ does not
contain any pole or zero of $f$, then each $\zeta \in (\pi
-z_0+f(z))\cap S_{f(z)}$ has the form $\zeta =f(\chi (t))$ and the
mapping $[0,2\pi )\ni t\mapsto \zeta \in (\pi -z_0+f(z))\cap
S_{f(z)}$ is bijective, hence $\Delta _{\chi }Arg_nf \ne 0$. This is
possible, since the set $A$ of poles or zeros of $f$ in $W$ consists
of isolated points, $\forall z_0\in A$: $\min \{ |y-z_0|: y\in
A\setminus \{ z_0 \} \} >0 $. In view of Condition $(Q7)$ for
$a_n=2\pi M$ and $\xi (t) = \exp (2\pi Mt)$ in Formulas $(1,2)$ we
have $\Delta _{\chi }Arg_nf=2\pi u k M$, where $|M|=1$, $Re (M)=0$,
$M\in {\cal A}_b$, $k$ is the winding number of $\chi $, which we
take equal to $1$, $u$ is the sum of orders of all either poles or
zeros from $W$ encompassed by $S_z$ (see also Theorem 16). For
sufficiently small $|Im (z)|-|Im (z_0)| = \epsilon >0$ the number
$u$ is equal to the sum of orders of all either poles or zeros
belonging to $S=S_{z_0}$ and to $B({\cal A}_r,Re (z_0),|Im (z_0|) :=
\{ y\in {\cal A}_r: |y-Re (z_0)|\le |Im (z_0)| \} $. \par Then $u\ne
0$, since all $z_l$ belonging to $S$ are simultaneously zeros or
poles together with $z_0$ (see above). It is possible to take
$\gamma =\chi $, that gives $u=u_{\chi }$. If take $\gamma $
consisting of $\gamma _+$ and $\gamma _-$ and $\psi $ as above, then
$u=u_{\gamma _+} - u_{\gamma _-} >0$, since $S_{z_0}\subset B({\cal
A}_b, Re (z_0),\rho _+)\setminus B({\cal A}_b,Re(z_0),\rho _-)$. On
the other hand, $\Delta _{\omega }Arg_1f=2\pi v N$, where $|N|=1$,
$Re (N)=0$, $N\in {\cal A}_r$, $v$ is the order of $z_0$, $v\ge 1$
for zero, $v\le -1$ for pole. Thus $1\le p=u/v\in \bf Q$, also
$M=KN$ for $K=MN^*$ due to the alternativity of ${\cal A}_r$ for
$2\le r\le 3$, since $u$ and $v$ are of the same sign and $|u|\ge
|v|$, $|M|^2=MM^*=M^*M=-M^2$ for purely imaginary $M\in {\cal A}_b$,
and inevitably the statement of this theorem follows.
\par {\bf 16.2. Remark.} If suppositions of Theorem 16.1 are satisfied
besides that $S_{z_0}$ contains both zeros and poles then it may
happen, that $p=0$ due to $N-P=0$, where $N=\sum_kN_k$ is a number
of zeros and $P=\sum_kP_k$ is a number of poles belonging to $S\cap
{\cal A}_r$, where each zero and pole is counted in accordance with
its order $N_k$ and $P_k$ respectively.

\par {\bf 17. Theorem.} {\it Let $f$ be a $(1,b)$-quasi-integral
function such that $f({\tilde z}) = {\tilde f}(z)$ for each $z\in
{\cal A}_b$, also $f^s(p)=f^s(-p)$ for each $p\in {\cal A}_b$ with
$Re (p)\ne 0$, where $f^s$ is a $(1,b)$-quasi-integral function in
spherical ${\cal A}_b$-coordinates with $f=f^s\circ E^{-1}$ (see
Definitions 2 and 2.1), $2\le b\le 3$, $0<q<\infty $, every zero
$z_0$ of $f|_{\bf C}$ may be only in the strip $\{ z\in {\bf C}: -q
\le Re (z) \le q \} $, $f(z)$ has not any real zero. Then all zeros
of the restriction $f|_{\bf C}$ of $f$ on the complex field $\bf C$
are complex and belong to the line $Re (z)=0$.}
\par {\bf Proof.} Let $z_0$ be a complex zero of $f|_{\bf C}$,
$f(z_0)=0$, then $-q\le Re (z_0)\le q$ and $Im (z_0)\ne 0$ by the
supposition of this theorem. Put $v_0 := Re (z_0)$. In the case
$v_0=0$ there is nothing to prove. So suppose that $v_0\ne 0$, then
$f(-z_0)=0$, $f({\tilde z}_0)=0$ and $f(-{\tilde z}_0)=0$ due to the
symmetry properties of $f$, since $f^s(p)=f^s(-p)$ for each $p\in
{\cal A}_b$ with $Re (p)\ne 0$, $E(y)=y$ for each $y\in \bf C$,
$f=f^s\circ E^{-1}$. Thus without loss of generality consider $0<
v_0 \le q$. Hence zero surfaces of the $(1,b)$-quasi-conformal
extension of $f(z)$ are: $S^{f}_{z_0}=:S_{z_0}$ and
$S^{f}_{-z_0}=:S_{-z_0}$ (see Lemma 12), since $z_0, {\tilde z}_0\in
S^{f}_{z_0}$ and $-z_0, -{\tilde z}_0\in S^{f}_{-z_0}$. These
surfaces $S_{z_0}$ and $S_{-z_0}$ are symmetric relative to the
hyperplane $\pi _0 := \{ z\in {\cal A}_b: Re (z)=0 \} $. Without
loss of generality put $im (z_0)>0$, where $im (z_0)=i_1^*Im (z_0)$,
$Im (z_0)=z_0-Re (z_0)=i_1 im (z_0)$.
\par  Since the quaternion skew field ${\bf H}={\cal A}_2$
has the natural embedding into the octonion algebra ${\bf O}={\cal
A}_3$, then it is sufficient to prove this theorem for $b=2$.
Mention, that $Ln$ has the countable number of branches with the
noncommutative Riemann surface $\cal R$ given in \S 10.2. Therefore,
$Ln_n(z_1,...,z_{n-1},1;z_n)$ has the noncommutative Riemann surface
embedded into ${\cal R}^n$, since $z\mapsto a^{-1}z$ is the
pseudo-conformal mapping by $z$ for $a\ne 0$ in ${\cal A}_b$,
$z\mapsto 1/z$ is also pseudo-conformal for $z\ne 0$ (see Corollary
2.7 \cite{lusmnfgpcd}), where $z_j\ne 0$ and $z_j\in {\cal A}_b$ for
each $j=1,...,n$, consequently, $E^{-1}$ has the noncommutative
Riemann surface ${\cal R}_E$ of dimension $2^b$ over $\bf R$
embedded into ${\cal R}^n$ (see Equations 2.1$(1,2)$).
\par Take $b=2$ and consider the loops
\par $(1)$ $q_j(t)=v_j + \rho _j
\exp (\pi K_j\exp (2\pi N_jt)/2)$ \\
encompassing $S_{z_j}$ parameterized by $t\in [0,1)\subset {\bf R}$,
where $j=0$ or $j=1$, $v_1= -v_0$, $K=i_1={\bf i}$, $N=i_3$,
$|N|=1$, $N$ is a marked purely imaginary quaternion orthogonal to
$K$, $Re (NK^*)=0$, $0<\rho _j -|Im (z_0)|$ is sufficiently small,
$\rho _0=\rho _1$ (see Theorem 16.1). For $j=0$ take $K_0=K$ and
$N_0=N$ and for $j=1$ take $K_1=-K$ and $N_1=-N$. Consider spheres
$S({\bf H},v_0,|Im (z_0)|)$ and $S({\bf H},-v_0,|Im (z_0)|)$, where
$S({\cal A}_b,x,R):= \{ y\in {\cal A}_b: |y-x| = R \} $, $0<R<\infty
$. \par In accordance with Theorem I.20.2 \cite{lavrshab} if $D$ is
an open connected domain in $\bf C$ and functions $f_1$ and $f_2$
are holomorphic on $D$ such that $f_1(x_n)=f_2(x_n)$ for each $n\in
\bf N$ and there exists a limit point $x\in D$ of a sequence $ \{
x_n: n\in {\bf N} \} \subset D$, then $f_1(y)=f_2(y)$ for each $y\in
D$. The function $f$ is holomorphic on $\bf C$ and
$(1,b)$-quasi-integral, hence by the latter theorem zeros of $f$ and
$f'$ in $\bf C$ are isolated. Therefore, there exist $0<\rho <\infty
$ and $0<\delta <|Im (z_0)|$ such that for each other complex zero
$z_2\notin \{ z_0, -z_0, {\tilde z}_0, - {\tilde z}_0 \} $ not
belonging to $S({\bf C},v_0,|Im (z_0)|)\cup S({\bf C},-v_0,|Im
(z_0)|)$, either $ \{ z_2, -z_2, {\tilde z}_2, -{\tilde z}_2 \}
\subset B({\cal A}_b,v_0,|Im (z_0)|-\delta )\cup B({\cal
A}_b,-v_0,|Im (z_0)|-\delta )$ or $ \{ z_2, -z_2, {\tilde z}_2,
-{\tilde z}_2 \} \subset {\cal A}_b\setminus [B({\cal A}_b,v_0,\rho
+\delta )\cup B({\cal A}_b,-v_0,\rho +\delta )]$. Take $\rho _+ =
\rho _0$ and $\rho _- = \rho -\delta $ in Formulas 16.1$(1,2)$ with
$v_0=Re (z_0)$ and $v_1= -v_0=Re (-z_0)$ as data in place of $Re
(z_0)$ for $\gamma $ there and get two loops $\gamma _0$ and $\gamma
_1$ corresponding to that of given by Equations $(1,2)$ for $\rho
_+$ and $\rho _-$, denote them by $\gamma _{j,+}$ and $\gamma
_{j,-}$ for $j=0, 1$ respectively. Then $\gamma _0$ and $\gamma _1$
join by a rectifiable curve $\eta $ not containing any zero of $f$
such that $\eta $ is gone twice in one and the opposite direction.
\par There is the identity \par  $K\exp (2\pi Nt)= K\cos (2\pi t) +
K N \sin (2\pi t)= \exp (\pi K\exp (2\pi Nt)/2)$, \\ since $|K\cos
(2\pi t) + KN \sin (2\pi t)|=1$ and $e^M=\cos (|M|) +M \sin
(|M|)/|M|$ for each $M\in {\cal I}_b\setminus \{ 0 \} $, $\sin (\pi
/2)=1$, consequently, $q_j(t)$ is orthogonal to $\bf R$ in $\bf H$
relative to the scalar product $(z,\xi ):= Re (z\xi ^*)$, moreover,
$q_0(0)=z_0+ (\rho _0-\rho ')i_1$, $q_1(0)=-z_0- (\rho _0-\rho
')i_1$, since without loss of generality put $im (z_0)=\rho '>0$,
where $z_0=v_0+ im (z_0)i_1$. Consider also circles $q_{2+j}(t) =
v_j + \rho _2\exp (\pi K_j \exp (2\pi N_jt)/2)$. It is supposed that
$\gamma $ is written in the $z$-representation with the help of
formulas 2.$(2-5)$.
\par  Put $p_w(t) := v_0 +\rho _w i_1 - 2\pi
i_2t$, where $\rho _0:=\rho _+$ and $\rho _2:=\rho _-$, $w=0$ and
$w=2$, then $E_2(p_w(t)) = q_w(t)$ and $E_2(-p_w(t))=q_{1+w}(t)$ for
each $t\in [0,1]$ (see Definition 2.1 and Formula 2.1$(1)$), since
$v_1=-v_0\ne 0$.

Therefore, \par $(2)$ $f(q_w(t)) =f(q_{1+w}(t))$\\
for each $t\in [0,1]$ and for $w=0, 2$, since $f^s(p) = f^s(-p)$ for
each $p\in {\cal A}_b$ with $Re (p)\ne 0$, $f=f^s\circ E_2^{-1}$ by
the conditions of this theorem.

\par Consider the loop $\gamma $ consisting of $q_j(t)$ and twice
gone paths joining them such that $\gamma $ is gone clock-wise by
$q_1$ and $q_2^-$ and counter-clock-wise by $q_0$ and $q_3^-$ in the
planes $v_j + K{\bf R}\oplus KN{\bf R}$ as seen from the negative
axis of $((-\infty ,0)N)$ perpendicular to the real axis: $\gamma
(t)=\gamma _0(4t)$ for $0\le t<1/4$, $\gamma (t)=\gamma _2(
4(t-1/4))$ for $1/4\le t<1/2$, $\gamma (t)=\gamma _1(4(t-1/2))$ for
$1/2\le t<3/4$ and $\gamma (t)=\gamma _2(1-4(t-3/4))$ for $3/4\le
t<1$, where $\gamma _0(t)$ and $\gamma _1(t)$ are composed from
$q_0, q_2^-$ and $q_1, q_3^-$ for $S_{z_0}$ and $S_{-z_0}$
respectively and the joining them paths gone twice in one and the
opposite directions as in the proof of Theorem 16.1, the rectifiable
curve $\{ \gamma _2(t): 0\le t\le 1 \} $ joins $q_0(1)$ with
$q_1(0)$ such that $\gamma ([0,1])\subset V:= \{ z\in {\bf H}: -q\le
Re (z)\le q \} $, where $\gamma (0)=\gamma (1)$, $\gamma
_{0,+}=q_0$, $\gamma _{1,+}= q_1$, $\gamma _{0,-}=q_2^-$, $\gamma
_{1,-}=q_3^-$, since $K(KN)= -N$. Instead of talking about
orientations it is sufficient to write analytic formulas for curves,
that is done in this section.
\par Since $q$ and $|z|$ are finite, then the curve $\gamma $ can be
chosen rectifiable. If $z_2$ is some other zero in the circles
$S({\bf C},v_0,|Im (z_0)|)\cup S({\bf C},-v_0,|Im (z_0)|)$, then
$S_{z_2}$ and $S_{-z_2}$ are encompassed by $\gamma $ as well in the
sense of Theorem 16. Their additional value to $p$ of $\Delta
_{\gamma }Arg_2f$ will be positive together with $z_0$ and $-z_0$ in
accordance with Theorem 16.1. Let $p_0$ be that part of $p$, which
corresponds to $z_2\in S({\bf C},v_0,|Im (z_0)|)\cup S({\bf
C},-v_0,|Im (z_0)|)$ with $Re (z_2)=0$. Denote the set of such zeros
of $f$ by $Z$, $Z := \{ z\in S({\bf C},v_0,|Im (z_0)|)\cup S({\bf
C},-v_0,|Im (z_0)|): Re (z)=0, f(z)=0 \} $ . Then $Z$ is finite and
may happen to be empty, since $Re (z_0)=v_0\ne 0$ by the supposition
made above. If $z\in \bf C$ and $f(z)=0$, then the imaginary part of
$z$ is nonzero, $Im (z)\ne 0$, since $f$ has not any real zero by
the supposition of this theorem. \par There is not any zero of $f$
outside the band $-q\le Re (z)\le q$ in $\bf C$, hence there are not
any chains of crossing spheres around zeros of $f$ of the type
considered above besides may be pairs of spheres with centers at
$v_0$ and $-v_0$ with $0<v_0\le q$. Indeed, for $z_0\in \bf C$ with
$|Im(z_0)|>q$ it may be only two such spheres with a given $z_0$,
$Re (z_0)=v_0$. In the domain $V_q := \{ x\in {\bf C}: |Re (x)|\le
q$  $\mbox{and}$  $|im(x)|\le q \} $ only finite number of zeros of
$f$ may be and the consideration reduces to pairs of spheres if
there exists a zero $z_0\in V_q$ of $f$. If $z_3\in \bf C$ is a zero
of $f$, then there exists an open neighborhood $W$ of $z_3$ such
that $W$ can intersect no more, than a finite family of circles
$S({\bf C},v_j,|im (z_0)|)$, where $v_0=Re (z_0)$, $v_1=-v_0$, $j=0$
or $j=1$, $z_0$ is a zero of $f|_{\bf C}$ different from four zeros
corresponding to $z_3$, $z_0\notin \{ z_3, -z_3, {\bar z}_3, - {\bar
z}_3 \} $. Therefore, the claimed loop $\gamma $ exists for each
complex zero $z_0$ of $f$.
\par For each zero $z_2\in
S({\bf C},v_0,|im (z_0)|)\cup S({\bf C},-v_0,|im (z_0)|)\setminus Z$
if it exists and for the marked $z_0$ of $f$ the symmetry imposed on
$f$ leads to the contradiction, when $Re (z_0) \ne 0$. To
demonstrate this denote the family of such zeros $z_2$ and $z_0$ of
$f$ by $Y$, $Y:= \{ z\in S({\bf C},v_0,|im (z_0)|)\cup S({\bf
C},-v_0,|im (z_0)|): Re (z)\ne 0, f(z)=0 \} $. For $z\in Y$ let
$k(z)\in \bf N$ denotes its order. Then $Y$ is a finite set, since
$f$ is the nontrivial integral function on $\bf C$. If $z\in Y$,
then $-z, {\bar z}, -{\bar z}\in Y$ also and $k(z)=k(-z)=k({\bar
z})=k(-{\bar z})$. Then $p=p_0+p_Y$, where $p_Y\ge 1$ and $p_0\ge 0$
($p_0\ge 1$ when $Z\ne \emptyset $) correspond to the sets $Y$ and
$Z$ respectively. Applying Theorems 16 and 16.1 for $r=1$ and $b=2$
we get for $\gamma $ that
\par $(3)$ $\Delta _{\gamma }Arg_2f=\Delta _{\gamma _0}Arg_2f+\Delta
_{\gamma _1}Arg_2f=0$, \\
since $\Delta _{\gamma _0}Arg_2f$ and $\Delta _{\gamma _1}Arg_2f$
are independent from $a_1=K$ and $a_1=-K$ respectively, but
$N_0=-N_1$. On the other hand,
\par $(4)$ $\int_{\gamma _0}dLn_2(\pi K_0/2,1;f(z)) = \int_{\gamma
_1} dLn_2(\pi K_1/2,1;f(z))$, \\
since $\int_{q_j}dLn_2(\pi K_j/2,1;f(z))=\int_0^1dLn_2(\pi
K_j/2,1;f(q_j(t))$ and the Equality $(2)$ is satisfied. Therefore,
the application of Theorems 16 and 16.1 to $\gamma _0$ only due to
Equation $(4)$ gives \par $(5)$ $|\Delta _{\gamma }Arg_2f| = k\pi $,
\\ where $k=k_0+k_Y$, $k_0\ge 0$ corresponds to zeros from $Z$, while
$k_Y=\sum_{z\in Y}k(z) \ge 4$ corresponds to zeros of $f$ from $Y$.
This gives the contradiction of Equation $(3)$ with $(5)$,
consequently, all complex zeros of $f$ may lie only on the line $Re
(z)=0$.
\par {\bf 17.1. Remarks.} Examples of quasi-regular and
quasi-integral functions can be provided with the help of
Proposition 9.1, Corollary 9.2 and Theorem 9.4 and Section 9.5.
Certainly each pseudo-conformal function on $U\setminus A$ or ${\cal
A}_b\setminus A$ besides a set $A$ consisting of isolated points of
zeros of its derivative $A := \{ z\in U: f'(z)=0 \} $ is at the same
time quasi-regular on $U$ or quasi-integral on ${\cal A}_b$
respectively (see about pseudo-conformal functions in
\cite{ludagpd,lusmnfgpcd}).
\par On the other hand, if $a>0$, $q>0$, then put
$P(x)= (x-a-qi) ({\bar x}-a+qi) (x+a-qi) ({\bar x}+a+qi)$ for $x\in
\bf C$. The polynomial $P$ satisfies the necessary symmetry
properties on $\bf C$, but it has not a quasi-regular extension on
$U$ open in ${\cal A}_b$ with $W=U\cap {\bf C}$ open in $\bf C$ for
$2\le b\le 3$, since the left and the right sides of Formula $(Q7)$
for ${\hat R}_{z,y}P$ differ on terms such as $-q^2(Mv+vM) (z+a-qM)
({\tilde z} +a+qM)- (z-a-qM)({\tilde z} - a +qM)q^2(Mv+vM) +
q(zv-v{\tilde z})(z+a-qM) ({\tilde z} +a+qM)+ q(z-a-qM)({\tilde z} -
a +qM) (zv-v{\tilde z})$, where $z\in {\cal A}_b\setminus \bf C$,
$z-Re (z) := Im (z)\ne 0$, $M=Im (z)/|Im (z)|$, $v\ne 0$, $v \| M$,
$v\in {\cal I}_b\setminus \bf C$, $Re (v)=0$, which follows from the
$z$-representation with the help of Formulas 2.$(1-9)$. Then
functions of the type $(f_1P)f_2$ also generally need not satisfy
conditions of Theorem 17 even when $f_1$ and $f_2$ are
quasi-integral. Therefore, the class of functions satisfying
conditions of Theorem 17 is rather narrow. The graph $\{ (z,f(z)):
z\in {\cal A}_b \} \subset {\cal A}_b^2$ of $f$ satisfying
conditions of Theorem 17 has the natural interpretation. On the
other hand, in view of Corollary 9.2 the power $f^n$ of $f$
satisfies conditions of Theorem 17 for each $n=2,3,4,...$ if $f$
satisfies them.

\par The existence
problem of functions satisfying conditions of Theorem 17 is
considered in the next section together with properties of their
noncommutative integral transformations of the Laplace and Mellin
types.
\par {\bf 18. Theorem.} {\it Let $f(z,t)$ be an ${\cal A}_b$ valued
function on $W := U\times [a,\infty )$ and there exists $(\partial
f(z,t)/\partial z).h$ continuous on $W\times B({\cal A}_b,0,1)$,
where $U = \{ z\in {\cal A}_b: z_j\in [a_j,b_j], j=0, 1,..., 2^b-1
\} $, $z=z_0i_0+...+ z_{2^b-1}i_{2^b-1}$, $a_j<b_j$, $B({\cal
A}_b,y,R) := \{ z\in {\cal A}_b: |z-y|\le R \} $. Let $x\in U$ be
such that $F(x) := \int_a^{\infty }f(x,t)dt$ converges, while the
improper integral depending on the parameter $z\in U$:
$G(z,h):=\int_a^{\infty } (\partial f(z,t)/\partial z).hdt$
converges uniformly on $U\times B({\cal A}_b,0,1)$. Then the
improper integral $F(z) :=\int_a^{\infty }f(z,t)dt$ depending on the
parameter $z\in U$ converges uniformly on $U$ and for each $z\in U$
there exists $D_zF(z).h = \int_a^{\infty } (\partial f(z,t)/\partial
z).hdt = G(z,h)$ for each $h\in B({\cal A}_b,0,1)$.}
\par {\bf Proof.} Write $(\partial f(z,t)/\partial z).h$ in the form
$(\partial f(z,t)/\partial z).h=g_0i_0+...+g_{2^b-1}i_{2^b-1}$,
where $g: W\times {\cal A}_b\to \bf R$. Since $(\partial
f(z,t)/\partial z).h$ is continuous on $W\times B({\cal A}_b,0,1)$
and the improper integral $G(z,h) := \int_a^{\infty } (\partial
f(z,t)/\partial z).hdt$ converges uniformly on $U\times B({\cal
A}_b,0,1)$, then
\par $\int_x^z G(y,h)dy = \int_x^z (\int_a^{\infty }
(\partial f(y,t)/\partial y).hdt) dy$
\par $=\int_a^{\infty } dt (\int_x^z(\partial f(y,t)/\partial y).hdy) $ \\
and the improper integral on the right converges uniformly on
$U\times B({\cal A}_b,0,1)$, since $t$ is the real parameter and
$\bf R$ is the center of ${\cal A}_b$. Take $h=w$, $|w|=1$, such
that $z=x+vw$, where $v>0$. Then $\int_a^{\infty } dt
(\int_x^z(\partial f(y,t)/\partial y).wdy) =\int_a^{\infty }
f(z,t)dt - \int_a^{\infty } f(x,t)dt= F(z)-F(x)$, consequently, the
improper integral $F(z) := \int_a^{\infty }f(z,t)dt$ converges
uniformly on $U$ and $F(z) - F(x) = \int_x^zG(y,w)dy$. Using the
additivity of the integral we have $F(z)- F(\eta )=\int_{\eta
}^zG(y,h)dy$ for each $z, \eta \in U$, where $z - \eta =vh$, $v>0$.
Thus $\partial z(v)/\partial v=h$ and $\int_{\eta }^zG(y,h)dy =
\int_0^vG(\eta +qh,h)dq$, where $q\in [0,v]$. Therefore, $D_zF(z).h
= G(z,h)$ for each $z\in U$ and $h\in {\cal A}_b$, since $(\partial
f(z,t)/\partial z).(sh)= s (\partial f(z,t)/\partial z).h$ for each
$s\in \bf R$, every $h\in {\cal A}_b$ and each $t\in [a,\infty )$
(see for comparison the commutative case in \S IV.2.4 \cite{kamyn}).
\par {\bf 19. Definition.} Define a $n$-residue operator of a function
$f$ holomorphic in $U\setminus \{ z_0 \} $ for some open ball $U$
with the center $z_0$ in ${\cal A}_b$ as
\par $res_n(z_0;f).M := (2\pi )^{-1} (\int_{\gamma _n}
f(z_0+Ln_{n-1}(a_1,...,a_{n-1};(z-z_0))dLn_{n-1}(a_1,...,a_{n-1};
(z-z_0))$, \\
where $\gamma _n(t)=z_0+ \rho Exp_n(a_1,...,a_{n-1},1;2\pi Mt)$,
$M\in {\cal I}_b$, ${\cal I}_b := \{ M\in {\cal A}_b: Re (M)=0 \} $,
$|M|=1$, $2\le b\le 3$, $0<\rho <\infty $, $0\le t\le 1$; $1\le
n<2^b$, $Exp_0(z) := id(z)=z$, $Ln_0(z):=id(z)=z$, $Exp_n$ and
$Ln_n$ are given by Definition 10; $a_1,...,a_{n-1}\in {\cal A}_b$,
$a_1\ne 0$,..., $a_{n-1}\ne 0$, $\gamma ([0,1])\subset U$, where
$Ln_{n-1}(a_1,...,a_{n-1};Exp_{n-1}(a_1,...,a_{n-1};z))=z$ is such
branch of $Ln_{n-1}$. Then extend this operator by the formula
$res_n(z_0;f)M := [res_n(z_0;f).(M/|M|)]|M|$ for $M\ne 0$ in ${\cal
I}_b$, $res_n(z_0;f).0 := 0$.
\par For $n=1$ this coincides with the usual definition of
$res(z_0;f)$ \cite{luoyst,lusmfcdv,luhcnfcdv}.
\par {\bf 20. Proposition.} {\it Let $g$ be a pseudo-conformal function
in $V\setminus \{ z_0 \} $, where $V$ is an open ball with the
center $z_0$ in ${\cal A}_r$, $z_0\in {\cal A}_r$, $1\le r\le 2$.
Suppose that a $(r,b)$-quasi-conformal function $f$ in $U\setminus
\{ z_0 \} $ is obtained from $g$ in accordance with Conditions
$(Q1-Q7)$, where $U$ is an open ball with the center $z_0$ in ${\cal
A}_b$, $r<b\le 3$, $2\le n<2^b$. Then $res_n(z_0;f).M=res(z_0;g).M$
for each $M\in {\cal I}_r$.}
\par {\bf Proof.} Take $M\in {\cal I}_b$ with $|M|=1$.
For the fixed branch of the logarithmic function consider the
integral $\int_{\gamma _n}
f(z_0+Ln_{n-1}(a_1,...,a_{n-1};(z-z_0)))dLn_{n-1}(a_1,...,a_{n-1};
(z-z_0))$, where $Ln_{n-1}(a_1,...,a_{n-1};
Exp_{n-1}(a_1,...,a_{n-1};z))=z$ is such branch of $Ln_{n-1}$, $\exp
(2\pi Mt)$ is periodic. We have $Ln_{n-1}(a_1,...,a_{n-1};\gamma
_n(t)-z_0)= exp (2\pi Mt)$ for each $t\in \bf R$, when $a_1\ne
0,...,a_{n-1}\ne 0$. Choose $\rho >0$ sufficiently small such that
$\gamma ([0,1])\subset V$. Therefore, $\int_{\gamma _n}
f(z_0+Ln_{n-1}(a_1,...,a_{n-1};(z-z_0)))dLn_{n-1}(a_1,...,a_{n-1};
(z-z_0))=\int_{\gamma } f(z)dz$, where $\gamma (t) := z_0+\rho \exp
(Mt)$ for each $t\in \bf R$. We have that $\gamma ([0,1])\subset
{\cal A}_r$ and $res_n(z_0;f)0=0$ and $res(z_0;g)0=0$. Since for
$M\in {\cal I}_r$ and such $\gamma $ there is the equality:
\par $\int_{\gamma } g(z)dz=\int_{\gamma } f(z)dz$, then
\par $\int_{\gamma } g(z)dz= \int_{\gamma _n}
f(z_0+Ln_{n-1}(a_1,...,a_{n-1};(z-z_0))dLn_{n-1}(a_1,...,a_{n-1};
(z-z_0))$, \\
consequently, $res_n(z_0;f)|_{{\cal I}_r}= res(z_0;g)|_{{\cal
I}_r}$, since $\gamma ([0,1])$ is the loop gone around once and
$res_n(z_0;f)M := [res_n(z_0;f).(M/|M|)]|M|$ for each $M\ne 0$ in
${\cal I}_b$.
\par Mention that $res_n(z_0;f)$ is independent of
$a_1\ne 0,...,a_{n-1}\ne 0$ in ${\cal A}_b$, since \\
$Ln_{n-1}(a_1,...,a_{n-1};Exp_{n-1}(a_1,...,a_{n-1};z))=z$ for each
$z\in {\cal A}_b$ for a branch of $Ln$ such that $Ln (1)=0$. Using
the homotopy theorem for ${\cal A}_b$ line integrals we can consider
more general $U$ than balls (see Theorems 2.11 and 3.9 in
\cite{luoyst,lusmfcdv,luhcnfcdv}).

\section{Noncommutative integral transformations over
$\bf H$ and $\bf O$.}
\par {\bf 1. Definitions.} A function $f: {\bf R}\to {\cal A}_b$ we call
function-original, where ${\cal A}_b$ is the Cayley-Dickson algebra,
which may be, in particular, ${\cal A}_2=\bf H$ over the quaternion
skew field or ${\cal A}_3=\bf O$ the octonion algebra, if it
satisfies the following conditions $(1-3)$:
\par $(1)$ $f(t)$ satisfies the H\"older condition: $|f(t+h)-f(t)|
\le A |h|^{\alpha }$ for each $|h|<\delta $ (where $0<\alpha \le 1$,
$A=const >0$, $\delta >0$ are constants for a given $t$) everywhere
on $\bf R$ may be besides points of discontinuity of the first type.
On each finite interval in $\bf R$ the function $f$ may have only
the finite number of points of discontinuity of the first kind.
Remind, that a point $t_0$ is called the point of discontinuity of
the first type, if there exist finite left and right limits
$\lim_{t\to t_0, t<t_0} f(t) =: f(t_0-0)\in {\cal A}_b$ and
$\lim_{t\to t_0, t>t_0} f(t) =: f(t_0+0)\in {\cal A}_b$.
\par $(2)$ $f(t)=0$ for each $t<0$.
\par $(3)$ $f(t)$ increases not faster, than the exponential function,
that is there exist constants $C=const>0$, $s_0=s_0(f)\in \bf R$
such that $|f(t)|<C \exp (s_0t)$ for each $t\in \bf R$.
\par If there exists an original \par $(4)$ $F(p;q):=F(p)
:=\int_0^{\infty } f(t)e^{-u(p,t;q)}dt$, \\
then $F(p)$ is called the Laplace transformation at a point $p\in
{\cal A}_b$ of the function-original $f(t)$, where either
$u(p,t;q)=pt+q$ or $u(p,t;q)=E(pt+q)$ (see Definition 2.2.1), $p\in
{\cal A}_b$, $q\in {\cal A}_b$ is a parameter. It is supposed that
the integral for $F(p;q)$ is written in the $(p,q)$-representation
with the help of Formulas 2.2.$(2-5)$. For $q=0$ it can be omitted
from $u$ putting $u(p,t;0)=u(p,t)$. If $q$ is specified, then it can
also be written shortly $F(p)$ instead of $F(p;q)$.
\par It can be taken the automorphism of the algebra ${\cal A}_b$
and instead of the standard generators $ \{ i_0,...,i_{2^b-1} \} $
use new generators $ \{ N_0,...,N_{2^b-1} \} $. Provide also
$u(p,t;q )=u_N(p,t;q)$ relative to a new basic generators, where
$2\le b\le 3$. In this more general case an image we denote by
$\mbox{ }_NF(p)$ for the original $f(t)$ or in more details we
denote it by $\mbox{ }_N {\cal F}(f;p;q)$ or $\mbox{ }_NF_u(p;q)$.
\par Let $\gamma : (-\infty ,\infty )\to {\cal A}_b$ be a path such that
the restriction $\gamma_l :=\gamma |_{[-l,l]}$ is rectifiable for
each $0<l\in \bf R$ and put by the definition \par $(5)$
$\int_{\gamma }f(z)dz=\lim_{l\to \infty }\int_{\gamma _l}f(z)dz$, \\
where ${\cal A}_b$ integrals by rectifiable paths were defined in \S
2.5 \cite{luoyst,lusmfcdv,luhcnfcdv}. So they are defined along
curves also, which may be classes of equivalence of paths relative
to increasing piecewise smooth mappings $\tau : [a,b]\to [a_1,b_1]$
realizing reparametrization of paths. Then we shall talk, that an
improper integral (5) converges. \par Consider now a function
$f(z,y)$ defined for all $z$ from a domain $U$ and for each $y$ in a
neighborhood $V$ of a curve $\gamma $ in ${\cal A}_b$. The integral
$G(z) := \int_{\gamma }f(z,y)dy$ converges uniformly in a domain
$U$, if for each $\epsilon >0$ there exists $l_0>0$ such that
\par $(6)$ $|\int_{\gamma }f(z,y)dy - \int_{\gamma _l}
f(z,y)dy|<\epsilon $ for each $z\in U$ and $l>l_0$. Analogously is
considered the case of unbounded $\gamma $ in one side with
$[0,\infty )$ instead of $(-\infty ,\infty )$.
\par {\bf 2. Theorem.} {\it Let $V$ be a bounded neighborhood
of a rectifiable curve $\gamma $ in ${\cal A}_r$, and a sequence of
functions $f_n: V\to {\cal A}_r$ be uniformly convergent on $V$,
where $2\le r<\infty $, then there exists the limit
\par $(i)$ $\lim_{n\to \infty }\int_{\gamma } f_n(z)dz=\int_{\gamma
}\lim_{n\to \infty }f_n(z)dz$.}
\par {\bf Proof.} For a given $\epsilon >0$ in view of
the uniform convergence of the sequence $f_n$ there exists $n_0\in
\bf N$ such that $|f_n(z)-f(z)|<\epsilon /l$ for each $n>n_0$, where
$0<l<\infty $ is the length of the rectifiable curve $\gamma $,
$f(z) := \lim_{n\to \infty }f_n(z)$. In view of the Inequality
2.7(4) \cite{luoyst,lusmfcdv,luhcnfcdv} there are only two positive
constants $C_1>0$ and $C_2>0$ such that $|\int_{\gamma }f(z)dz -
\int_{\gamma }f_n(z)dz|< (\epsilon /l) lC_1 \exp (C_2R^s)=\epsilon
C_1\exp (C_2R^s)$, where $s=2^r+2$, $0<R<\infty $, such that $V$ is
contained in the ball $B({\cal A}_r,z_0,R)$ in ${\cal A}_r$ of the
radius $R$ with the center at some point $z_0\in \bf K$. This means
the validity of Equality $(i)$.
\par {\bf 3. Theorem.} {\it If a function $f(z,y)$ holomorphic by
$z$ is piecewise continuous by $y$ for each $z$ from a simply
connected (open) domain $U$ in ${\cal A}_r$ with $2\le r<\infty $
and for each $y$ from a neighborhood $V$ of the path $\gamma $,
where $\gamma _l$ is rectifiable for each $0<l<\infty $, and the
integral $G(z):=\int_{\gamma }f(z,y)dy$ converges uniformly in the
domain $U$, then it is  the holomorphic function in $U$.}
\par {\bf Proof.} For each $0<l<\infty $ the function
$\int_{\gamma _l} f(z,y)dy =: G_l(z)$ is continuous by $z$ in view
of Theorem 2, that together with 1(6) in view of the triangle
inequality gives the continuous function $G(z)$ on $U$. In view of
Theorems 2.16 and 3.10 \cite{luoyst,lusmfcdv,luhcnfcdv} the integral
holomorphicity of the function $G(z)$ implies its holomorphicity.
But the integral holomorphicity is sufficient to establish in the
interior $Int (B({\cal A}_r,z_0,R))$ of each ball $B({\cal
A}_r,z_0,R)$ contained in $U$.  Let $\psi $ be a rectifiable path
such that $\psi \subset Int (B({\cal A}_r,z_0,R))$. Therefore,
$\int_{\psi } G(z)dz = \int_{\psi } (\int_{\gamma }f(z,y)dy)dz$.
With the help of the  proof of Theorem 2.7
\cite{luoyst,lusmfcdv,luhcnfcdv} these integrals can be rewritten in
the real coordinates and with the generators $i_0,...,i_{2^r-1}$ of
the Cayley-Dickson algebra ${\cal A}_r$, since $f=\sum_{k=0}^{2^r-1}
f_ki_k$, where $f_k\in \bf R$ for each $k$. In view of the uniform
convergence $G(z)$ and the Fubini Theorem it is possible to change
the order of the integration and then $\int_{\psi } G(z)dz =
\int_{\gamma } (\int_{\psi }f(z,y)dz )dy=0$, since $\int_{\psi
}f(z,y)dz=0$.
\par {\bf 4. Theorem.} {\it For each original $f(t)$ its
image $F(p)$ is defined in the half space $\{ p\in {\cal A}_r: Re
(p)>s_0 \} $, moreover, $F(p)$ is holomorphic by $p$ in this half
space, where $2\le r\le 3$ and the indicator of the growth of $f(t)$
is not greater, than $s_0$.}
\par {\bf Proof.} Integral 1(4) is absolutely convergent
for $ Re (p)>s_0$, since it is majorized by the converging integral
\par $|\int_0^{\infty }f(t)\exp (-u(p,t;q))dt| \le \int_0^{\infty }
C\exp (-(s-s_0)t)dt = C(s-s_0)^{-1}$, \\
since $|e^z|=\exp (Re (z))$ for each $z\in {\cal A}_r$ in view of
Corollary 3.3 \cite{luoyst,lusmfcdv,luhcnfcdv}, where $s=Re (p)$,
$C>0$ is independent from $p$ and $t$. While an integral, produced
from the integral 1(4) differentiating by $p$ (see Theorem 2.18)
converges also uniformly:
\par $(i)$ $|\int_0^{\infty }f(t)[\partial \exp (-u(p,t;q))/
\partial p].hdt|
\le |h| \int_0^{\infty } C t \exp (-(s-s_0)t)dt =
|h|C(s-s_0)^{-2}$ \\
for each $h\in {\cal A}_r$, since each $z\in {\cal A}_r$ can be
written in the form $z=|z|\exp (M)$ in accordance with Proposition
3.2 \cite{luoyst,lusmfcdv,luhcnfcdv}, where $|z|^2=z{\tilde z}\in
[0,\infty )\subset \bf R$, $M\in {\cal A}_r$, $Re (M):= (M+{\tilde
M})/2=0$ . Therefore,
\par $[\partial (\int_0^{\infty }f(t)
\exp (- u(p,t;q)dt)/\partial {\tilde p}].h=0$ \\ for each $h\in
{\cal A}_r$, since $u(p,t;q)$ is written in the
$(p,q)$-representation. In view of convergence of integrals given
above $F(p)$ is (super)differentiable by $p$, moreover, $\partial
F(p)/\partial {\tilde p}=0$ in the considered $p$-representation,
consequently, $F(p)$ is holomorphic by $p\in {\cal A}_b$ with $Re
(p)>s_0$ due to Theorem 3.

\par {\bf 5. Theorem.} {\it If a function $f(t)$ is an original
(see Definition 1), such that $\mbox{ }_NF_u(p;q) :=
\sum_{j=0}^{2^r-1} \mbox{ }_NF_{u,j}(p;q)N_j$ is its image,
where the function $f$ is written in the form \\
$f(t)=\sum_{j=0}^{2^r-1} f_j(t)N_j$, $f_j: {\bf R}\to \bf R$ for
each $j=0,1,...,2^r-1$, $f({\bf R})\subset {\cal A}_r$, $2\le r\le
3$. Then at each point $t$, where $f(t)$ satisfies the H\"older
condition there is accomplished the equality:
\par $(i)$ $f(t) = (2\pi N_1)^{-1} Re (S{\tilde N}_1) \sum_{j=0}^{2^r-1}
(\int_{a-S\infty }^{a+S\infty }\mbox{ }_NF_{u,j}(p;q)
\exp (u(p,t;q))dp)N_j$, \\
where either $u(p,t;q)=pt+q_0$ with $S=N_1$ and $Im (q)=0$, or
$u(p,t;q)=E(pt+q)$, the integral is taken along the straight line
$p(\tau )=a+S\tau \in {\cal A}_r$, $\tau \in \bf R$, $S\in {\cal
A}_r$, $Re (S)=0$, $|S|=1$, $Re (S{\tilde N}_1)\ne 0$ is non-zero,
$Re (p) = a>s_0$ and the integral is understood in the sense of the
principal value.}
\par {\bf Proof.} In view of the decomposition of a function $f$
in the form $f(t)=\sum_{j=0}^{2^r-1} f_j(t)N_j$ it is sufficient to
consider the inverse transformation of the real valued function
$f_j$, which we denote for simplicity by $f$. Since $t\in \bf R$,
then $\int_0^{\infty }f(\tau )d\tau $ is the Riemann integral. If
$w$ is a holomorphic function of the Cayley-Dickson variable, then
locally in a simply connected domain $U$ in each ball $B({\cal
A}_r,z_0,R)$ with the center at $z_0$ of radius $R>0$ contained in
the interior $Int (U)$ of the domain $U$ there is accomplished the
equality
\par $(\partial \int_{z_0}^zw(\zeta )d\zeta /\partial z).1=w(z)$, \\
where the integral depends only on an initial $z_0$ and a final $z$
points of a rectifiable path in $B({\cal A}_r,z_0,R)$. On the other
hand, along the straight line $a+S\bf R$ the restriction of the
antiderivative has the form $\int_{\theta _0}^{\theta }w(a+S\tau
)d\tau $, since \par $\int_{z_0=a+S\theta _0}^{z=a+S\theta }w(\zeta
)d\zeta =\int_{\theta _0}^{\theta } {\hat
w}(a+S\tau ).Sd\tau $, \\
while $\partial f(z)/\partial \theta =(\partial f(z)/\partial z).S$
for super-differentiable by $z\in U$ function $f(z)$, moreover, the
antiderivative is unique up to a constant from ${\cal A}_r$ with the
given representation of the function and the branch of the
noncommutative line integral (for example, specified with the help
of the left or right algorithm) \cite{luoyst,lusmfcdv,luhcnfcdv}.
\par The integral $g_B(t) := \int_{a-SB}^{a+SB} \mbox{ }_NF_{u,j}(p;q
)\exp (u(p,t;q))dp$ for each $0<B<\infty $ with the help of
generators of the algebra ${\cal A}_r$ and the Fubini Theorem for
real valued components of the function can be written in the form:
\par $g_B(t) = (2\pi N_1)^{-1} Re (S{\tilde N}_1)
\int_0^{\infty }f(\tau )d\tau \int_{a-SB}^{a+SB} \exp (u(p,t;q))$
\par $\exp (-u(p,\tau ;q))dp$, \\
since the integral $\int_0^{\infty }f(\tau )\exp (-u(p,\tau
;q))d\tau $ is uniformly converging relative to $p$ in the half
space $Re (p)>s_0$ in ${\cal A}_r$ (see also Proposition 2.18
\cite{lusmtslt,luhcnlt}). In view of the alternativity of the
algebra ${\cal A}_r$ use the automorphism $v$ from Lemma 2.17
\cite{lusmtslt,luhcnlt}. This gives the change of the basis of
generators, hence instead of $\mbox{ }_NF_u(p;q)$ consider $\mbox{
}_KF_u(p;q)$, where $K_j=v(N_j)$ is the new basis of generators of
${\cal A}_r$, $j=1,...,2^r-1$, $N_0=K_0=1$. Then with such $v$ the
function $u_K(p,t;q)=v(E(pt+q))$ has the form given by the formulas:
\par $(1)$ $v(u_N(p,t;q))= u_K(p,t;q)=
(p_0t+q_0)+ ({p'}_1t + {q'}_1)K$, where
\par $(2)$ $K= [K_1 \cos
({q'}_2) + K_2 \sin ({q'}_2) \cos ({q'}_3) +K_3 \sin ({q'}_2) \sin
({q'}_3)]$ for quaternions;
\par $(3)$ $K= [ K_1 \cos ({q'}_2) + K_2 \sin ({q'}_2)$
$\cos ({q'}_3) +...+ K_6 \sin ({q'}_2) ...\sin ({q'}_6) \cos
({q'}_7)$\\  $+ K_7\sin ({q'}_2)...\sin ({q'}_6) \sin ({q'}_7)]$ for
octonions, where $p_0,{p'}_1, q_0, {q'}_1,...,{q'}_{2^r-1}\in \bf
R$, $t\in \bf R$, $K_1,...,K_{2^r-1}\in {\cal A}_r$ are new
generators with $Re (K_j)=0$ for each $j=1,...,2^r-1$, where
$K_0=N_0=1$, ${p'}_0=p_0$ and ${q'}_0=q_0$,
$p=p_0N_0+p_1N_1+...+p_{2^r-1}N_{2^r-1}=p_0K_0+{p'}_1K_1+...+
{p'}_{2^r-1}K_{2^r-1}$, $q=q_0N_0+q_1N_1+...+q_{2^r-1}N_{2^r-1}=
q_0K_0+{q'}_1K_1+...+ {q'}_{2^r-1}K_{2^r-1}$, since $v(1)=1$ and,
consequently, $v(t)=t$ for each $t\in \bf R$. Formula $(i)$ is
satisfied if and only if it is accomplished after application of the
automorphism $v$ to both parts of the Equality, since $v(z)=v(\zeta
)$ for $z, \zeta \in {\cal A}_r$ is equivalent to that $z=\zeta $.
\par Then up to an automorphism of the algebra ${\cal A}_r$
the proof reduces to the case $p=(p_0,p_1,0,...,0)$,
$N=(N_0,N_1,N_2,...,N_{2^r-1})$, where $N_0=1$, since $\bf R$ is the
center of the algebra ${\cal A}_r$. But this gives $p_1=p_1(t)=Re
(S{\tilde N}_1)t$ for each $t\in \bf R$. For $u(p,t;q_0)=pt+q_0$
with $Im (q)=0$ take simply $S=N_1$. Consider the particular case $c
:= Re (S{\tilde N}_1)\ne 0$, then the particular case $Re (S{\tilde
N}_1)=0$ is obtained by taking the limit when $Re (S{\tilde N}_1)\ne
0$ tends to zero. Thus,
\par $g_B(t) = (2\pi N_1)^{-1} c \int_0^{\infty }f(\tau )d\tau
\int_{a-SB}^{a+SB} \exp (
at+c(q_1+ t)K_1) \exp (-(a\tau +c(q_1+ \tau )K_1)dp$, \\
since $q_0, a\in \bf R$, where $K_1$ is ether given by Formulas
$(2,3)$ for $u_N(p,t;q)=E(pt+q)$ or $K_1=N_1=S_1$ can be taken for
$u_N(p,t;q)=pt+q_0$. Then
\par $g_B(t) = (\pi N_1)^{-1} c
\int_0^{\infty }f(\tau )e^{a(t-\tau )}[\sin (Bc(t-\tau ))] (ct-c\tau
)^{-1}$ \par  $ =(\pi )^{-1} e^{at} \int_{-t}^{\infty } f(\zeta +t)
e^{-a(\zeta +t)}[\sin (B\zeta )]\zeta ^{-1}d\zeta $, \\
where it can be used the substitution $\tau -t= \zeta $. Put
$w(t):=f(t)e^{-at}$, where $w(t)=0$ for each $t<0$. Therefore,
\par $g_B(t)=(\pi )^{-1}e^{at}\int_{-\infty }^{\infty }[w(\zeta +t)-
w(t)]\zeta ^{-1}\sin (B\zeta )d\zeta +(\pi )^{-1}f(t)\int_{-\infty
}^{\infty }{\zeta }^{-1}\sin (B\zeta )d\zeta $. The integral in the
second term is known as the Euler integral: $\int_{-\infty }^{\infty
}{\zeta }^{-1}\sin (B\zeta )d\zeta =\pi $ for each $B>0$,
consequently, the second term is equal to $f(t)$. It remains to
prove, that $\lim_{B\to \infty }\int_{-\infty }^{\infty }[w(\zeta
+t)- w(t)]\zeta ^{-1}\sin (B\zeta )d\zeta =0$, that follows from the
subsequent lemma.

\par {\bf 6. Lemma.} {\it If a function $\psi (y)$ with values in
the Cayley-Dickson algebra ${\cal A}_r$ is integrable on the segment
$[\alpha ,\beta ]\subset \bf R$, then \par $\lim_{b\to \infty }
\int_{\alpha }^{\beta }\psi (y)\sin (by)dy =0$.}
\par {\bf Proof.} If $\psi $ is continuously differentiable on the
segment $[\alpha ,\beta ]$, then the result of the integration by
parts is:
\par $\int_{\alpha }^{\beta }\psi (y)\sin (by)dy
=-\psi (y)\cos (by)b^{-1}|_{\alpha }^{\beta }+\int_{\alpha }^{\beta
}\psi '(y)\cos(by)b^{-1}dy$ and, consequently,
\par $\lim_{b\to \infty } \int_{\alpha }^{\beta }\psi (y)\sin
(by)dy =0$. If $\psi (y)$ is an arbitrary integrable function, then
for each $\epsilon >0$ there exists a continuous differentiable
function $\psi _{\epsilon }(y)$ such that $\int_{\alpha }^{\beta
}|\psi (y)- \psi _{\epsilon }(y)|dy <\epsilon /2$. Then
$\int_{\alpha }^{\beta }\psi (y)\sin (by)dy =\int_{\alpha }^{\beta }
[\psi (y)-\psi _{\epsilon }(y)]\sin (by)dy +\int_{\alpha }^{\beta
}\psi _{\epsilon }(y)\sin (by)dy $, where $|\int_{\alpha }^{\beta }
[\psi (y)-\psi _{\epsilon }(y)]\sin (by)dy|<\epsilon /2$ for each
$b$, since $|\sin (by)|\le 1$ and the second term tends to zero:
$\lim_{b\to \infty }\int_{\alpha }^{\beta }\psi _{\epsilon }(y)\sin
(by)dy =0$ by the one proved above.

\par The final part of the {\bf Proof} of Theorem 5. For a fixed
$\epsilon >0$ there is an equality:
\par $\int_{-\infty }^{\infty }[w(y +t)-w(t)]y^{-1}\sin
(By)dy =\int_{-B}^B [w(y +t)-w(t)]y^{-1}\sin (By)dy +
\int_{|y|>B}w(y +t)y^{-1}\sin (By)dy - w(t)\int_{|y|>B}\sin (By )y
^{-1}dy$. The second and the third terms are converging integrals
and therefore for sufficiently large $B>0$ by the absolute value
they are smaller than $\epsilon /3$. In view of the H\"older
condition $|[w(y+t)-w(t)] y^{-1}|\le A|y|^{1-c}$, where $c>0$,
$A>0$, $y$ is in a neighborhood of zero. Then in view of Lemma 6
above there exists $B_0>0$ such that
\par $|\int_{-B}^B [w(y+t)-w(t)]y^{-1}\sin (By)dy|<\epsilon /3$
for each $B>B_0$. Thus,
\par $\lim_{B\to \infty } \int_{-\infty }^{\infty } [w(y
+t)-w(t)]{y}^{-1}\sin (By)dy =0$.
\par This theorem for the general function
$u_N(p,t;q)=E(pt+q)$ in the basis of generators $\{
N_0,...,N_{2^r-1} \} $ follows also directly by the calculation of
appearing integrals by real variables $t$ and $\tau $ using Lemma 6
with the help of integrals evaluated in \cite{lusmtslt,luhcnlt}.

\par {\bf 7. Theorem.} {\it An original $f(t)$ with $f({\bf R})\subset
{\cal A}_r$ for $r=2, 3$ is completely defined by its image $\mbox{
}_NF(p)$ up to values at points of discontinuity.}
\par {\bf Proof.} In view of Theorem 5 the value $f(t)$
at each point $t$ of continuity of $f(t)$ is expressible throughout
$\mbox{ }_NF(p)$ by Formula $5(i)$. At the same time values of the
original at points of discontinuity do not influence on the image
$\mbox{ }_NF(p)$, since on each bounded interval a number of points
of discontinuity is finite.

\par {\bf 8. Theorem.} {\it If a function $\mbox{
}_NF_u(p)$ is analytic by the variable $p\in {\cal A}_r$ in the half
space $W := \{ p\in {\cal A}_r: Re (p)>s_0 \} $, where $2\le r\le
3$, $f({\bf R})\subset {\cal A}_r$, either $u(p,t)=pt$ or $u(p,t) =
E(pt) $, moreover, for each $a>s_0$ there exist constants $C_a
>0$ and $\epsilon _a >0$ such that \par $(i)$ $|\mbox{ }_NF_u(p)|\le
C_a\exp (-\epsilon _a |p|)$ for each $p\in {\cal A}_r$ with $Re
(p)\ge a$, where $s_0$ is fixed, the integral \par $(ii)$
$\int_{a-S\infty }^{a+S\infty }\mbox{ }_NF_u(p)dp$ is absolutely
converging, where $S\in {\cal I}_r$, $|S|=1$. Then $\mbox{
}_NF_u(p)$ is the image of the function
\par $(iii)$ $f(t)=(2\pi )^{-1}{\tilde S}\int_{a-S\infty
}^{a+S\infty }\mbox{ }_NF_u(p)\exp (u(p,t))dp$.}
\par {\bf Proof.} The case $u(p,t) = pt$ follows from
$u(p,t) := E(pt)$, when $p=(p_0,p_1,0,...,0)$, but the integral
along the straight line $a+St$, $t\in \bf R$, with such $p$ in the
basis of generators $(N_0,...,N_{2^r-1})$ can be obtained from the
general integral by an automorphism $v$, $z\mapsto v(z)$, of the
algebra ${\cal A}_r$, $2\le r\le 3$. That is, as in the proof of
Theorem 5 it is sufficient to prove the equality of the type $(iii)$
for $\mbox{ }_KF_u(p)$ after the action of the automorphism $v$.

\par Let $Re (p)=a>s_0$, then \par $|\int_{a-S\infty }^{a+S\infty }
\exp (u(p,t))\mbox{ }_NF_u(p)dp|$  $\le  \int_{-
\infty }^{\infty } |\mbox{ }_NF_u(a+S\theta )|d\theta $. \\
In view of the supposition of this theorem this integral converges
uniformly relative to $t\in \bf R$. For $f(t)$ given by the Formula
$(iii)$ for $Re (\eta ) =: \eta _0>s_0$ and $(\eta -Re (\eta ))=: Im
(\eta )$ parallel to $S$, we get
\par $\int_0^{\infty }f(t)\exp (-\eta t)dt $ \\  $= (2\pi
)^{-1}{\tilde S}\sum_{j=0}^{2^r-1} \int_0^{\infty }(\int_{a-S\infty
}^{a+S\infty }\mbox{
}_NF_{u,j}(p) \exp (u(p,t))dp)\exp (-\eta t)(dt)N_j$, \\
in which it is possible to change the order of the integration,
since $t\in \bf R$. Then
\par $\int_0^{\infty }f(t)\exp
(-\eta t)dt = (2\pi )^{-1}{\tilde S}\sum_{j=0}^{2^r-1}
\int_{a-S\infty }^{a+S\infty }(\int_0^{\infty }\mbox{ }_NF_{u,j}
(p)\exp ((p-\eta )t)dt)(dp)N_j$, \\
since $e^v\in \bf R$ for each $v\in \bf R$,
$e^{aM}e^{bM}=e^{(a+b)M}$ for each $a, b\in \bf R$. In view of
$a<\eta _0$ and
\par $\int_0^{\infty }e^{(p-\eta )t}dt= - (p-\eta )^{-1}$, \\
then \par $\int_0^{\infty }f(t)\exp (-\eta t)dt = - (2\pi
)^{-1}{\tilde S}\sum_{j=0}^{2^r-1} (\int_{a-S\infty }^{a+S\infty
}\mbox{ }_NF_{u,j}(p)(p-\eta )^{-1}dp)N_j$ \par $= - (2\pi
)^{-1}{\tilde S} (\int_{a-S\infty }^{a+S\infty }\mbox{ }_NF_u(p
)(p-\eta )^{-1}dp$.
\par To finish the proof it is necessary the following analog
of the Jordan lemma.

\par {\bf 9. Lemma.} {\it Let a function $F$ of the variable $p$ from
the Cayley-Dickson algebra ${\cal A}_r$ with $2\le r\in \bf N$
satisfy Conditions $(1-3)$: \par $(1)$ the function $F(p)$ is
continuous by the variable $p\in {\cal A}_r$ in an open domain $W$
of the half space $\{ p\in {\cal A}_r: Re (p)>s_0 \} $, moreover for
each $a>s_0$ there exist constants ${C'}_a
>0$ and $\epsilon _a >0$ such that
\par $(2)$ $|F(p)|\le {C'}_a\exp (-\epsilon _a |p|)$ for each
$p\in S_{R_n}$, $S_R := \{ z\in {\cal A}_r: |z|=R, Re (z)\ge a \} $,
$0<R_n<R_{n+1}$ for each $n\in \bf N$, $\lim_{n\to \infty }R_n =
\infty $, where $s_0$ is fixed, the integral
\par $(3)$ $\int_{a-S\infty }^{a+S\infty }F(p)dp$ is absolutely
converging. Then
\par $(4)$ $\lim_{n\to \infty }\int_{\gamma _n} F(p)\exp
(- u(p,t;q))dp=0$ \\
for each $t>0$ and each sequence of rectifiable curves $\gamma _n$
contained in $S_{R_n}\cap W$, moreover either $F(p)$ is holomorphic
in $W$, which is $(2^r-1)$-connected open domain in ${\cal A}_r$
(see \cite{span}), such that the projection $\pi _{s,p,t}(W)$ is
simply connected in $s{\bf R}\oplus p\bf R $ for each $s=i_{2k}$,
$p=i_{2k+1}$, $k=0,1,...,2^{r-1}-1$ for each $t\in {\cal
A}_{r,s,p}:={\cal A}_r\ominus s{\bf R}\ominus p\bf R$ and $u\in
s{\bf R}\oplus p\bf R$, for which there exists $z=t+u\in {\cal
A}_r$; or there exists a constant ${C'}_V>0$ such that the
variations (lengths) of curves are bounded $V(\gamma _n)\le {C'}_V
R_n$ for each $n$, where $n\in \bf N$, either $u(p,t;q)=pt+q$ or
$u(p,t;q)=E(pt+q)$.}
\par {\bf Proof.} If $0< \epsilon \le \min
(a-s_0,\epsilon _a)$, then in view of Condition $(2)$ there exists a
constant  $C' >0$ such that
\par $(5)$  $|F(p)|\le C' e^{-\epsilon |p|}$, \\
for each $p\in {\cal A}_r$ with $Re (p)\ge a>s_0$. If $U$ is a
domain in ${\cal A}_r$ of the diameter not greater than $\rho $,
then in view of (4) from the proof of Theorem 7
\cite{luoyst,lusmfcdv,luhcnfcdv} there is accomplished the
inequality:
\par $\sup_{p\in U} \| {\hat F}(p) \| \le \sup_{p\in U} |F(p)|
C_1\exp (C_2\rho ^n)$, \\
where $C_1$ and $C_2$ are positive constants independent from $F$,
$n=2^r+2$, $2\le r\in \bf N$. In particular, as $U$ it is possible
to take the interior of the parallelepiped with ribs of lengths not
greater than $\rho /2^{r/2}$. Then the path of integration can be
covered by a finite number of such parallelepipeds. In the case of
the circle of radius $R$ a number of necessary parallelepipeds is
not greater, than $2^{1+r/2}\pi R/\rho +1$. There exists $R_0>0$
such that for each $R>R_0$ there is accomplished the inequality
$2^{1+r/2}\pi R/\rho <\exp (C_2\rho ^{n-1}(R-\rho ))$. Therefore, in
$\rho $ neighborhood $C_R^{\rho }$ of the circle $C_R$ of radius $R$
and with the center at zero with $R>R_0$ there is accomplished the
inequality:
\par $\sup_{p\in C_R^{\rho }} \| {\hat F}(p) \| \le
\sup_{p\in C_R^{\rho }} |F(p)| C_1\exp (C_2\rho ^{n-1}R)$, \\
where $A^{\rho } := \{ z\in {\cal A}_r: \inf_{a\in A} |a-z|< \rho \}
$ for a subset $A$ in ${\cal A}_r$. Since $\rho >0$ can be taken
arbitrary small, then there exists $\rho _0>0$ such that for each
$0<\rho <\rho _0$ there is accomplished the inequality $C_2\rho
^{n-1}<\epsilon $, consequently,
\par $\sup_{p\in C_R^{\rho }, Re (p)\ge a} \| {\hat F}(p) \| \le
C' C_1\exp ((C_2\rho ^{n-1}-\epsilon )R)\le C' C_1\exp (-\delta R)$ \\
in view of the condition imposed on $F$, where $C$ is the positive
constant for the given $F$, $\delta = \epsilon - C_2\rho ^{n-1}
>0$. With this the length of the path of integration does not exceed
$2\pi R$ and $\lim_{R\to \infty } C' C_1 2\pi R \exp (-\delta R)=0$.
The function $F(p)$ is continuous by $p$, hence it is integrable
along each rectifiable curve in the domain $W$ of the half space $\{
p\in {\cal A}_r: Re (p)>s_0 \} $. \par If $F(p)$ is holomorphic,
then in view of Theorem 2.11 \cite{luoyst,lusmfcdv,luhcnfcdv} \\
$\int_{\gamma _n} F(p)\exp ( -u(p,t;q)) dp$ is independent from the
type of the curve and it is defined only by the initial and final
its points. If $V(\gamma _n)\le {C'}_V R_n$ for each $n$, then it is
sufficient to prove the statement of this Lemma for each subsequence
$R_{n(k)}$ with $R_{n(k+1)}\ge R_{n(k)}+1$ for each $k\in \bf N$.
Denote for the simplicity such subsequence by $R_n$. Each
rectifiable curve can be approximated by the converging sequence of
rectifiable of polygonal type line composed of arcs of circles. If a
curve is displayed on the sphere, then these circles can be taken
with the common center with the sphere. Condition $(2)$ in each
plane ${\bf R}\oplus N\bf R$, where $N\in {\cal A}_r$, $Re (N)=0$,
$|N|=1$, is accomplished, moreover, uniformly relative to a
directrix $N$ and it can be accomplished a diffeomorphism $g$ in
${\cal A}_r$, such that $g(W)=W$, $g(S_{R_n})=S_{R_n}$ for each
$n\in \bf N$, and an image of a $C^1$ curve from $W$ is an arc of a
circle, since $0<R_n+1<R_{n+1}$ for each $n\in \bf N$ and
$\lim_{n\to \infty }R_n=\infty $. \par The function $(F,\gamma
)\mapsto \int_{\gamma }F(p)dp$ is continuous from $C^0(V,{\cal
A}_r)\times \Gamma $ into ${\cal A}_r$, where $V$ is the compact
domain in ${\cal A}_r$, $\Gamma $ is the family of rectifiable
curves in $V$ with the metric $\rho (v,w):= \max (\sup_{z\in
v}\inf_{\zeta \in w}|z-\zeta |, \sup_{z\in w}\inf_{\zeta \in
v}|z-\zeta |)$ (see Theorem 2.7 \cite{luoyst,lusmfcdv,luhcnfcdv}).
The space $C^1$ of all continuously differentiable functions of the
real variable is dense in the space of continuous functions $C^0$ in
the compact-open topology in the case of a finite number of
variables. In addition a rectifiable curve is an uniform limit of
$C^1$ curves, since each rectifiable curve is continuous. Therefore,
consider $\gamma _n=\psi _n\cap \{ p\in {\cal A}_r: Re (p)>a \} \cap
W$, where $\psi _n$ is a curve in $S_{R_n}$ corresponding to $\gamma
_n$. Consequently,
\par $\lim_{n\to \infty }\int_{\gamma _n} F(p)\exp
( -u(p,t;q)) dp=0$, \\
since this is accomplished for $\gamma _n=\pi _n\cap C_{R_n}$ and
hence for general $\gamma _n$ with the same initial and final
points, where $\pi _n$ are two dimensional over $\bf R$ planes in
${\cal A}_r$.
\par The continuation of the {\bf Proof} of Theorem 8.
In view of Lemma 9
\par $ |\int_{\psi _R} F(p)(p-\eta )^{-1}dp|\le u(R)\pi R/(R-|\eta
|)$, \\
where $0<u(R)<\infty $ and there exists $\lim_{R\to \infty }
u(R)=0$, while $\psi _R$ is the arc of the circle $|p|=R$ in the
plane ${\bf R}\oplus S{\bf R}$ with $Re (p)>a$, consequently,
\par $ \lim_{R\to \infty } \int_{\psi _R} F(p)(p-\eta )^{-1}dp=0$, \\
since $u(R)\le u_0\exp (-\delta R)$ for each $R>R_0$, where
$u_0=const >0$. \par Then the straight line $a+S\theta $ with
$\theta \in \bf R$ can be substituted by the closed contour $\phi
_R$ composed from $\psi _R$ and the segment $[a+Sb,a-Sb]$ passed
from the above to the bottom. Thus,
\par $\int_0^{\infty }f(t)\exp (- \eta t)dt = (2\pi )^{-1}
{\tilde S}\int_{\phi _R}F(p)(p-\eta )^{-1}dp,$ \\
where the sign in front of the integral is changed due to the change
of the pass direction of the loop $\phi _R$. Recall, that in the
case of the  Cayley-Dickson algebra ${\cal A}_r$ the residue of a
function is the operator $\bf R$-homogeneous and ${\cal
A}_r$-additive by the argument $L\in {\cal A}_r$ with $Re (L)=0$,
where the residue is naturally dependent on a function and a point.
In the domain $ \{ p\in {\cal A}_r: Re (p)\ge a, |p| \le R \} $ the
analytic function $F(p)$ has only one point of singularity $p=\eta
$, which is the pole of the first order with the residue $res (\eta
; (p-\eta )^{-1}F(p)).L = L F (\eta )$ for each $L\in {\cal A}_r$
with $Re (L)=0$, consequently,
\par $\int_0^{\infty }  f(t)\exp (-\eta t)dt = F(\eta )$, since
$L=S$ in the given case and $S{\tilde S}=1$.
\par For $t<0$ in view of the aforementioned ${\cal A}_r$ Lemma 9
we get, that \par $\lim_{R\to \infty }\int_{\phi _R}
F(p)e^{u(p,t)}dp=0$, \\
since $Re (p)=a>0$, consequently, the straight line $a+S\theta $,
$\theta \in \bf R$, can be substituted by the loop $\phi _R$ as
above. Then for $t<0$ we get:
\par $f(t) = (2\pi )^{-1}\int_{\phi _R}F(p)e^{u(p,t)}dp=0$, \\
since $F(p)$ is analytic by $p$ together with $e^{u(p,t)}$ in the
interior of the domain $\{ p: p\in {\cal A}_r; |p|\le R', Re (p)>s_0
\} $, $a>s_0$, $0<R<R'\le \infty $. Then the condition 2 for the
original is accomplished. On the other hand,
\par $|f(t)| \le (2\pi )^{-1} e^{at} \int_{-\infty }^{\infty }
|F(a+S\theta )|d\theta = Ce^{at}$, \\
where $C=(2\pi )^{-1}\int_{-\infty }^{\infty }|F(a+S\theta )|d\theta
<\infty $, consequently, Condition (3) is satisfied. As well as
$f(t)$ is continuous, since the function $F(p)$ in the integral is
continuous satisfying Conditions $(i,ii)$ and
\par $\lim_{R\to \infty } \int_{\gamma (\theta ): |\theta |\ge R }
F(p)dp=0$. Moreover, the integral \par $\int_{a-S\infty }^{a+S\infty
}\mbox{ }_NF(p) [\partial \exp
(u(p,t))/\partial t]dp$ \\
converges due to Conditions $(i,ii)$ and the proof above,
consequently, the function $f(t)$ is differentiable and hence
satisfies the Lipschitz condition.

\par {\bf 10. Note.} In Theorem 8 Condition (i) can be replaced on
\par $\lim_{n\to \infty }\sup_{p\in S_{R(n)}} \| {\hat F}(p) \|
=0$ \\
by the sequence $S_{R_n} := \{ z\in {\cal A}_r: |z| =R_n, Re (z)>s_0
\} $, where $R_n<R_{n+1}$ for each $n$, $\lim_{n\to \infty }
R_n=\infty $, since this leads to the accomplishment of the ${\cal
A}_r$ analog of the Jordan Lemma for each $r\ge 2$ (see also Note 36
\cite{lusmtslt,luhcnlt}). But in Theorem 8 itself it is essential
the alternativity of the algebra, that is, in it in general are
possible only $r=2$ or $r=3$ for $f({\bf R})\subset {\cal A}_r$.

\par {\bf 11. Definition.} Consider function-originals, satisfying
conditions $(1-3)$ below:
\par $(1)$ $f(t)$ satisfies the H\"older condition: $|f(t+h)-f(t)|
\le A |h|^{\alpha }$ for each $|h|<\delta $ (where $0<\alpha \le 1$,
$A=const >0$, $\delta >0$ are constants for a given $t$) everywhere
on $\bf R$ may be besides points of discontinuity of the first kind.
On each finite interval in $\bf R$ a function $f$ may have only a
finite number of points of discontinuity and of the first kind only.
\par $(2)$ $|f(t)| < C_1 \exp (-s_1t)$ for each $t<0$,
where $C_1=const >0$, $s_1=s_1(f)=const \in \bf R$.
\par $(3)$ $|f(t)|<C_2 \exp (s_0t)$ for each $t\ge 0$, that is,
$f(t)$ is growing not faster, than the exponential function, where
$C_2=const>0$, $s_0=s_0(f)\in \bf R$.
\par The two-sided Laplace transformation over the Cayley-Dickson
algebras ${\cal A}_r$ with $2\le r\le 3$ is defined by the formula:
\par $(4)$ ${\cal F}^s(f;p) :=
\int_{-\infty }^{\infty }f(t)\exp ( -u(p,t;q))dt$ \\
for all numbers $p\in {\cal A}_r$, for which the integral exists,
where $q\in {\cal A}_r$ is a parameter, either $u(p,t;q)=pt+q$ or
$u(p,t;q)=E(pt+q)$ (see Definition 2.2.1). Denote for short ${\cal
F}^s(f;p)$ through $F^s(p)$. For a basis of generators $\{
N_0,...,N_{2^r-1} \} $ in ${\cal A}_r$ we shall write in more
details $\mbox{ }_N{\cal F}^s(f;p)$ or $\mbox{ }_NF^s_u(p)$ in the
case of necessity.
\par {\bf 12. Note.} Naturally, that the two-sided Laplace integral
can be considered as the sum of two one-sided integrals \\ $(1)$
$\int_{-\infty }^{\infty }f(t)\exp (-u(p,t;q))dt= \int_{-\infty
}^0f(t)\exp (-u(p,t;q))dt +
\int_0^{\infty }f(t)\exp (-u(p,t;q)dt$ \\
$=\int_0^{\infty }f(-t)\exp (-u(p,-t;q)))dt +
\int_0^{\infty }f(t)\exp (-u(p,t;q))dt$. \\
The second integral converges for $Re (p)>s_0$. Since $u(p, -t;q) =
u(-p,t;q)$, then the first integral converges for $Re (-p)> -s_1$,
that is, for $Re (p)< s_1$. Then there is a region of convergence
$s_0 < Re (p) < s_1$ of the two-sided Laplace integral. For
$s_1=s_0$ the region of convergence reduces to the vertical
hyperplane in ${\cal A}_r$ over $\bf R$. For $s_1<s_0$ there is no
any common domain of convergence and $f(t)$ can not be transformed
with the help of the two-sided transformation 1(4).
\par {\bf 13. Example.} ${\cal F}^s(\exp (-\alpha t^2); p)=
\int_{-\infty }^{\infty } \exp (-\alpha t^2-pt)dt= (\pi/\alpha
)^{1/2}\exp (p^2/(4\alpha ))$, where $\alpha >0$, since
$\int_{-\infty }^{\infty }\exp (-t^2)dt=(\pi )^{1/2}$. For
comparison the one-sided Laplace transformation gives:
\par ${\cal F}(\exp (-\alpha t^2)Ch_{[0,\infty )}; p)=
\int_0^{\infty } \exp (-\alpha t^2-pt)dt$ \par  $= (\alpha
)^{-1/2}\exp (p^2/(4\alpha ))\int_{p/(2(\alpha )^{1/2})}^{\infty
}\exp (-t^2)dt$ \\ $=2^{-1}(\pi /\alpha )^{1/2}\exp (p^2/(4\alpha ))
Erf (p/(2(\alpha )^{1/2}))$ \\ (see also \cite{lusmtslt,luhcnlt},
where numerous examples of calculations of noncommutative Laplace
transformations of such and more general type over $\bf H$ and $\bf
O$ and their applications to super-differential equation were
given).

\par The application of Theorem 4 to
$\int_0^{\infty }f(-t)\exp (- u(-p,t;q))dt$ and \\ $\int_0^{\infty
}f(t)\exp (-u(p,t;q))dt$ gives.
\par {\bf 14. Theorem.} {\it If an original $f(t)$ satisfies
Conditions 11$(1-3)$, and moreover, $s_0<s_1$, then its image ${\cal
F}^s(f;p)$ is holomorphic by $p$ in the domain $\{ z\in {\cal A}_r:
s_0< Re (z)<s_1 \} $, where $2\le r\le 3$.}

\par {\bf 15. Examples.} 1. There may be cases, when a domain of
convergence for a sum is greater, then for each additive. For
example, ${\cal F}^s(\exp (at) U(t);p)=(p-a)^{-1}$, also ${\cal
F}^s((\exp (at) -1) U(t);p)=ap^{-1}(p-a)^{-1}$ for $Re (p)>a$ in
both cases, when $a\in \bf R$, $U(t):=1$ for $t>0$, $U(0)=1/2$,
while $U(t)=0$ for $t<0$. But $(p-a)^{-1} -
ap^{-1}(p-a)^{-1}=p^{-1}$ and ${\cal F}^s(U(t);p)=p^{-1}$ for each
$Re (p)>0$.
\par It is necessary to note, that the two-sided Laplace
transformation of the function $t^n$ does not exist, but the
one-sided transformation was elucidated in examples 2.30.1 and 2.33
\cite{lusmtslt,luhcnlt}.
\par 2. ${\cal F}^s(\exp (-\alpha |t|)/2;p)=\alpha
(\alpha ^2-p^2)^{-1}$ in the domain $-\alpha <Re (p)<\alpha $ for
$\alpha >0$.
\par 3. Consider the two-sided transformation
\par ${\cal F}^s((e^t+1)^{-1};p) = \int_{-\infty }^{\infty }
(e^t+1)^{-1}\exp (-pt)dt$ \\ in the domain $-1<Re (p)<0$. Make the
substitution $v=(e^t+1)^{-1}$, then the integral reduces to the
Euler integral of the first kind $\int_0^1v^p(1-v)^{-p-1}dv = -\pi
/\sin (\pi p)$ (see Proposition 4.6, Definition 4.14 and Theorem
4.17 in \cite{lusmfcdv}).

\par {\bf 16. Theorem.} {\it If a function $f(t)$ is an original
such that \par $\mbox{ }_N{\cal F}^s(f;p;q) := \sum_{j=0}^{2^r-1}
\mbox{ }_NF^s_{u,j}(p;q)N_j$ is its image, where a function $f$ is
written in the form \par $f(t) = \sum_{j=0}^{2^r-1} f_j(t)N_j$,
$f_j: {\bf R}\to \bf R$ for each $j=0,1,...,2^r-1$, $f({\bf
R})\subset {\cal A}_r$ for $2\le r\le 3$,
\par $\mbox{ }_NF^s_{u,j}(p;q) := \int_{-\infty }^{\infty
}f_j(t)\exp (-u(p,t;q))dt$. Then at each point $t$, where $f(t)$
satisfies the H\"older condition there is true the equality:
\par $(i)$ $f(t) = (2\pi N_1)^{-1} Re (S{\tilde N}_1) \sum_{j=0}^{2^r-1}
(\int_{a-S\infty }^{a+S\infty }\mbox{ }_NF^s_{u,j}(p;q)
\exp (u(p,t;q))dp)N_j$ \\
in the domain $s_0(f)< Re (p) <s_1(f)$, where either $u(p,t;q
)=pt+q$ with $S=N_1$ and $Im (q)=0$, or $u(p,t;q )=E(pt+q)$ and the
integral is taken along the straight line $p(\tau )=a+S\tau \in
{\cal A}_r$, $\tau \in \bf R$, $S\in {\cal A}_r$, $Re (S)=0$,
$|S|=1$, $Re (S{\tilde N}_1)\ne 0$ is non-zero, while the integral
is understood in the sense of the principal value.}
\par {\bf Proof.} The two-sided transformation in the basis of generators
\\ $N = \{ N_0,N_1,...,N_{2^r-1} \} $ can be written in the form
\par ${\cal F}^s(f;p;q) := \int_{-\infty }^{\infty }f(t)\exp
(-u(p,t;q))dt = {\cal F}^s(fU(t);p;q) + {\cal F}^s(f(1-U)(t);p;q) $,
\\ where the index $N$ is omitted,
$U(t)=1$ for $t>0$, $U(0)=1/2$, $U(t)=0$ for $t<0$, also
\par ${\cal F}^s(f(1-U)(t);p;q) = \int_0^{\infty
}f(-t)U(t)\exp (-u(-p,t;q))dt,$ \\  since $u(p,-t;q)= u(-p,t;q)$,
where $|f(-t)|\le C_1\exp (s_1t)$ for each $t>0$. The common domain
of the existence \\ $\int_0^{\infty }f(-t)U(t)\exp (-u(-p,t;q))dt$
and $\int_0^{\infty }f(t)U(t)\exp (-u(p,t;q))dt$ is $s_0(f)<Re
(p)<s_1(f)$, since the inequality $Re (-p)> - s_1(f)$ is equivalent
to the inequality $Re (p)<s_1(f)$. Then the application of Theorem 5
twice to $f(t)U(t)$ and to $f(-t)U(t)$ gives the statement of this
theorem.
\par {\bf 17. Theorem.} {\it If a function $\mbox{ }_NF^s_u(p)$
is analytic by the variable $p\in {\cal A}_r$ in the domain $W := \{
p\in {\cal A}_r: s_0< Re (p) <s_1 \} $, where $2\le r\le 3$, $f({\bf
R})\subset {\cal A}_r$, either $u(p,t)=pt$ or $u(p,t) := E(pt)$. Let
also $\mbox{ }_NF^s_u(p)$ can be written in the form $\mbox{
}_NF^s_u(p)=\mbox{ }_NF^{s,0}_u(p) + \mbox{ }_NF^{s,1}_u(p)$, where
$\mbox{ }_NF^{s,0}_u(p)$ is holomorphic by $p$ in the domain $s_0<Re
(p)$, also $\mbox{ }_NF^{s,1}_u(p)$ is holomorphic by $p$ in the
domain $Re (p)<s_1$, $S\in {\cal I}_r$, $|S|=1$, moreover, for each
$a>s_0$ and $b<s_1$ there exists constants $C_a>0$, $C_b>0$ and
$\epsilon _a
>0$ and $\epsilon _b>0$ such that
\par $(i)$ $|\mbox{ }_NF^{s,0}_u(p)|\le C_a\exp (-\epsilon _a |p|)$
for each $p\in {\cal A}_r$ with $Re (p)\ge a$,
\par $(ii)$ $|\mbox{ }_NF^{s,1}_u(p)|\le C_b\exp (-\epsilon _b |p|)$
for each $p\in {\cal A}_r$ with $Re (p)\le b$, where $s_0$ and $s_1$
are fixed, also the integral
\par $(iii)$ $\int_{w-S\infty }^{w+S\infty }\mbox{
}_NF^{s,k}_u(p)dp$ \\ converges absolutely for $k=0$ and $k=1$ for
$s_0<w<s_1$. Then $\mbox{ }_NF^s_u(p)$ is the image of the function
\par $(iv)$ $f(t)=(2\pi )^{-1}{\tilde S}\int_{w-S\infty
}^{w+S\infty }\mbox{ }_NF^s_u(p)\exp (u(p,t))dp$.}
\par {\bf Proof.} For the function
$\mbox{ }_NF^{s,1}_u(p)$ we consider the substitution of the
variable $p=-g$, $-s_1<Re (g)$. In view of Theorem 8 there exist
originals $f^0$ and $f^1$ for functions $\mbox{ }_NF^{s,0}_u(p)$ and
$\mbox{ }_NF^{s,1}_u(p)$ while a choice of $w\in \bf R$ in the
common domain $s_0<Re (p)<s_1$, that is, $s_0<w<s_1$. At the same
time the supports of the functions $f^0$ and $f^1$ are contained in
$[0,\infty )$ and $(-\infty ,0]$ respectively. Then $f=f^0+f^1$ is
the original for $\mbox{ }_NF^s_u(p)$ while $q=0$, since
\par $f(t)=f^0(t)+f^1(t)=(2\pi )^{-1}{\tilde S}\int_{w-S\infty
}^{w+S\infty }\mbox{ }_NF^{s,0}_u(p)\exp (u(p,t))dp+ $\\  $(2\pi
)^{-1}{\tilde S}\int_{w-S\infty }^{w+S\infty }\mbox{
}_NF^{s,1}_u(p)\exp (u(p,t))dp= (2\pi )^{-1}{\tilde
S}\int_{w-S\infty }^{w+S\infty }\mbox{ }_NF^s_u(p)\exp (u(p,t))dp$
\\ due to the distributivity of the multiplication in the algebra
${\cal A}_r$.

\par {\bf 18. Note.} While the definition of the one- and two-sided
Laplace transformations over the Cayley-Dickson algebras above the
Riemann integral of the real variable was used, while for the
inverse transformation the noncommutative integral along paths over
${\cal A}_r$ from the works \cite{luoyst,lusmfcdv}. It can be
considered also a generalization of the direct transformation with
the Riemann-Stieltjes integral as the starting point. For a function
$\alpha (t)$ with values in ${\cal A}_r$ of the variable $t\in \bf
R$ such that $\alpha (t)$ has a bounded variation on each finite
segment $[a,b]\subset \bf R$, we consider the Stieltjes integral
\par $\int_{-\infty }^{\infty }(d\alpha (t))\exp (-u(p,t;q)):=$
\\ $\lim_{b\to \infty } \int_0^b(d\alpha (t))\exp (-u(p,t;q)) +
\lim_{b\to \infty } \int_{-b}^0(d\alpha (t))\exp (-u(p,t;q)) $,
where
\par $\int_a^b(d\alpha (t))f(t)=\sum_{v,w}(\int_a^bf_v(t)d\alpha
_w(t)) (i_wi_v)$, \\
$i_0,i_1,...,i_{2^r-1}$ are generators of the Cayley-Dickson algebra
${\cal A}_r$, $f_v$ and $\alpha _w$ are real-valued functions such
that $\alpha =\sum_w\alpha _wi_w$ and $f=\sum_vf_vi_v$, also
$\int_a^b f_v(t)d\alpha _w(t)$ is the usual Stieltjes integral over
the field of real numbers on a finite segment $[a,b]$. Under
imposing the condition $\alpha (t)\exp (-p_0t)|_{-\infty }^{\infty
}=0$ the integration by parts gives the relation
\par $\int_{-\infty }^{\infty }(d\alpha (t))\exp (-u(p,t;q))=
 - \int_{-\infty }^{\infty }\alpha (t)d[\exp (-u(p,t;q))] $. In view
of the associativity of $\bf H$ and the alternativity of $\bf O$ it
gives:
\par $(1)$ $\int_{-\infty }^{\infty }(d\alpha (t))\exp ( -pt)=
{\cal F}^s(\alpha (t)p;p)$. \\
For $u(p,t;q)=E(pt+q)$ there is the formula
\par $(2)$  $\int_{-\infty }^{\infty }(d\alpha (t))\exp (-u(p,t;q))
=p_0{\cal F}^s(\alpha (t); p;q) $ \\
$+p_1{\cal F}^s(\alpha (t); p;q - i_1\pi /2)+...+p_{2^r-1}{\cal
F}^s(\alpha (t); p;q- i_{2^r-1}\pi /2)$ \\
over ${\cal A}_r$ with $2\le r\le 3$.

\par Thus, the Laplace transformation over ${\cal A}_r$
can be spread on a more general class of originals. For this it is
used instead of an ordinary notion of convergence of improper
integrals their convergence by Cesaro of order $p>0$:
\par $(C,p) \int_0^{\infty }f(t)dt := \lim_{b\to \infty }
\int_0^bf(t)(1-t/b)^pdt$. \\
If this integral converges by Cesaro for some $p>0$, then it
converges for each $q>p$, moreover,
\par $(C,p) \int_0^{\infty }f(t)dt = (C,q)\int_0^{\infty }f(t)dt$. \\
That is, with the growth or the order $p$ a family of functions
enlarges for which an improper integral converges. The limit case of
the limit by Cesaro is the Cauchy limit:
\par $(C) \int_0^{\infty }f(t)dt := \lim_{\epsilon \to +0}
\int_0^{\infty }f(t)\exp (-\epsilon t)dt$. \\
For two-sided integrals convergence of improper integrals by Cesaro
of order $p$ is defined by the equality:
\par $(C,p) \int_{-\infty }^{\infty }f(t)dt := \lim_{b\to \infty }
\int_{-b}^bf(t)(1-|t|/b)^pdt$ \\
and by Cauchy:
\par $(C) \int_{-\infty }^{\infty }f(t)dt := \lim_{\epsilon \to +0}
\int_{-b}^bf(t)\exp (-\epsilon |t|)dt$, \\
when these limits exist.

\par The noncommutative Mellin transformation is based on the two sided
transformation of the Laplace type over Cayley-Dickson algebras,
which was presented above.
\par {\bf 19. Remark.} If $f$ is an original function of the two-sided
Laplace transformation over the Cayley-Dickson algebra ${\cal A}_r$
and $g(\tau )=f(\ln \tau )$ for each $0<\tau <\infty $, then
Conditions 11$(1-3)$ for $f$ are equivalent to the following
Conditions M(1-3):
\par $M(1)$ $g(\tau )$ satisfies the H\"older condition:
$|g(\tau +h)-g(\tau )| \le A |h|^{\alpha }$ for each $|h|<\delta $
(where $0<\alpha \le 1$, $A=const >0$, $\delta >0$ are constants for
a given $\tau $) everywhere on $\bf R$ may be besides points of
discontinuity of the first kind. On each finite segment $[a,b]$ in
$(0, \infty )$ a function $g$ may have only a finite number of
points of discontinuity and of the first kind only.
\par $M(2)$ $|g(\tau )| < C_1 \tau ^{s_0}$ for each $0<\tau <1$,
where $C_1=const >0$, $s_0=s_0(g)=const \in \bf R$.
\par $M(3)$ $|g(\tau )|<C_2 \tau ^{-s_1}$ for each $\tau \ge 1$, that is,
$g(\tau )$ is growing not faster, than the power function, where
$C_2=const>0$, $s_1=s_1(g)\in \bf R$.
\par This is because of the fact that the logarithmic function
$\ln : (0,\infty )\to (-\infty ,\infty )$ is the diffeomorphism.
\par {\bf 20. Definition.} In the two-sided integral transformation
of the Laplace type substitute variables $p$ on $-p$ and $t$ on
$\tau =e^t$ and $q$ on $-q$, then the formula takes the form:
\par $(1)$ ${\cal M}(g;p) := \int_0^{\infty }
f(\ln \tau )\exp (- u(-p,\ln \tau ;-q))\tau ^{-1}d\tau $, \\
where $f$ is an original function, $g(\tau )=f(\ln \tau )$ for each
$0<\tau <\infty $ (see also Definition 11). \par For a specified
basis $ \{ N_0, N_1,..., N_{2^r-1} \} $ of generators of the
Cayley-Dickson algebra ${\cal A}_r$ we can write the notation in
more details $\mbox{ }_N{\cal M}_u(g;p;q)$ if necessary.
\par {\bf 21. Theorem.} {\it If an original function $g(\tau )$ satisfies
Conditions 20.M(1-3), where $s_0<s_1$, then its image ${\cal
M}(g;p;q)$ is holomorphic by $p$ in the domain $\{ z\in {\cal A}_r:
s_0< Re (z)<s_1 \} $, where $2\le r\le 3$.}
\par {\bf Proof.} The application of Theorem 4 to $g(\tau )=
f(\ln \tau )$ gives the statement of this theorem.

\par {\bf 22. Theorem.} {\it Let $g(\tau )$ be an original function
such that \par $\mbox{ }_N{\cal M}(g;p;q) := \sum_{j=0}^{2^r-1}
\mbox{ }_NG_{u,j}(p;q)N_j$ be its image, where a function $g$ is
written in the form \par $g(\tau ) = \sum_{j=0}^{2^r-1} g_j(\tau
)N_j$, $g_j: (0,\infty )\to \bf R$ for each $j=0,1,...,2^r-1$,
$g((0,\infty ))\subset {\cal A}_r$ for $2\le r\le 3$,
\par $\mbox{ }_NG_{u,j}(p;q) := \int_0^{\infty
}g_j(\tau )\exp (-u(-p,\ln \tau ; - q))\tau ^{-1}d\tau $. Then at
each point $\tau $, where $g(\tau )$ satisfies the H\"older
condition the equality is accomplished:
\par $(i)$ $g(\tau ) = (2\pi N_1)^{-1} Re (S{\tilde N}_1)
\sum_{j=0}^{2^r-1} (\int_{a-S\infty }^{a+S\infty }\mbox{
}_NG_{u,j}(p;q) \exp (u(-p,\ln \tau ; - q))dp)N_j$ \\
in the domain $s_0(g)< Re (p) <s_1(g)$, where either $u(p,t;q)=pt+q$
with $S=N_1$ for $Im (q)=0$, or $u(p,t;q)=E(pt+q)$ (see \S 2.2.1)
and the integral is taken along the straight line $p(\theta
)=a+S\theta \in {\cal A}_r$, $\theta \in \bf R$, $S\in {\cal A}_r$,
$Re (S)=0$, $|S|=1$, $Re (S{\tilde N}_1)\ne 0$ is non-zero, while
the integral is understood in the sense of the principal value.}
\par {\bf Proof.} Putting $t=\ln \tau $ and substituting $p$ on $-p$
and applying Theorem 16 we get the statement of this theorem.

\par {\bf 23. Theorem.} {\it If a function $\mbox{
}_NG_u(p)$ is analytic by the variable $p\in {\cal A}_r$ in the
domain $W := \{ p\in {\cal A}_r: s_0< Re (p) <s_1 \} $, where $2\le
r\le 3$, $g((0,\infty ))\subset {\cal A}_r$, either $u(p,t)=pt$ or
$u(p,t) :=E(pt)$. Let also $\mbox{ }_NG_u(p)$ can be written in the
form $\mbox{ }_NG_u(p)=\mbox{ }_NG^0_u(p) + \mbox{ }_NG^1_u(p)$,
where $\mbox{ }_NG^0_u(p)$ is holomorphic by $p$ in the domain
$s_0<Re (p)$, also $\mbox{ }_NG^1_u(p)$ is holomorphic by $p$ in the
domain $Re (p)<s_1$, $S\in {\cal I}_r$, $|S|=1$, moreover, for each
$a>s_0$ and $b<s_1$ there exists constants $C_a>0$, $C_b>0$ and
$\epsilon _a
>0$ and $\epsilon _b>0$ such that
\par $(1)$ $|\mbox{ }_NG^0_u(p)|\le C_a\exp (-\epsilon _a |p|)$
for each $p\in {\cal A}_r$ with $Re (p)\ge a$,
\par $(2)$ $|\mbox{ }_NG^1_u(p)|\le C_b\exp (-\epsilon _b |p|)$
for each $p\in {\cal A}_r$ with $Re (p)\le b$, where $s_0$ and $s_1$
are fixed, while the integral
\par $(3)$ $\int_{w-S\infty }^{w+S\infty }\mbox{
}_NG^k_u(p)dp$ \\ converges absolutely for $k=0$ and $k=1$ for
$s_0<w<s_1$, then $\mbox{ }_NG_u(p)$ is the image of the function
\par $(4)$ $g(\tau )=(2\pi )^{-1}{\tilde S}\int_{w-S\infty
}^{w+S\infty }\mbox{ }_NG_u(p)\exp (u(-p,\ln \tau ))dp$.}
\par {\bf Proof.} The change of the variable $p$ on $-p$ and the
substitution $t=\ln \tau $ for $\tau >0$ with the help of Theorem 17
gives the statement of this Theorem.

\par {\bf 24. Theorem.} {\it Let $f$ be an original function
from Definition either 1 or 11 or 20. Suppose that $F$ is its image
function of the noncommutative either Laplace or two-sided Laplace
or Mellin transformation for either $u(p,t)=pt$ or $u(p,t)=E(pt)$ in
the domain $V:= \{ z\in {\cal A}_b: s_0< Re (z) <s_1 \} $, $b=2$ or
$b=3$, where $s_1=\infty $ for the noncommutative one-sided Laplace
transformation. Then $F$ is either $(1,b)$-quasi-regular or
$(1,b)$-quasi-regular in spherical ${\cal A}_b$-coordinates
respectively in $V$ with $y_0=0$ if and only if its original is real
$f(t)\in \bf R$ for each continuity point $t$ of $f$ either in
$[0,\infty )$ or $\bf R$ or $(0,\infty )$ respectively.}
\par {\bf Proof.} Since the Mellin transformation is obtained from
the two-sided Laplace transformation with the help of smooth change
of real variables and the one-sided Laplace transformation is the
particular case of that of two-sided, then it is sufficient to prove
this theorem for the two-sided noncommutative Laplace
transformation. Thus consider $F(p)= \int_{-\infty }^{\infty
}f(t)\exp (-u(p,t))dt$, since $q=0$, $y_0=0$ by the conditions of
this theorem. We have that ${\hat R}_{p,x}=R_{w(p),w(x)}$, where $w$
is a pseudo-conformal diffeomorphism of $V$, $R_{p,x}$ is given in
Examples 2.2 and 2.9.5.2. To each automorphism ${\hat R}_{p,x}$ the
operator belonging to the Lie group $SO_{\bf R}(2^b,{\bf R})$ on the
real shadow ${\bf R}^{2^b}$ corresponds. Therefore,
\par $(1)$ ${\hat R}_{p,y} F_u(y) = \int_{-\infty }^{\infty }[
R_{w(p),w(y)}f(t)]\exp (- (R_{w(p),w(y)} y)t) dt$ for $u=pt$, for
each $p\in V$ and every $y\in V\cap \bf C$ such that $Re (p)=Re (y)$
and $R_{w(p), w(y)} y=p$, since ${\hat R}_{p,y}|_{\bf R}=id$. Then
\par $(2)$ ${\hat R}_{E(p),E(y)} F_u(y) = \int_{-\infty }^{\infty }
[R_{w(E(p)), w(E(y))} f(t)] \exp (- R_{w(E(tp)), w(E(ty))} E(ty)) dt$ \\
for $u(p,t) =E(pt)$ respectively due to Formula 9.5.2(2) for each
$p\in V$ and every $y\in V\cap \bf C$ such that $Re (E(p))=Re
(E(y))$ and $R_{w(E(p)), w(E(y))} E(y)=E(p)$. Thus
\par $(3)$ $F_u(p)={\hat R}_{p,y}F_u(y)=
\int_{-\infty }^{\infty }f(t)\exp
(- pt) dt= \int_{-\infty }^{\infty }[{\hat R}_{p,y}f(t)]\exp (- pt)
dt$
\\
for $u=pt$, for each $p\in V$ and every $y\in V\cap \bf C$ such that
$Re (p)=Re (y)$ and $R_{w(p), w(y)} y=p$. On the other hand,
\par $(4)$ $R_{w(E(p)), w(E(y))}\exp (E(ty))=R_{w(E(tp)), w(E(ty))}
\exp (E(ty))=\exp (E(tp))$ \\
for $u(p,t) =E(pt)$, for each $Re (E(p))=Re (E(y))$ with
$R_{w(E(p)), w(E(y))} E(y)=E(p)$, since $R_{z,x}(tx)=tR_{z,x}x$ for
each $t\in \bf R$ and $E(y)=y$ for each $y\in \bf C$. Consequently,
\\ $(5)$ $F_u(p)={\hat R}_{E(p),E(y)} F_u(y)=
\int_{-\infty }^{\infty }f(t)\exp (- E(tp))dt =\int_{-\infty
}^{\infty } [{\hat R}_{E(p),E(y)}f(t)]\exp (- E(tp))dt$ \\
for $u(p,t)=E(pt)$. Particularly, $w=id$ can also be taken.
\par The two-sided Laplace transformation
is injective such that ${\cal F}^s(f_1;z)={\cal F}^s(f_2;z)$ for
each $z\in V$ if and only if $f_1(t)=f_2(t)$ at each point $t$ in
$\bf R$ where $f_1(t)$ and $f_2(t)$ are continuous (see Theorems 5,
7, 8, 16, 17, 22 and 23). Thus due to Formulas $(3,4)$ $F(z)$ is
either $(1,b)$-quasi-regular or $(1,b)$-quasi-regular in spherical
${\cal A}_b$-coordinates correspondingly if and only if either
${\hat R}_{p,y}f(t) = f(t)$ or ${\hat R}_{E(p),E(y)}f(t) = f(t)$
respectively for each $t\in \bf R$ a point of continuity of $f$ and
each $p\in V$ and every $y\in V\cap \bf C$ such that $Re (p)=Re (y)$
and either ${\hat R}_{p,y}(y)=p$ or ${\hat
R}_{E(p),E(y)}(E(y))=E(p)$ correspondingly. This means that $f(t)\in
\bf R$, since if $Im (s)\ne 0$ for some $s\in {\cal A}_b$, then
there exist $p\in V\setminus \bf C$ and $y\in V\cap \bf C$ such that
either $({\hat R}_{p,y} s) \ne s$ or $({\hat R}_{E(p),E(y)} E(s))
\ne E(s)$ respectively (see 2.1$(Q2-Q5)$). Since $f$ is continuous
besides points of discontinuity of the first kind, then using limits
from the left or from the right redefine $f$ at points of
discontinuity such that $f$ will be real everywhere on $\bf R$.
\par If $f$ is real-valued, then $F_u$ satisfies Conditions
2.1$(Q1-Q5)$ by the construction of $F_u$. Then $F_u$ satisfies
2.1$(Q6,Q7)$ as well due to Theorem 2.18, since the function
$e^{ap}=v(p)$ is pseudo-conformal for $a\ne 0$ on ${\cal A}_b$,
$p\in {\cal A}_b$.

\par {\bf 25. Theorem.} {\it Let suppositions of Theorem 24 be
satisfied for the noncommutative two-sided Laplace or Mellin
transformation. If $f$ is real-valued, then
\par $(1)$ $F_u({\tilde p})={\tilde F}_u(p)$ for $u(p,t) =pt$ or
\par $(1')$ $F_u(p_0-p_1i_1+p_2i_2+...+p_{2^b-1}i_{2^b-1})=
{\tilde F}_u(p)$ for $u(p,t) =E(pt)$ respectively for each $p\in V$.
Moreover, either $f(t)=f(-t)$ is even for each $t\in \bf R$ or
$f(t)=f(1/t)$ for each $t>0$ at each point of continuity of $f$ if
and only if its noncommutative two-sided Laplace or Mellin
transformation $F_u(p)$ for $u(p,t)=pt$ or $u(p,t)=E(pt)$ satisfies
the condition:
\par $(2)$ $F_u(-p)=F_u(p)$ for each $p\in V$ for both types of $u$.}
\par {\bf Proof.} If an original $f$ is real-valued, then
\par $[\int_{-\infty }^{\infty }f(t)\exp (-u(p,t))dt]^*=
\int_{-\infty }^{\infty }[\exp (-u(p,t))]^* [f(t)]^*dt$ \\
$= \int_{-\infty }^{\infty } f(t) \exp (-[u(p,t)]^*)dt$, \\ but
$[u(p,t)]^*=u(p^*,t)$ for $u=pt$ and $[u(p,t)]^*= u((p_0-
p_1i_1+p_2i_2+...+p_{2^b-1}i_{2^b-1}),t)$ for $u(p,t)=E(pt)$ (see
Formulas 2.2.1$(1,2)$ or 16.1$(3,5)$), where
$p=p_0+p_1i_1+...+p_{2^b-1}i_{2^b-1}$, $p_j\in \bf R$ for each
$j=0,...,2^b-1$. Therefore, either $(1)$ or $(1')$ respectively is
satisfied.
\par An original $f$ is even on $\bf R$ if and only if
\par $\int_{-\infty }^{\infty } f(t) \exp (-u(p,t))dt=
\int_{-\infty }^{\infty } f(-t)\exp (-u(p,t))dt$ \\
$= - \int_{ \infty }^{- \infty } f(t)\exp (-u(p, -t))dt =
\int_{-\infty }^{\infty } f(t) \exp (-u(-p,t))dt$ \\ for both
variants $u(p,t)=pt$ and $u(p,t)=E(pt)$, since $u(p,-t)=u(-p,t)$
while the two-sided Laplace transformation is injective such that
${\cal F}^s(f_1;z)={\cal F}^s(f_2;z)$ for each $z\in V$ if and only
if $f_1(t)=f_2(t)$ at each point $t$ in $\bf R$ where $f_1(t)$ and
$f_2(t)$ are continuous (see Theorems 16, 17). Consequently,
Condition $(2)$ is equivalent to $f(t)=f(-t)$ for each $t\in \bf R$
for the noncommutative two-sided Laplace transformation.
 \par Substituting $t$ on $\ln (\tau)$ and $p$ on $-p$ gives
that $(2)$ is equivalent to $f(t)=f(1/t)$ for each $t>0$ for the
noncommutative Mellin transformation due to Theorems 22 and 23.

\par {\bf 25.1. Proposition.} {\it Let $f$ be either a
$(1,b)$-quasi-regular or $(1,b)$-quasi-regular in spherical ${\cal
A}_b$-coordinates function on a domain $V$, $f(z)\ne 0$ for each
$z\in V$, where $2\le b\le 3$. Then $1/f(z)$ is either a $(1,b)$
quasi-regular or $(1,b)$-quasi-regular in spherical ${\cal
A}_b$-coordinates function respectively on $V$.}
\par {\bf Proof.} Take without loss of generality $y_0=0$.
Since ${\hat R}_{z,x}$ and ${\hat R}_{E(z),E(x)}$ are automorphisms
of ${\cal A}_b$, then $1/f$ or $1/f\circ E^{-1}$ respectively
satisfies Conditions 2.1$(Q1-Q6)$ on $V$ (see also Definition
2.2.1). Since $f$ is ${\cal A}_b$ holomorphic, then $1/f$ is also
${\cal A}_b$ holomorphic, $(\partial (1/f(z))/\partial {\tilde
z}).h=0$ for each $h\in {\cal A}_b$ and all $z\in V$ (see
\cite{luoyst,lusmfcdv,luhcnfcdv}). On the other hand,
$f(z)[1/f(z)]=1$ for each $z\in V$, hence \par $(1)$
$[\partial (1/f(x))/\partial x].h = - f(x)[(f'(x).h)(1/f(x))]$ \\
for each $h\in {\cal A}_b$ and every $x\in V$, since ${\bf O}={\cal
A}_3$ is alternative, ${\bf H}={\cal A}_2$ is associative, where
$f'(z).h=(\partial f(z)/\partial z).h$. Acting on both sides of
Equation $(1)$ by either ${\hat R}_{z,x}$ or ${\hat R}_{E(z),E(x)}$
gives $(Q7)$ for $1/f(z)$ or $1/f\circ E^{-1}$ respectively, since
$f(z)$ or $f\circ E^{-1}$ correspondingly satisfies $(Q1-Q7)$.

\par {\bf 26. Examples. 1.} Consider now the zeta function
on ${\cal A}_b$ (see Example 9.5.2). In view of Theorem 2.1
\cite{titchm} the zeta function $\zeta (s)$ has the holomorphic
extension in ${\bf C}\setminus \{ 1 \} $ with the pole at $s=1$ with
residue $1$, moreover, it satisfies the functional equation $\zeta
(s) = 2^s\pi ^{s-1}\sin (s\pi /2) \Gamma (1-s)\zeta (1-s)$.
\par Construct for $\zeta (s)$ $(1,b)$-quasi-conformal in
Examples 1 and 2 and $(1,b)$-quasi-conformal in spherical ${\cal
A}_b$-coordinates in Example 3 extensions in ${\cal A}_b\setminus \{
1 \} $. For this put $z=x+yM$, where $x, y\in \bf R$, $Re (M)=0$,
$|M|=1$, $r=1$, $y_0=0$. Then $z$ is obtained from $s=x+{\bf i}y$ by
the automorphism ${\hat R}_{z,s}$ such that ${\hat R}_{z,s}(i)=M$,
where ${\bf i}=i_1$. Then ${\bf R}\oplus M{\bf R} =: {\bf C}_M$ is
the subalgebra in ${\cal A}_b$ isomorphic with $\bf C$. Let $a$ and
$q$ be positive integers, $q>a$, $z\ne 1$, then
\par $(1)$ $\sum_{n=a+1}^q n^{-z}= (q^{1-z} - a^{1-z} )/(1-z)
- z \int_a^q (x- [x] -1/2)x^{-z-1}dx + (q^{-z} - a^{-z})/2$, \\
where $[x]$ denotes the greatest integer not exceeding $x$. For $Re
(z)=:\sigma >1$ and $a=1$ consider $q\to \infty $, then from Formula
$(1)$ we get
\par $(2)$ $\zeta (z)=z \int_1^{\infty } ([x]- x +1/2) x^{-z-1}
dx +1/(z-1) +1/2$. \\
The function $[x] -x+1/2$ is bounded, consequently, this integral
converges for $\sigma >0$ and uniformly converges in the domain
$\sigma >\delta $ in ${\cal A}_b$, where $\delta >0$ is the
constant. Therefore, this integral defines a holomorphic function of
$z$ $(1,b)$-quasi-regular for $\sigma >0$, $z\ne 1$, due to Theorem
24. The right hand side of Equation $(2)$ thus provides the ${\cal
A}_b$ holomorphic continuation of $\zeta (z)$ up to $\sigma =0$,
while there is a simple pole at $z=1$ with residue $1$.
\par For $0<\sigma <1$ Formula $(2)$ may be written as
\par $\zeta (z) = z\int_0^{\infty } ([x] -x)x^{-z-1}dx$,\\
since $\int_0^1([x]- x)x^{-z-1} dx= - \int_0^1 x^{-z} dx= 1/(z-1)$
and $z\int_1^{\infty }x^{-z-1}dx/2=1/2$. Consider $f(x)=[x]-x+1/2$,
$f_1(x)=\int_1^xf(y)dy$, then $f_1(y)$ is bounded, since
$\int_k^{k+1}f(y)dy=0$ for each integer $k$. Consequently,
$\int_{x_1}^{x_2}f(x)x^{-z-1}dx =f_1(x)x^{-z-1} |_{x_1}^{x_2} +
(z+1) \int_{x_1}^{x_2}f_1(x)x^{-z-2}dx$, which tends to zero as
$x_1\to \infty $ and $x_2\to \infty $, while $\sigma >-1$.
Therefore, the integral in $(2)$ is convergent for $\sigma >-1$,
hence $(2)$ gives the holomorphic continuation of $\zeta (z)$ for
$\sigma >-1$. Since $z\int_0^1([x]-x+1/2)x^{-z-1}dx=1/(z-1)+1/2$ for
$\sigma <0$. Hence
\par $(3)$ $\zeta (z)=z\int_0^{\infty }([x]-x+1/2)x^{-z-1}dx$
for $-1<\sigma <0$. In view of Proposition 2.9.1 and Theorem 24 and
Formulas $(2,3)$ and using the continuous extension from $\{ z\in
{\cal A}_b: -1<Re (z)<0$ $\mbox{or}$ $0<Re (z) \} $ the function
$\zeta (z)$ is $(1,b)$-quasi-regular in the domain $ \{ z\in {\cal
A}_b: -1<Re (z), z\ne 1 \} $.
\par Consider $\int_R^{\infty }\sin (2\pi nx) x^{-z-1}dx =$ \par
$[-\cos (2\pi nx )/(2\pi n x^{z+1})]|_R^{\infty } - (z+1)(2\pi
n)^{-1}\int_R^{\infty }\cos (2\pi nx)x^{-z-2}dx$
\par $= O(1/(nR^{\sigma +1})) +
O(n^{-1}\int_R^{\infty }x^{-\sigma -2}dx)=O(1/(n R^{\sigma +1}))$,
where $R>0$, consequently, $\lim_{R\to \infty }\sum_{n=1}^{\infty
}n^{-1}\int_R^{\infty }\sin (2\pi nx)x^{-z-1}dx=0$ for $-1<\sigma
<0$. Since there is the Fourier series expansion: $[x]-x+1/2
=\sum_{n=1}^{\infty }\sin (2\pi nx)(\pi n)^{-1}$ for non-integer
real $x$, then integrating in $(3)$ term by term series we obtain
\par $(4)$ $\zeta (z) = z(\pi )^{-1} \sum_{n=1}^{\infty } n^{-1}
\int_0^{\infty } \sin (2\pi nx)x^{-z-1}dx=$
\par $z\pi ^{-1} \sum_{n=1}^{\infty }(2\pi n)^zn^{-1}
\int_0^{\infty }\sin (y) y^{-z-1}dy$
\par $= z\pi ^{-1}(2\pi )^z \{ - \Gamma
(-z) \} \sin (z\pi /2) \zeta (1-z)$, \\
where for $\Gamma (z)$ the $(1,b)$-quasi-conformal extension of
Example 2.9.5.3 is used. Formula $(4)$ is initially valid for $-1<
\sigma <0$, but the right-hand side of $(4)$ is true also for each
$\sigma <0$, where $\sigma =Re (z)$. Thus this provides the
$(1,b)$-quasi-regular extension of $\zeta (z)$ on ${\cal
A}_b\setminus \{ 1 \} $ and the following formula is satisfied:
\par $(5)$ $\zeta (1-z) = 2^{1-z} \pi ^{-z} \cos (z\pi /2)
\Gamma (z)\zeta (z)$. \\
Equation $(5)$ transforms into
\par $(6)$ $\zeta (z)=\chi (z)\zeta (1-z),$ where \par
$\chi (z)=2^z\pi ^{z-1}\sin (\pi z/2)\Gamma (1-z)$\\
by changing $z$ into $1-z$. Then $\chi (z)= \pi ^{z-1/2} \Gamma (1/2
-z/2)/\Gamma (z/2)$, hence $\chi (z)\chi (1-z)=1$. Then $\xi (z)=
\xi (1-z)$ for each $Re (z)\ne 1/2$, where $\xi (z)=z(z-1)\pi
^{-z/2} \Gamma (z/2) \zeta (z)/2$, consequently,
\par $(7)$ $\Upsilon (z) = \Upsilon (-z)$ for each $Re (z)\ne 0$, \\
where $\Upsilon (z) = \xi (z +1/2)$. Since $(2^z)^* = 2^{z^*}$,
$(\pi ^{z-1})^* = \pi ^{z^*-1}$, $\sin (\pi z^*/2)=(\sin (\pi
z/2))^*$, $\Gamma (1-z^*)=(\Gamma (1-z))^*$ for each $z\in {\cal
A}_b$, then
\par $(8)$ $(\Upsilon (z))^*=\Upsilon (z^*)$ for each $z\in {\cal
A}_b$, \\ where $z^* := {\tilde z}$.

\par For $\sigma >0$ we have $\int_0^{\infty }x^{z-1}e^{-nx}dx=
n^{-z}\int_0^{\infty }y^{z-1}e^{-y}dy=n^{-z}\Gamma (z)$, since $n$
and $y$ are real and $\eta ^{z-1}$ is defined as $Exp ((z-1)Ln (\eta
))$ with the branch of the logarithm $Ln (R)$ real for $R>0$ so that
$n^{-z}$ and $y^{z-1}$ commute. For $\sigma >1$ we have the
convergent series $\sum_{n=1}^{\infty } \int_0^{\infty } x^{\sigma
-1} e^{-nx}dx= \Gamma (\sigma )\zeta (\sigma )$. Therefore, $\Gamma
(z)\zeta (z)=\sum_{n=1}^{\infty } \int_0^{\infty } x^{z-1}e^{-nx}dx$
\par $=\int_0^{\infty } x^{z-1}
\sum_{n=1}^{\infty } e^{-nx}dx=\int_0^{\infty
}x^{z-1}(e^x-1)^{-1}dx$.
\par Consider the integral
$J(z)=\int_C\eta ^{z-1}(e^{\eta }-1)^{-1}d\eta $, where the contour
$C$ starts at infinity on the positive real axis, encircles the
origin in the plane ${\bf R}\oplus M\bf R$ in the positive direction
besides the points $2\pi Mk$, where $0\ne k\in \bf Z$ and returns to
positive infinity. Therefore, $Arg (Ln (\eta ))$ varies from $0$ to
$2\pi M$ round the contour. So we take $C$ consisting of the real
axis from $\infty $ to $0<R<2\pi $, the circle $|z|=R$, and the real
axis from $R$ to $\infty $. Thus on the circle $|\eta ^{z-1}| = \exp
((\sigma -1) ln |\eta | - t~ arg (\eta ))\le |\eta |^{\sigma -1}
\exp (2\pi |t|)$ and $|Exp (\eta )-1|>A|\eta |$, where $arg (\eta )
= M^* Arg (\eta )$, $z=\sigma +tM$, $\sigma =Re (z)$, $t\in \bf R$.
Consequently, the integral round the circle tends to zero while
$R\to 0$ for $\sigma
>1$. Taking the limit with $R\to 0$ gives $J(z) = -\int_0^{\infty
}x^{z-1}(e^x - 1)^{-1}dx + \int_0^{\infty } (xExp (2\pi
M))^{z-1}(e^x-1)^{-1}dx=(exp (2\pi Mz) -1)\Gamma (z)\zeta (z) =2\pi
MExp (\pi zM)(\Gamma (1-z))^{-1}\zeta (z)$, hence
\par $(9)$ $\zeta (z)= \Gamma (1-z)Exp (-\pi Mz)(2\pi )^{-1}M^*
\int_C\eta ^{z-1} (Exp (\eta )-1)^{-1}d\eta $. \\
The latter formula has been proved for $\sigma >1$. But the integral
$J(z)$ is uniformly convergent in $G_M$ for any bounded region $G_M$
of the ${\bf R}\oplus M\bf R$ plane and uniformly by purely
imaginary $M\in {\cal A}_b$, $Re (M)=0$, $|M|=1$, where $G_M={\hat
R}_{M,{\bf i}}G_{\bf i}$. Thus Formula $(9)$ defines the
$(1,b)$-quasi-regular function on ${\cal A}_b\setminus \{ 1 \} $.
\par Formulas $(4,9)$ have been obtained by the same family $R_{z,x}$
of Example 2.2. If a pole of a complex meromorphic function is at
the real axis, then for its quasi-conformal extension with a marked
point $y_0=0$ its pole will remain the same real pole, since the
rotation axis is $\bf R$. Thus the only possible singularities of
$\zeta (z)$ may be poles of $\Gamma (1-z)$, $z=1,2,3,...$. In view
of $(4)$ $\zeta (z)$ is regular at $z=2,3,...$, more exactly $J(z)$
vanishes at these points (see \cite{titchm} and Theorem 2.11
\cite{luoyst,lusmfcdv}). At $z=1$ we have $J(1)= \int_C(Exp
(z)-1)^{-1}dz=2\pi M$ and $\Gamma (1-z)=-(z-1)^{-1} +...$, hence the
residue at this pole is $1$.

\par {\bf 2.} For the logarithmic derivative $\psi (1+z)=d Ln \Gamma (1+z)
/dz$ of the gamma function there is the expression $\psi
(1+z)=-C-\sum_{k=1}^{\infty } ((z+k)^{-1}-k^{-1})$ (see Formula
VII.89(9) in \cite{lavrshab}). Hence it is valid for its
$(1,b)$-quasi-meromorphic extension with the operators ${\hat
R}_{z,y}$ as in Example 2, where $y_0=0$, $2\le b\le 3$. Take
$-1<a<0$, then in view of the noncommutative ${\cal A}_b$ analog of
the Jordan Lemma 9 and Notes 10 above, 2.47 \cite{lusmtslt,luhcnlt}
with $-W:= \{ z: Re (z)<s_0 \} $ instead of $W$ and with $a<s_0<0$
and Theorem 3.9 about residues \cite{lusmfcdv,luoyst,luhcnfcdv} we
have
\par $(1)$ $\zeta (z)=\exp (M\pi z)(2\pi )^{-1}M^*\int_{a-M\infty
}^{a+M\infty } \{ \psi (1+\eta ) - Ln (\eta ) \} \eta ^{-z}d\eta $
for each $\sigma >1$, where $M\in {\cal A}_b$, $Re (M)=0$, $|M|=1$,
$z\in {\cal A}_b$, $\sigma =Re (z)$, $z=\sigma +M v$, $\sigma , v\in
\bf R$,  $ -1< a <0$. \par The function $\{ \psi (1+\eta ) - Ln
(\eta ) \} \eta ^{-z}$ is $O(|\eta |^{-1-\sigma })$, consequently,
the integral in $(1)$ is convergent and Formula $(1)$ is valid by
analytic continuation for $\sigma >0$. Again using the
noncommutative analog of the Jordan lemma transform the integral in
$(1)$ to
\par $(2)$ $\zeta (z)= - \sin (\pi z) \pi ^{-1}
\int_0^{\infty } \{ \psi (1+x) - ln (x) \} x^{-z}dx$ for each
$0<\sigma <1$.
\par The function $\{ \psi (1+\eta ) - Ln (\eta ) \} \eta ^{-z}$ is
real on $(0,\infty )= \{ x\in {\bf R}: 0<x \} $, where the branch of
$Ln$ is such that $Ln|_{\bf R}=ln: (0,\infty )\to \bf R$. In view of
the theorems about uniqueness and inversion of the noncommutative
version of the Mellin transformation the $(1,b)$-quasi-regular
extension of $\zeta (z)$ coincides with the noncommutative version
of the Mellin transform $(2)$, when $z\in {\cal A}_b$ with $0<Re
(z)<1$. Then ${\hat R}_{z,y}g(y)=g(z)$ for each $y\in \bf C$ and
$z\in {\cal A}_b$ with $0<Re (y)=Re (z)<1$ such that ${\hat
R}_{z,y}y=z$, where $g(z) := \int_0^{\infty } \{ \psi (1+x) - ln (x)
\} x^{-z}dx$, $y_0=0$. In view of Theorem 2.18 there exists
$(\partial \int_0^{\infty } \{ \psi (1+x) - ln (x) \}
x^{-z}dx/\partial z).h =(\int_0^{\infty } \{ \psi (1+x) - ln (x) \}
ln (x) x^{-z}dx).h$ for each $h\in {\bf R}\oplus M\bf R$ and every
$0<\sigma <1$, where $z=\sigma + Mv$. Thus $g(z)$ satisfies
$(Q1,Q6)$ and $g'(z)$ satisfies $(Q7)$ when $g'(z)\ne 0$,
consequently, $g(z)$ is the $(1,b)$-quasi-regular function.\par In
view of 1$(6)$ there is the symmetry relation: $g(z)=-(\sin (\pi
z))^{-1}\pi \zeta (z)=- (\sin (\pi (1-z))^{-1} \pi \chi (z)\zeta
(1-z)$, since $\sin (\pi -\phi )=\sin (\phi )$ for each $\phi \in
{\cal A}_b$, where \\ $\chi (z):=2^z\pi ^{z-1} \sin (\pi z/2) \Gamma
(1-z)$, $\chi (z)\chi (1-z)=1$. But $|2^z|=2^{\sigma }$, $|\pi
^{z-1}|=\pi ^{\sigma -1}$, $\sin (\pi z/2)=0$ if and only if $z=2k$
with $k\in \bf Z$, $\sin (\pi z/2)$ has not poles, $\Gamma (1-z)$
has not zeros, $\Gamma (1-z)$ has a pole at $z$ if and only if
$z=1,2,3,...$, consequently, $\chi (z)$ has not any zero or pole in
the domain $V:=\{ z\in {\cal A}_b: 0<Re (z)<1 \} $. At the same time
the multiplier $(\sin (\pi z))^{-1}\pi |_V$ has not any pole or zero
in $V$.
\par {\bf 3.} Consider now new type of an extension in spherical
${\cal A}_b$-coordinates. Let \par $(1)$ $\psi (x) :=
\sum_{n=1}^{\infty } \exp (-n^2\pi x)$,\\ where $x>0$, then
\par $(2)$ $\zeta (y)=\pi ^{y/2} [\Gamma (y/2)]^{-1} \int_0^{\infty }
x^{y/2-1}\psi (x)dx$ \\ for $\sigma =Re (y)>1$, $y\in \bf C$. It is
known that \par $(3)$ $2\psi (x)+1= [2\psi (1/x) +1]/(x)^{1/2}$ for
each $x>0$. Therefore, from $(2,3)$ it follows, that
\par $(4)$ $\pi ^{-y/2} \Gamma (y/2) \zeta (y) = \int_0^1
x^{y/2-1} \psi (x) dx +\int_1^{\infty } x^{y/2-1}\psi (x)dx$ \\ $=
\int_0^1x^{y/2-1}[\psi (1/x) (x)^{-1/2} +(x)^{-1/2}/2 -1/2]dx +
\int_1^{\infty } x^{y/2-1}\psi (x)dx$ \\  $= 1/(y-1)-1/y
+\int_0^1x^{y/2-3/2}\psi (1/x)dx+ \int_1^{\infty }x^{y/2-1}\psi
(x)dx$ \\  $= 1/[y(y-1)] + \int_1^{\infty } (x^{-y/2-1/2} +
x^{y/2-1}) \psi (x)dx$. \\
The last integral is convergent for all values of $y\in \bf C$, so
Formula $(4)$ holds for all values of $y$ by analytic continuation
(see Formulas 2.6.1-4 in \S 2.6 \cite{titchm}). Write the term
$1/[y(y-1)]$ in the form:
\par $(5)$ $w(q) := 1/(y-1) -1/y= - [\int_0^{\infty }[\exp (-ty) +
\exp (-t(1-y))]dt$ \\  $= -[\int_0^{\infty } \exp (-t/2) [\exp (-tq)
+ \exp (tq)]dt= - [\int_{-\infty }^{\infty }
\exp (-|t|/2) \exp (-tq)dt$, \\
which converges in the strip $-1/2<Re (q)<1/2$, where $q=y-1/2$.
Then the term $\int_1^{\infty } (x^{-y/2-1/2} + x^{y/2-1}) \psi
(x)dx$ putting $q=y-1/2$ and then $x=e^t$ write in the form:
\par $(6)$ $\int_1^{\infty } (x^{-y/2-1/2} + x^{y/2-1}) \psi (x)dx$ \\
$= \int_1^{\infty }(x^{-3/4-q/2} + x^{-3/4 +q/2})\psi (x)dx $
$=\int_0^{\infty } \exp (-3t/4) [\exp (-tq/2) + \exp (tq/2)] \psi
(e^t) e^tdt $ \\  $= \int_0^{\infty } \exp (t/4) \psi (e^t)\exp
(-tq/2)dt+\int_{-\infty }^0\exp (-t/4) \psi (e^{-t})\exp (-tq/2)dt$
\\  $=\int_{-\infty }^{\infty }\exp (|t|/4) \psi (\exp (|t|))
\exp (-tq/2)dt$. \\
Therefore, Formulas $(4-6)$ give:
\par $(7)$ $\pi ^{-q/2-1/4} \Gamma (q/2+1/4) \zeta (q+1/2) =$ \\
$ \int_{-\infty }^{\infty }[ - \exp (-|t|/2)+ 2 \exp (|t|/2) \psi
(\exp (2|t|))] \exp (-tq)dt$ \\
valid on $\bf C$ by analytic continuation. Then \par $(8)$ $\xi (y)=
y(y-1) [\pi ^{-y/2} \Gamma (y/2) \zeta (y)]/2 = [w(y-1/2)]^{-1}[\pi
^{-y/2} \Gamma (y/2) \zeta (y)]/2$ \\ is the integral function on
$\bf C$.
\par Take the family ${\hat R}_{E(z),E(x)}$ satisfying Condition
2.9.5.2$(2)$. In view of Theorems 14 and 24 $w(y)$ has the
$(1,b)$-quasi-regular extension $w^s(p)$ in spherical ${\cal
A}_b$-coordinates, where $2\le b\le 3$, $w^s(p) := - [\int_{-\infty
}^{\infty } \exp (-|t|/2) \exp (-E(tp))dt$. In accordance with
Proposition 25.1 the function $1/w^s(p)$ is $(1,b)$-quasi-regular in
spherical ${\cal A}_b$-coordinates in the domain $ -1/2<Re (p)<1/2$.
By Corollary 9.2 the product of $(1,b)$-quasi-regular functions with
the same family ${\hat R}_{z,x}$ is $(1,b)$-quasi-regular. Then from
Definition 2.2.1 it follows, that the product $f^s_1f^s_2$ of
$(1,b)$-quasi-regular functions $f^s_1$ and $f^s_2$ in spherical
${\cal A}_b$-coordinates with the same family ${\hat R}_{E(z),E(x)}$
is $(1,b)$-quasi-regular in spherical ${\cal A}_b$-coordinates,
since $f_1=f^s_1\circ E^{-1}$ and $f_2=f^s_2\circ E^{-1}$ are
$(1,b)$-quasi-regular.
\par  On the other hand, the right side of Equation $(8)$
gives the $(1,b)$-quasi-meromorphic in spherical ${\cal
A}_b$-coordinates extension \par $\int_{-\infty }^{\infty }[ - \exp
(-|t|/2)+ 2 \exp (|t|/2) \psi (\exp (2|t|))] \exp (-E(tp))dt =:
g^s(p)$ in accordance with Theorems 14, 24. Put $\Omega (q) := \xi
(q+1/2)$. Then $\Omega (q)$ has the $(1,b)$-quasi-integral in
spherical ${\cal A}_b$-coordinates extension $\Omega ^s(p)=
[w^s(p)]^{-1} g^s(p)$. The function $f(t) := - \exp (-|t|/2)+ 2 \exp
(|t|/2) \psi (\exp (2|t|))$ is real-valued and even on $\bf R$. In
view of Theorems 24 and 25 $\Omega ^s(p)$ has the symmetry
properties 25$(1',2)$.  This also can be seen from Equations
$(7,8)$. The symmetry property 25$(1')$ for $f^s$ implies 25$(1)$
for $f=f^s\circ E^{-1}$, since if $z =E(p)$, then the adjoint number
is ${\tilde z} = E(p_0 -p_1i_1 +p_2i_2+...+p_{2^b-1}i_{2^b-1})$ in
accordance with Formulas 2.2.1$(1,2)$.
\par It is known that $\zeta (z)$ has no any poles in $\bf
C$ besides $z=1$, that is, $\zeta (z)$ contains only complex zeros
in the domain $0<Re (y)<1$ in $\bf C$. It is well-known that all
complex zeros of $\zeta (z)$ are in the complex strip $0<\sigma <1$
and they form a discrete set in $\bf C$ without finite accumulation
points \cite{titchm}. Thus the function $f=f^s\circ E^{-1}$ with
$f^s(p)=\Omega ^s(p)$ satisfies conditions of Theorem 2.17, since
$\xi (z)$ has not any real zeros and all its complex zeros are in
the strip $0\le Re (z)\le 1$ (see page 30 \cite{titchm}). On the
other hand, $\zeta (z)$ and $\xi (z)$ have common all complex zeros
and $E(y)|_{\bf C} =y$ for each $y\in \bf C$. Thus in view of
Theorem 2.17 it is proved the following.

\par {\bf 27. Theorem.} {\it All complex zeros of the $\zeta $
function lie only on the line $Re (z)=1/2$.}

\par {\bf 27.1. Remark.} This is not so surprising, since by
Theorem 2.13 \cite{titchm} each meromorphic function
$f(s)=G(s)/P(s)$, where $G$ is an integral function of finite order
and $P$ is a polynomial on $\bf C$, and $f$ is satisfying the
symmetry property 26.1(6) and having the series expansion
$f(s)=\sum_{n=1}^{\infty } a_nn^{-s}$ absolutely convergent for
$\sigma >1$, then $f(s)$ is $c\zeta (s)$, where $c=const$, $a_n\in
\bf C$ is a constant for each $n\in \bf N$. On the other hand, the
class of $(r,b)$-quasi-conformal functions is more narrow and
specific in comparison with the class of ${\cal A}_b$ holomorphic
functions, where $1\le r<b\le 3$ (see also Notes 2.13 and 2.17.1).
Moreover, the class of $(1,b)$-quasi-integral functions is more
narrow than that of $(1,b)$-quasi-regular which in its turn is
restricted by Theorem 24. Mention, that if $z =z_0 +z_1i_1 +z_2i_2
+z_3i_3 =E_2(p)$, then $E_2( -p) = -z_0 -z_1i_1 +z_2i_2 -z_3i_3$ in
accordance with Formulas 2.2.1$(1)$ and 16.1$(3)$, where $z_0, z_1,
z_2, z_3\in \bf R$, $z, p\in \bf H$. Consequently, $E_2$ and $E_6$
are neither even nor odd functions. More narrow class is that of
satisfying symmetry properties 25$(1',2)$, which need to be met for
using Theorem 2.17. For example, the Dirichlet function does not
satisfy conditions of Theorem 2.17 (see \S 10.25 \cite{titchm}).
\par Consider the identity $\int_{ -\infty }^{\infty }
f(t)\exp (p_0+{\bf i} p_1t) dt= g(p)\int_{-\infty }^{\infty }
f(t)\exp ((1-p_0)t - {\bf i}p_1t)dt$ for an original nonzero
function $f: {\bf R}\to \bf R$ and a meromorphic function $g$ in
$\bf C$ such that $g$ is holomorphic and without zeros in the band
$G := \{ p\in {\bf C}: s_0<p_0<s_1, s_0<1-p_0<s_1 \} $, where $g$
may have only isolated poles in $\bf C$, $p_0, p_1\in \bf R$,
$p=p_0+{\bf i}p_1$, $0<s_0(f)<s_1(f)<1$, $|f(t)|<C_1\exp (-s_1t)$
for each $t<0$, $|f(t)|<C_2\exp (s_0t)$ for each $t\ge 0$,
$s_0=s_0(f)$, $s_1=s_1(f)$. Then $g(p)g(1-p)=1$ and ${\bar
g}(p)=g({\bar p})$ in $G$ besides poles and $g(1/2)=1$, and ${\bar
F}(p)=F({\bar p})$ in $G$, where $F(p)= \int_{ -\infty }^{\infty }
f(t)\exp (pt) dt$. In particular, $g'(1/2+{\bf i}p_1)=g'(1/2-{\bf
i}p_1)$ for each $p_1\bf R$ besides poles of $g$. Since the
two-sided Laplace transformation of $f(t)$ is holomorphic in the
band $\{ p\in {\bf C}: s_0<p_0<s_1 \} $, then the differentiation of
this identity by $p$ in $G$ gives:
\par $\int_{-\infty }^{\infty }f(t)t\exp (pt)dt=
g'(p)(\int_{-\infty }^{\infty } f(t)\exp ((1-p)t)dt -
g(p)\int_{-\infty }^{\infty } f(t)t\exp ((1-p)t)dt$ (see also
Theorem 2.18). Therefore, the class of such functions $F(p)$ is
narrow.

\par {\bf Acknowledgement.} The author is sincerely grateful
to Prof. Fred van Oystaeyen and Prof. Jan van Casteren for helpful
discussions on noncommutative analysis and geometry over quaternions
and octonions at Mathematical Department of Antwerpen University in
2002 and 2004 and for hospitality.


\begin{thebibliography}{99}
\bibitem{baez} J.C. Baez. "The octonions". Bull. Amer.
Mathem. Soc. {\bf 39: 2} (2002), 145-205.
\bibitem{berez} F.A. Berezin. "Introduction to
superanalysis" (D. Reidel Publish. Comp., Kluwer group: Dordrecht,
1987).
\bibitem{bourbal} N. Bourbaki. "Algebra. Algebraic structures.
Linear and polylynear algebra" (Moscow: Fizmatgiz, 1962).
\bibitem{brdeso} F. Brackx, R. Delanghe, F. Sommen.
"Clifford analysis" (London: Pitman, 1982).
\bibitem{connes} A. Connes. "Noncommutative geometry"
(Academic Press: San Diego, 1994).
\bibitem{dewitt} B. DeWitt. "Supermanifolds"
2d ed. (Cambridge Univ. Press: Cambridge, 1992).
\bibitem{emch} G. Emch. "M$\grave e$chanique quantique quaternionienne et
Relativit$\grave e$ restreinte". Helv. Phys. Acta {\bf 36}, 739-788
(1963).
\bibitem{eng} R. Engelking. "General topology"
(Heldermann: Berlin, 1989).
\bibitem{guetze} F. G\"ursey, C.-H. Tze.
"On the role of division, Jordan and related algebras in particle
physics" (World Scientific Publ. Co.: Singapore, 1996).
\bibitem{grauert} H. Grauert, I. Lieb, W. Fischer. "Differential-
und Integralrechnung" (Springer-Verlag: Berlin, 1967-1968).
\bibitem{hamilt} W.R. Hamilton. "Selected papers. Optics. Dynamics.
Quaternions" (Nauka: Moscow, 1994).
\bibitem{kansol} I.L. Kantor, A.S. Solodovnikov.
"Hypercomplex numbers" (Berlin: Springer, 1989).
\bibitem{kamyn} L.I. Kamynin. "Course of Mathematical Analysis"
(Moscow State Univ. Press: Moscow, 1995).
\bibitem{khren} A. Khrennikov. "Superanalysis",
(Series "Mathem. and its Applic."; V. {\bf 470}; Kluwer: Dordrecht,
1999).
\bibitem{lavrshab} M.A. Lavrentjev, B.V. Shabat.
"Methods of theory of functions of the complex variable"
(Moscow: Nauka, 1987).
\bibitem{lawmich} H.B. Lawson, M.-L. Michelson. "Spin geometry"
(Princeton: Princ. Univ. Press, 1989).
\bibitem{luoyst} S.V. L\"udkovsky, F. van Oystaeyen.
"Differentiable functions of quaternion variables". Bull. Sci. Math.
(Paris). Ser. 2. {\bf 127} (2003), 755-796.
\bibitem{lusmfcdv} S.V. Ludkovsky. "Differentiable functions of
Cayley-Dickson numbers and line integration". J. Math Sci. {\bf 141:
3} (2007), 1231-1298 (previous version: Los Alam. Nat. Lab. {\bf
math.NT/0406048}; {\bf math.CV/0406306}; {\bf math.CV/0405471}).
\bibitem{luhcnfcdv} S.V. Ludkovsky. "Differentiable
functions of Cayley-Dickson numbers". Hypercomplex numbers in Geometry
and Physics. {\bf 3: 1} (2005), 93-140.
\bibitem{lufscdvm} S.V. Ludkovsky.
"Functions of several Cayley-Dickson variables and manifolds over
them". J. Mathem. Sci. {\bf 141: 3} (2007), 1299-1330 (previous
variant: Los Alamos Nat. Lab. {\bf math.CV/0302011}).
\bibitem{ludanavf} S.V. Ludkovsky. "Algebras of vector fields over the
quaternion skew field". Dokl. Akad. Nauk. {\bf 403: 3} (2005), 309-312.
\bibitem{lusmalop} S.V. Ludkovsky.  "Algebras of operators in Banach
spaces over the quaternion skew field and the octonion algebra".
Sovrem. Mathem. i ee Pril. {\bf 35} (2005) (previous variant: Los Alam. Nat. Lab.
math.OA/0603025).
\bibitem{ludagpd} S.V. Ludkovsky. "Groups of pseudoconformal
diffeomorphisms of quaternion variables". Dokl. Akad. Nauk. {\bf 408: 5}
(2006), 587-590 (Doklady Mathematics {\bf 73: 3} (2006), 403-406).
\bibitem{lusmnfgpcd} S.V. Ludkovsky.
Sovrem. Mathem. Fundam. Napravl. "Normal families of functions and
groups of pseudoconformal diffeomorphisms of quaternion and octonion
variables", {\bf 18} (2006), 101-164 (previous variant: Los Alam.
Nat. Lab. math.DG/0603006).
\bibitem{lusmgdlcm} S.V. Ludkovsky.  "Stochastic processes on geometric
loop groups and diffeomorphism groups of connected manifolds,
associated unitary representations". J. Math. Sci. {\bf 141: 3}
(2007), 1331-1384 (previous version: Los Alam. Nat. Lab.
math.AG/0407439, July 2004).
\bibitem{lusmtslt} S.V. Ludkovsky.  "The two-sided Laplace
transformation over the Cayley-Dickson algebras and its
applications". Sovrem. Mathem. Fundam. Napravl. {\bf 122} (2006)
(see also: {\bf math.CV/0612755}).
\bibitem{luhcnlt} S.V. Ludkovsky. "The Laplace transformation
over the Cayley-Dickson algebras". Hypercomplex numbers in Geometry
and Physics. {\bf 5: 1} (2006), 67-99.
\bibitem{oystaey} F. van Oystaeyen. "Algebraic geometry
for associative algebras" (Series "Lect. Notes in Pure and Appl.
Mathem."; V. {\bf 232}; Marcel Dekker: New York, 2000).
\bibitem{razmus} Y.P. Razmyslov. "Identities of algebras
and their representations" (Nauka: Moscow, 1989).
\bibitem{rothe} H. Rothe. "Systeme Geometrischer Analyse"
in: "Encyklop\"adie der Mathematischen Wissenschaften. Band 3.
Geometrie", 1277-1423 (Leipzig: Teubner, 1914-1931).
\bibitem{shabat} B.V. Shabat. "Introduction into complex analysis"
(Moscow: Nauka, 1985).
\bibitem{span} E.H. Spanier. "Algebraic topology"
(Acad. Press: New York, 1966).
\bibitem{titchm} E.C. Titchmarsh. "The theory of the Riemann
zeta-function" (Oxford: Clarendon Press, 1988).
\bibitem{zorich} V.A. Zorich. "Mathematical Analysis"
(Nauka: Moscow, 1984).

\end{thebibliography}
\end{document}